\ifx\shlhetal\undefinedcontrolsequence\let\shlhetal\relax\fi
\documentstyle[11pt]{article}

\input amssym.def
\input amssym.tex

\newcount\skewfactor
\def\mathunderaccent#1#2 {\let\theaccent#1\skewfactor#2
\mathpalette\putaccentunder}
\def\putaccentunder#1#2{\oalign{$#1#2$\crcr\hidewidth
\vbox to.2ex{\hbox{$#1\skew\skewfactor\theaccent{}$}\vss}\hidewidth}}
\def\name{\mathunderaccent\tilde-3 }

\newcommand{\forces}{\Vdash}

\newcommand{\V}{{\bf V}} 
\newcommand{\lesdot}{\mathrel{\mathord{<}\!\!\raise 
0.8 pt\hbox{$\scriptstyle\circ$}}} 
\newcommand{\comp}{\circ}   
\newcommand{\conc}{{}^\frown\!}
\newcommand{\lh}{\ell g\/} 
\newcommand{\rest}{{\restriction}}
\newcommand{\dom}{{\rm dom}} 

\newcommand{\cf}{{\rm cf}}
\newcommand{\bP}{{\Bbb P}}
\newcommand{\bQ}{{\Bbb Q}}
\newcommand{\bR}{{\Bbb R}}
\newcommand{\nbQ}{{\name{\Bbb Q}}}
\newcommand{\cB}{{\cal B}}
\newcommand{\QED}{\hfill\vrule width 6pt height 6pt depth 0pt 
\vspace{0.1in}} 
\newcommand{\Proof}{\noindent{\sc Proof} \hspace{0.2in}} 

\newcommand{\Ker}{{\rm Ker}}

\newcommand{\Hchi}{({\cal H}(\chi),\in,<^*_\chi)}
\newcommand{\mod}{{\rm mod}}
\newcommand{\UP}{{\bf UP}}
\newcommand{\Pro}{{\rm Pr}}
\newcommand{\otp}{{\rm otp}}
\newcommand{\Gra}{{\cal G}^\theta_S}

\newcommand{\innagra}{{\cal G}^*_{\bar{N}}}
\newcommand{\tagra}{{\cal G}^\heartsuit_{\bar{N},D,X,\bar{a}}}
\newcommand{\grapik}{{\cal G}^\spadesuit_{\bar{M},\langle\bar{N}^i:i<\delta
\rangle}}
\newcommand{\com}{{\rm COM}}
\newcommand{\inc}{{\rm INC}}
\newcommand{\family}{{\cal S}}
\newcommand{\dkl}{{\frak D}_{<\kappa,<\lambda}(\mu^*)}

\newcommand{\dak}{{\frak D}^\alpha_{<\kappa,<\lambda}(\mu^*)}
\newcommand{\dlk}{{\frak D}^\lambda_{<\kappa,<\lambda}(\mu^*)}
\newcommand{\ck}{{\frak C}_{<\kappa}(\mu^*)}
\newcommand{\ckm}{{\frak C}_{<\kappa}^-(\mu^*)}
\newcommand{\ckp}{{\frak C}_{<\kappa}^\spadesuit(\mu^*)}
\newcommand{\baza}{{\hat{\cal E}}}
\newcommand{\nacc}{{\rm nacc}}
\newcommand{\acc}{{\rm acc}}
\newcommand{\Tr}{{\rm Tr}}
\newcommand{\rk}{{\rm rk}}
\newcommand{\tree}{{\cal T}}
\newcommand{\ftr}{{\rm FTr}}
\newcommand{\ftrw}{{\rm FTr}_{wk}}
\newcommand{\proj}{{\rm proj}}
\newcommand{\inver}{\lim\limits^{\leftarrow}}

\newcommand{\Axkt}{{\rm Ax}^\kappa_\theta}
\newcommand{\HOM}{{\rm HOM}}
\newcommand{\cohen}{{\rm Cohen}}
\newcommand{\Lim}{{\rm Lim}}
\newcommand{\PART}[2]{

	\bigskip\bigskip\bigskip\bigskip\bigskip
	\begin{center} \LARGE  #2: #1 \end{center}

	\bigskip\medskip

	\renewcommand{\thesection}{#1.\arabic{section}}
}
\newtheorem{theorem}{Theorem}[section] 
\newtheorem{claim}{Claim}[theorem]
\newtheorem{fact}[theorem]{Fact}
 
\newtheorem{proposition}[theorem]{Proposition}

\newtheorem{pwa}[theorem]{Problem we address} 
\newtheorem{definition}[theorem]{Definition}
\newtheorem{notation}[theorem]{Notation}

\newtheorem{conclusion}[theorem]{Conclusion}
\newtheorem{remark}[theorem]{Remark}

\newtheorem{hypothesis}[theorem]{Hypothesis}
\newtheorem{ourassum}{Our Assumptions} 
\newtheorem{Gquest}[theorem]{G\"obel's question}
\title{Not collapsing cardinals $\leq\kappa$ in $(<\kappa)$--support
iterations: Part I}  
\author{{\bf Saharon Shelah}\thanks{\ \ Research supported by German-Israeli
Foundation for Scientific Research\ \&\ Development Grant No. G-294.081.06/93
and by The National Science Foundation Grant No. 144--EF67. Publication No
587.}\\ 
Institute of Mathematics\\
The Hebrew University of Jerusalem\\
91904 Jerusalem, Israel\\
and\\
Department of Mathematics\\
Rutgers University\\
New Brunswick, NJ 08854, USA\\
and\\
Mathematics Department\\
University of Wisconsin -- Madison\\
Madison, WI 53706, USA
}
\date{\today}

\begin{document} 
\maketitle 

\begin{abstract}
We deal with the problem of preserving various versions of completeness in
$(<\kappa)$--support iterations of forcing notions, generalizing the case
``$S$--complete proper is preserved by CS iterations for a stationary
co-stationary $S\subseteq\omega_1$''. We give applications to Uniformization
and the Whitehead problem. In particular, for a strongly inaccessible cardinal
$\kappa$ and a stationary set $S\subseteq\kappa$ with fat complement we can
have uniformization for $\langle A_\delta:\delta\in S'\rangle$, $A_\delta
\subseteq\delta=\sup A_\delta$, $\cf(\delta)=\otp(A_\delta)$ and a stationary
non-reflecting set $S'\subseteq S$. 
\end{abstract}
\vfill
\eject

\section*{Annotated Content}

\noindent {\bf Section 0:\quad Introduction}\qquad We put this work in a
context and state our aim.\\
--{\em 0.1\ Background: Abelian groups}\\
--{\em 0.2\ Background: forcing}\quad [We define $(<\kappa)$--support
iteration.]\\
--{\em 0.3\ Notation}
\medskip

{\sc CASE A}

Here we deal with Case A, say $\kappa=\lambda^+$, $\cf(\lambda)=\lambda$,
$\lambda=\lambda^{<\lambda}$. 
\medskip

\noindent {\bf Section A.1:\quad Complete forcing notions}\qquad We define
various variants of completeness and related games; the most important are 
{\em the strong $\family$--completeness} and {\em real $(\family_0,\hat{
\family}_1,D)$--completeness}. We prove that the strong
$\family$--completeness is preserved in $(<\kappa)$--support iterations
(\ref{iterstrong}) 
\smallskip

\noindent {\bf Section A.2:\quad Examples}\qquad We look at guessing clubs
$\bar{C}=\langle C_\delta:\delta\in S\rangle$. If $[\alpha\in\nacc(C_\delta)\
\Rightarrow\ \cf(\alpha)<\lambda]$ we give a forcing notion (in our context)
which adds a club $C$ of $\kappa$ such that $C\cap\nacc(C_\delta)$ is bounded
in $\delta$ for all $\delta\in S$. (Later, using a preservation theorem, we
will get the consistency of ``no such $\bar{C}$ guesses clubs''.) Then we deal
with uniformization (i.e.\ $\Pro_S$) and the (closely related) being Whithead.
\smallskip

\noindent {\bf Section A.3:\quad The iteration theorem}\qquad We deal
extensively with (standard) trees of conditions, their projections and inverse
limits. The aim is to build a $(\bP_\gamma,N)$--generic condition forcing
$\name{G}_\gamma\cap N$, and the trees of conditions are approximations to
it. The main result in the preservation theorem for our case
(\ref{iterweak}). 
\smallskip

\noindent {\bf Section A.4:\quad The Axiom}\qquad We formulate a Forcing Axiom
relevant for our case and we state its consistency.
\medskip

{\sc CASE B}

Here we deal with $\kappa$ strongly inaccessible, $S\subseteq\kappa$ usually a
stationary ``thin'' set of singular cardinals. There is no point to ask even
$\aleph_1$--completeness , so the completeness demands are only on sequences
of models. 
\medskip

\noindent {\bf Section B.5:\quad More on complete forcing notions}\qquad We
define completeness of forcing notions with respect to a suitable family
$\baza$ of increasing sequences $\bar{N}$ of models, say, such that
$\bigcup\limits_{j<\delta}N_j\cap\kappa\notin S$ for limit $\delta\leq
\lh(\bar{N})$. $S$ is the non-reflecting stationary set where ``something is
done''. The suitable preservation theorem for $(<\kappa)$--support iterations
is proved in \ref{nextiter}. So this $\baza$ plays a role of $\family_0$ of
Case A, and the presrvation will play the role of preservation of strong
$\family_0$--completeness. We end defining the version of completeness (which
later we prove is preserved; it is parallel to $(\family_0,\hat{\family}_1,
D)$--completeness of Case A).
\smallskip

\noindent {\bf Section B.6:\quad Examples for innaccessible}\qquad We present
forcing notion taking care of $\Pro_S$, at laest for cases which are locally
OK, say, $S\subseteq\kappa$ is stationary non-reflecting. WE show that it
satisfies the right properties (for iterating) for the naturally defined
$\baza_0,\baza_1$. Then we turn to the related problem of Whitehead group. 
\smallskip

\noindent {\bf Section B.7:\quad The iteration theorem for inaccessible
$\kappa$}\qquad We show that completeness for $(\baza_0,\baza_1)$is preserved
in $(<\kappa)$--support iterations (this covers the uniformization). Then we
prove the $\kappa^+$--cc for the simplest cases.
\smallskip

\noindent {\bf Section B.8:\quad The Axiom and its applications}\qquad We
phrase the axiom and prove its consistency. The main case is for a stationary
set $S\subseteq\kappa$ whose complement is fat, but checking that forcing
notions fit is clear for forcing notion related to non-reflecting subsets
$S'\subseteq S$, So $S$ can have stationary intersection with
$S^\kappa_\sigma$ for any regular $\sigma<\kappa$. The instance of
$S\cap\mbox{inaccessible}$ is not in our mind, but it is easier -- similar to
the successor case. Next we show the consistency of ``GCH + there are almost
free Abelian groups in $\kappa$, and all of them are Whitehead''. We start
with an enough indestructible weakly compact cardinal and a stationary
non-reflecting set $S\subseteq\kappa$, for simplicity $S\subseteq
S^\kappa_{\aleph_0}$, and then we force the axiom. Enough weak compactness
remains, so that we have:\quad every stationary set $S'\subseteq
\kappa\setminus S$ reflects in inaccessibles, hence ``$G$ almost free in
$\kappa$'' implies $\Gamma(G)\subseteq S\mod D_\lambda$'', but the axiom makes
all of them Whitehead.
\eject

\setcounter{section}{-1}

\section{Introduction} In the present paper we deal with the following
question of Theory of Forcing: 
\begin{pwa}
\label{pwa}
Iterate with $(<\kappa)$--support forcing notions not collapsing cardinals
$\leq\kappa$ preserving this property, generalizing ``$S$--complete proper is
preserved by CS iterations for a stationary co-stationary $S\subseteq
\omega_1$''.\\
But we concentrate on the ZFC case (or deal with cardinals which may exists in
${\bf L}$) and we demand that no bounded subsets of $\kappa$ are added. 
\end{pwa}
We use as our test problems instances of uniformization (see
\ref{uniformization} below) and Whitehead groups (see \ref{Wgroup} below),
but the need for \ref{pwa} comes from various questions of Set Theory. The
case of CS iteration and $\kappa=\aleph_1$ has gotten special attention (so we
generalize {\em no new real\/} case by $S$--completeness, see \cite[Ch
V]{Sh:f}) and is a very well understood case, but still with consequences in
CS iterations of $S$--complete forcing notions. This will be our starting
point.

One of the questions which caused us to look again in this direction was:
\begin{quotation}
\noindent is it consistent with ZFC + GCH that for some regular $\kappa$ there
is an almost free Abelian group of cardinality $\kappa$, but every such
Abelian group is a Whitehead one?
\end{quotation}
By G\"obel Shelah \cite{GbSh:579}, we have strong counterexamples for
$\kappa=\aleph_n$: an almost free Abelian group $G$ on $\kappa$ with $\HOM(G,
{\Bbb Z})=\{0\}$. The here idea is that we have an axiom for $G$ with
$\Gamma(G)\subseteq S$ (to ensure being Whitehead) and some reflection
principle gives  
\[\Gamma(G)\setminus S\mbox{ is stationary }\quad\Rightarrow\quad G\mbox{ is
not almost free in }\kappa,\]
(see \ref{aksitsapp}). This stream of investigations has a long history
already, one of the starting points was \cite{Sh:186} (see earlier references
there too), and later \cite{MkSh:274}, \cite{MkSh:313}.  

\begin{definition}
\label{uniformization}
Let $\kappa\geq\lambda$ be cardinals.
\begin{enumerate}
\item We let $S^\kappa_\lambda\stackrel{\rm def}{=}\{\delta<\kappa:
\cf(\delta)=\cf(\lambda)\}$.
\item {\em A $(\kappa,\lambda)$--ladder system} is a sequence $\bar{\cal A}=
\langle A_\delta: \delta\in S\rangle$ such that the set $\dom(\bar{\cal A})=S$
is a stationary subset of $S^\kappa_\lambda$ and 
\[(\forall\delta\in S)(A_\delta\subseteq\delta=\sup(A_\delta)\quad \&\quad
\otp(A_\delta)=\cf(\lambda)).\]
When we say that $\bar{\cal A}$ is a $(\kappa,\lambda)$--ladder system {\em is
on $S$} then we mean that $\dom(\bar{\cal A})=S$. 
\item  Let  $\bar{\cal A}$ be a $(\kappa,\lambda)$--ladder system. We say that
$\bar{\cal A}$ has {\em the $h^*$--{\bf U}niformization {\bf P}roperty} (and
then we may say that it has $h^*$--$\UP$) if $h^*:\kappa\longrightarrow\kappa$
and 
\begin{quotation}
\noindent {\em for every} sequence $\bar{h}=\langle h_\delta: \delta\in
S\rangle$, $S=\dom(\bar{\cal A})$, such that 
\[(\forall \delta\in S)(h_\delta:A_\delta\longrightarrow\kappa\ \ \&\ \
(\forall\alpha\in A_\delta)(h_\delta(\alpha)<h^*(\alpha))\]
{\em there is} a function $h:\kappa\longrightarrow\lambda$ with 
\[(\forall\delta\in S)(\sup\{\alpha\in A_\delta: h_\delta(\alpha)\neq
h(\alpha)\}<\delta\}).\]
If $h^*$ is constantly $\mu$ then we may write $\mu$--$\UP$; if $\mu=\lambda$
then we may omit it.
\end{quotation}
\item For a stationary set $S\subseteq S^\kappa_\lambda$, let $\Pro_{S,\mu}$
be the following statement
\begin{description}
\item[$\Pro_{S,\mu}\equiv$] each $(\kappa,\lambda)$--ladder system $\bar{\cal
A}$ on $S$ has the $\mu$--Uniformization Property. 
\end{description}
We may replace $\mu$ by $h^*$; if $\mu=\lambda$ we may omit it.
\end{enumerate}
\end{definition}
There are several works on the $\UP$, for example the author proved that it is
consistent with GCH that {\em there is} a $(\lambda^+,\lambda)$--ladder system
on $S^{\lambda^+}_{\lambda}$ with the Uniformization Property (see
\cite{SnKi:C3}, for more general cases see \cite{Sh:186}), but necessarily not
every such system has it (see \cite[AP, \S3]{Sh:f}). In the present paper we
are interested in a stronger statement: we want to have the $\UP$ for {\em
all} systems on $S$ (i.e. $\Pro_S$).

We work mostly without large cardinals. First we concentrate our attention on
the case when $\kappa=\lambda^+$, $\lambda$ a regular cardinal, and then we
deal with the related problem for inaccessible $\kappa$. The following four
cases should be treated somewhat separately. 
\begin{description}
\item[Case A:\quad] $\kappa=\lambda^+$, $\lambda=\lambda^{<\lambda}$, the set
$S^\kappa_\lambda\setminus S$ is stationary; 
\item[Case B:\quad] $\kappa$ is (strongly) inaccessible (e.g. the first one),
$S$ is ``thin'';
\item[Case C:\quad] $\lambda$ is singular, $S\subseteq S^{\lambda^+}_{
\cf(\lambda)}$ is a non--reflecting stationary set; 
\item[Case D:\quad] $\kappa$ is strongly inaccessible, the set
$\{\delta<\kappa: \delta\notin S$ and $\delta$ is strongly inaccessible$\}$\ \
\ is stationary; 
\item[Case E:\quad] $S=S^\kappa_\lambda$, $\kappa=\lambda^+$.
\end{description}
In the present paper we will deal with the first two (i.e. {\bf A} and
{\bf B}) cases. The other three cases will be considered in a subsequent paper
\cite{Sh:F259}.
 
Note that $\diamondsuit_S$ excludes the Uniformization Property for systems on
$S$. Consequently we have some immediate limitations and restrictions. Because
of a theorem of Jensen, in case {\bf B} we have to consider $S\subseteq\kappa$
which is not too large (e.g. not reflecting). In context of case {\bf C}, one
should remember that by Gregory \cite{Gre} and \cite{Sh:108}:\quad if
$\lambda^{<\lambda}=\lambda$ or $\lambda$ is strong limit singular, $2^\lambda
=\lambda^+$ and $S\subseteq\{\delta<\lambda^+: \cf(\delta)\neq\cf(\lambda)\}$
is stationary, then $\diamondsuit_S$ holds true.\\ 
By \cite[\S3]{Sh:186}, if $\lambda$ is a strong limit singular cardinal,
$2^\lambda=\lambda^+$, $\Box_\lambda$ and $S\subseteq\{\delta<\lambda^+:
\cf(\delta)=\cf(\lambda)\}$ reflects on a stationary set then $\diamondsuit_S$
holds; more results in this direction can be found in \cite{CDSh:571}.\\
In the cases {\bf A}, {\bf E} we will additionally assume that $\lambda^{
<\lambda}=\lambda$. We will start with the first (i.e. {\bf A}) case which
seems to be easier. The forcing notions which we will use will be quite
complete (see \ref{complete}, \ref{good}, \ref{really} below). Having this
amount of completeness we will be able to put weaker requirements on the
forcing notion for $S$. 

Finally note that we cannot expect here a full parallel of properness for
$\lambda=\aleph_0$, as even for $\lambda^+$--cc the parallel of {\em FS
iteration preserves ccc} fails.

We deal here with cases A and B, the other will appear in Part II,
\cite{Sh:F259}. For iterating $(<\lambda)$--complete forcing notions possibly
adding subsets to $\lambda$, $\kappa=\lambda^+$, see \cite{Sh:F249}. In
\cite{Sh:F259} we show a weaker $\kappa^+$--cc (parallel to pic, eec in
\cite[Ch VII, VIII]{Sh:f}) suffices. We also show that for a strong
limit singular $\lambda$ cardinal and a stationary set $S\subseteq S^{
\lambda^+}_{\cf(\lambda)}$, $\Pro_S$ (the uniformization for $S$) fails, but
it may hold for many $S$--ladder systems (so we have consequences for the
Whitehead groups).

This paper is based on my lectures in Madison, Wisconsin, in February and
March 1996, and was written up by Andrzej Roslanowski to whom I am greatly
indebted. 

\subsection{Background: Abelian groups}
We try to be self contained, but if you are lost see Eklof Mekler \cite{EM}.
\begin{definition}
\label{Wgroup}
1)\quad An Abelian group $G$ is {\em a Whitehead group} if for every
homomorphism $h:H\stackrel{\rm onto}{\longrightarrow} G$ from a group $H$ onto
$G$ such that $\Ker(h)\equiv {\Bbb Z}$ there is a lifting $g$ (i.e.~a
homomorphism $g:G\longrightarrow H$ such that $h\comp g={\rm id}_G$).

\noindent 2)\quad Let $h:H\longrightarrow G$ be as above, $G_1$ be a subgroup
of $G$. A homomorphism $g:G_1\longrightarrow H$ is a lifting for $G_1$ (and
$h$) if $h\comp g_1={\rm id}_{G_1}$.

\noindent 3)\quad We say that an Abelian group $G$ is a {\em direct sum} of
its subgroups $\langle G_i: i\in J\rangle$ (and then we write $G=
\bigoplus\limits_{i\in J} G_i$) if 
\begin{description}
\item[(a)] $G=\langle\bigcup\limits_{i\in J} G_i\rangle_G$ (where for a set
$A\subseteq G$, $\langle A\rangle_G$ is the subgroup of $G$ generated by $A$;
$\langle A\rangle_G=\{\sum\limits_{\ell<k}a_\ell x_\ell:k<\omega,\ a_\ell\in
{\Bbb Z},\ x_\ell\in A\}$), and
\item[(b)] $G_i\cap\langle \bigcup\limits_{i\neq j} G_j\rangle_G =\{0_G\}$ for
every $i\in J$.   
\end{description}
\end{definition}

\begin{remark}
{\em 
Concerning the definition of a Whitehead group, note that if $h:H
\longrightarrow G$ is a homomorphism such $\Ker(h)={\Bbb Z}$ and $H={\Bbb
Z}\oplus H_1$ then $h\restriction H_1$ is a homomorphism from $H_1$ into $G$
with kernel $\{0\}$ (and so it is one-to-one, and ``onto''). Thus $h\rest H_1$
is an isomorphism and $g\stackrel{\rm def}{=} (h\rest H_1)^{-1}$ is a required
lifting. 

Also conversely, if $g:G\longrightarrow H$ is a homomorphism such that $h\comp
g={\rm id}_G$ then $H={\Bbb Z}\oplus g[G]$.

The reader familiar with the Abelian group theory should notice that a group
$G$ is Whitehead if and only if ${\rm Ext}(G,{\Bbb Z})=\{0\}$.
}
\end{remark}

\begin{proposition}
\begin{enumerate}
\item If $h:H\longrightarrow G$ is a homomorphism, $G_1\oplus G_2\subseteq G$
and $g_\ell$ is a lifting for $G_\ell$ (for $\ell=1,2$) then there is a unique
lifting $g$ for $G_1\oplus G_2$ (called $\langle g_1,g_2\rangle$) extending
both $g_1$ and $g_2$ and such that $g(x_1+x_2)=g_1(x_1)+g_2(x_2)$ whenever
$x_1\in G_1$, $x_2\in G_2$. 
\item Similarly for $\bigoplus\limits_{i\in J}G_i$, $g_i$ a lifting for $G_i$.
\item If $\Ker(h)\cong {\Bbb Z}$ and $G_1\subseteq G$ is isomorphic to ${\Bbb
Z}$ then there is a lifting for $G_i$.
\end{enumerate}
\end{proposition}

\begin{definition}
\label{almfree}
Let $\lambda$ be an uncountable cardinal, $G$ be an Abelian group. 
\begin{description}
\item[(a)] $G$ is {\em free} if and only if $G=\bigoplus\limits_{i\in J} G_i$
where each $G_i$ is isomorphic to ${\Bbb Z}$. 
\item[(b)] $G$ is {\em $\lambda$--free} if every its subgroup of size
$<\lambda$ is free.  
\item[(c)] $G$ is {\em strongly $\lambda$--free} if for every $G'\subseteq G$
of size $<\lambda$ there is $G''$ such that 
\begin{description}
\item[$(\alpha)$] $G'\subseteq G''\subseteq G$ and $|G''|<\lambda$,
\item[$(\beta)$]  $G''$ is free,
\item[$(\gamma)$] $G/G''$ is $\lambda$--free.
\end{description}
\item[(d)] $G$ is {\em almost free in  $\lambda$} if it is strongly
$\lambda$--free of cardinality $\lambda$ but it is not free.
\end{description}
\end{definition}

\begin{remark}
{\em 
Note that the {\em strongly} in \ref{almfree}(d) does not have much
influence. In particular, for $\kappa$ strongly inaccessible,  ``strongly
$\kappa$--free'' is equivalent to ``$\kappa$--free''.
}
\end{remark}

\begin{proposition}
\label{getL}
Assume $G/G''$ is $\lambda$--free. Then for every $K\subseteq G$, $|K|<
\lambda$ there is a free group $L\subseteq G$ such that $K\subseteq G''\oplus
L\subseteq G$.
\end{proposition}

Suppose that $G$ is an almost free in $\kappa$ Abelian group of size $\kappa$,
and assume that $\kappa$ is a regular limit cardinal. Let $\langle
G_i:i<\kappa\rangle$ be a filtration of $G$, i.e.~$\langle
G_i:i<\kappa\rangle$ is an increasing continuous sequence of subgroups of $G$,
each of size less than $\kappa$. Let
\[\gamma(\bar{H})=\{i<\kappa: H/H_i \mbox{ is not $\kappa$--free}\},\]
and let $\Gamma[H]=\gamma(\bar{H})/D_\kappa$ for any filtration $\bar{H}$
(see \cite{EM}).\\

\begin{proposition}
Suppose that $G$, $\kappa$ and $\langle G_i:i<\kappa\rangle$ are as above.
\begin{enumerate}
\item $G$ is free if and only if $\Gamma(G)$ is not stationary. 
\item $\Gamma[G]$ cannot reflect in inaccessibles.
\end{enumerate}
\end{proposition}

The problem which was the {\em reason detre} of the paper is the following
question of G\"obel.
 
\begin{Gquest}
\label{Gques}
Is it consistent with {\em GCH} that for some regular cardinal $\kappa$ we
have
\begin{description}
\item[(a)] every almost free in $\kappa$ Abelian group is Whitehead, and 
\item[(b)] there are almost free in $\kappa$ Abelian groups ?
\end{description}
\end{Gquest}
 
\begin{remark}
{\em 
The point in \ref{Gques}(b) is that without it we have a too easy solution:
any weakly compact cardinal will do the job. This demand is supposed to be a
complement of G\"obel Shelah \cite{GbSh:579} which proves that, say for
$\kappa=\aleph_n$, there are (under GCH) almost free in $\kappa$ groups $H$
with $\HOM(H,{\Bbb Z})=\{0\}$.
}
\end{remark}

Now, our conclusion \ref{gobelcon} gives that 
\begin{description}
\item[(a)']  every almost free in $\kappa$ Abelian group $G$ with $\Gamma[G] 
\subseteq S/D_\kappa$ is Whitehead,
\item[(b)']  there are almost free in $\kappa$ Abelian groups $H$ with
$\Gamma[H]\subseteq S/D_\kappa$.
\end{description}
It can be argued that this answers the question if we understand it as whether
from an almost free in $\kappa$ Abelian group we can build a non--Whitehead
one, so the further restriction of the invariant to be $\subseteq S$.

However we can do better, starting with a weakly compact cardinal $\kappa$ we
can manage that in addition to {\bf (a)'}, {\bf (b)'} we have
\begin{description}
\item[(b)$^+$]  {\bf (i)}\ \ \ every stationary subset of $\kappa\setminus S$
reflects in inaccessibles, 

{\bf (ii)}\ \ \ moreover, for every almost free in $\kappa$ Abelian group $H$,
$\Gamma[H]\subseteq S/D_\kappa$
\end{description}
(in fact for an uncountable strongly inaccessible $\kappa$, {\bf (i)} implies
{\bf (ii)}). So we get a consistency proof for the original problem. This we
will do here.

We may ask, can we do it for small cardinals? Successor of singular? Successor
of regular? For many cardinals simultaneously? We may get consistency and
ZFC $+$ GCH information, but always the consistency strength is not
small. That is, we need a regular cardinal $\kappa$ and a stationary set
$S\subseteq\kappa$ such that we have enough uniformization on $S$. Now, for
Whitehead group: if $\bar{G}=\langle G_i:i<\kappa\rangle$ is a filtration of
$G$, $S=\gamma(\bar{G})$, $\lambda_i=|G_{i+1}/G_i|$ for $i\in S$, for
simplicity $\lambda_i=\lambda$, then we need a version of
$\Pro_{S,\lambda}$. We would like to have a suitable reflection (see Magidor
Shelah \cite{MgSh:204}); for a stationary $S'\subseteq\kappa\setminus S$ this
will imply $0^\#$. 

\subsection{Background: forcing} Let us review some basic facts concerning
iterated forcing and establish our notation. First remember that in forcing
considerations we keep the convention that
\begin{center}
{\em
a stronger condition (i.e. carrying more information) is the larger one.
}
\end{center}
For more background than presented here we refer the reader to either \cite[Ch
4]{J} or \cite{Sh:f}. 
\begin{definition}
Let $\kappa$ be a cardinal number. We say that $\bar{\bQ}$ is {\em a
$(<\kappa)$--support iteration of length $\gamma$ (of forcing notions
$\nbQ_\alpha$)} if $\bar{\bQ}=\langle\bP_\alpha,\nbQ_\beta:\alpha\leq\gamma,\
\beta<\gamma\rangle$  and for every $\alpha\leq\gamma$, $\beta<\gamma$:
\begin{description}
\item[(a)] $\bP_\alpha$ is a forcing notion,
\item[(b)] $\nbQ_\beta$ is a $\bP_\beta$--name for a forcing notion with the
minimal element ${\bf 0}_{\nbQ_\beta}$ 

[for simplicity we will assume that $\nbQ_\beta$ is a partial order on an
ordinal; remember that each partial order is isomorphic to one of this form],
\item[(c)] a condition $f$ in $\bP_\alpha$ is a partial function such that
$\dom(f)\subseteq\alpha$, $\|\dom(f)\|<\kappa$ and 
\[(\forall\xi\in\dom(f))(f(\xi)\mbox{ is a $\bP_\xi$--name and }\
\forces_{\bP_\xi} f(\xi)\in\nbQ_\xi)\]  

[we will keep a convention that if $f\in\bP_\alpha$, $\xi\in\alpha\setminus
\dom(f)$ then $f(\xi)={\bf 0}_{\nbQ_\xi}$; moreover we will assume that each
$f(\xi)$ is a canonical name for an ordinal, i.e. $f(\xi)=\{\langle
q_i,\gamma_i\rangle: i<i^*\}$ where $\{q_i:i<i^*\}\subseteq\bP_\xi$ is a
maximal antichain of $\bP_\xi$ and for every $i<i^*$: $\gamma_i$ is an ordinal
and $q_i\forces_{\bP_\xi}$``$f(\xi)=\gamma_i$''],
\item[(d)] the order of $\bP_\alpha$ is given by

$f_1\leq_{\bP_\alpha} f_2$\quad if and only if\quad $(\forall \xi\in\alpha)(
f_2\rest\xi\forces f_1(\xi)\leq_{\nbQ_\xi} f_2(\xi))$.
\end{description}
\end{definition}
Note that the above definition is actually an inductive one (see below too).

\begin{remark}
{\em 
The forcing notions which we will consider will satisfy {\em no new sequences
of ordinals of length $<\kappa$ are added}, or maybe at least {\em any new set
of ordinals of cardinality $<\kappa$ is included in an old one}. Therefore
there will be no need to consider the revised support iterations.
}
\end{remark}
Let us recall that:
\begin{fact}
Suppose $\bar{\bQ}=\langle\bP_\alpha,\nbQ_\beta:\alpha\leq\gamma,\ \beta<
\gamma\rangle$ is a $(<\kappa)$--support iteration, $\beta<\alpha\leq\gamma$.
Then 
\begin{description}
\item[(a)] $p\in\bP_{\alpha}$ implies $p\rest\beta\in\bP_\beta$,
\item[(b)] $\bP_\beta\subseteq\bP_\alpha$,
\item[(c)] $\leq_{\bP_\beta}=\leq_{\bP_\alpha}\rest \bP_\beta$,
\item[(d)] if $p\in\bP_\alpha$, $p\rest\beta\leq_{\bP_\beta} q\in\bP_\beta$
then the conditions $p,q$ are compatible in $\bP_\alpha$; in fact $q\cup
p\rest [\beta,\alpha)$ is the least upper bound of $p,q$ in $\bP_\alpha$,

consequently
\item[(e)] $\bP_\beta\lesdot\bP_\alpha$ (i.e. complete suborder). \QED
\end{description}
\end{fact}

\begin{fact}
\label{fakcik}
\begin{enumerate}
\item If $\bar{\bQ}=\langle\bP_\alpha,\nbQ_\beta:\alpha\leq\gamma,\
\beta<\gamma\rangle$ is a $(<\kappa)$--support iteration of length $\gamma$,
$\nbQ_\gamma$ is a $\bP_\gamma$--name for a forcing notion (on an ordinal)

then there is a unique $\bP_{\gamma+1}$ such that $\langle\bP_\alpha,
\nbQ_\beta:\alpha\leq\gamma+1,\ \beta<\gamma+1\rangle$ is a
$(<\kappa)$--support iteration. 
\item If $\langle\gamma_i: i<\delta\rangle$ is a strictly increasing continuous
sequence of ordinals with the limit $\gamma_\delta$, $\delta$ is a limit
ordinal, and for each $i<\delta$ the sequence $\langle \bP_\alpha,\nbQ_\beta:
\alpha\leq\gamma_i,\ \beta<\gamma_i\rangle$ is a $(<\kappa)$--support iteration

then there is a unique $\bP_{\gamma_\delta}$ such that $\langle\bP_\alpha,
\nbQ_\beta:\alpha\leq\gamma_\delta,\ \beta<\gamma_\delta\rangle$ is a
$(<\kappa)$--support iteration. \QED
\end{enumerate}
\end{fact}
Because of the above fact \ref{fakcik}(2) we may write $\bar{\bQ}=\langle
\bP_\alpha,\nbQ_\alpha:\alpha<\gamma\rangle$ when considering iterations (with
$(<\kappa)$--support), as $\bP_\gamma$ is determined by it (for $\gamma=\beta
+1$ essentially $\bP_\gamma=\bP_\beta*\nbQ_\beta$). For $\gamma'<\gamma$ and
an iteration $\bar{\bQ}=\langle\bP_\alpha,\nbQ_\alpha:\alpha<\gamma\rangle$ we
let 
\[\bar{\bQ}\rest\gamma'=\langle\bP_\alpha,\nbQ_\alpha:\alpha<\gamma'\rangle.\]

\begin{fact}
For every function ${\bf F}$ (even a class) and an ordinal $\gamma$ there is a
unique $(<\kappa)$--support iteration $\bar{\bQ}=\langle\bP_\alpha,
\nbQ_\alpha: \alpha<\gamma'\rangle$, $\gamma'\leq\gamma$ such that $\nbQ_\alpha
={\bf F}(\bar{\bQ}\rest\alpha)$ for every $\alpha<\gamma'$ and
\[\mbox{either }\quad \gamma'=\gamma\quad\mbox{ or }\quad {\bf F}(\bar{\bQ})\
\mbox{ is not of the right form.}\quad\QED\]
\end{fact}

For a forcing notion $\bQ$, the completion of $\bQ$ to a complete forcing will
be denoted by $\hat{\bQ}$ (see \cite[Ch XIV]{Sh:f}). Thus $\bQ$ is a dense
suborder of $\hat{\bQ}$ and in $\hat{\bQ}$ any increasing sequence of
conditions which has an upper bound has the least upper bound. In this context
note that we may define and proof by induction on $\alpha^*$ the following
fact. 

\begin{fact}
\label{uzupelnienie}
Assume $\langle\bP_\alpha',\nbQ_\alpha': \alpha<\alpha^*\rangle$ is a
$(<\kappa)$--support iteration. Let $\bP_\alpha$, $\nbQ_\alpha$ be such that
for $\alpha<\alpha^*$ 
\begin{enumerate}
\item $\bP_\alpha=\{f\in\bP_\alpha': (\forall\xi<\alpha)(f(\xi)\mbox{ is
a $\bP_\xi$--name for an element of }\nbQ_\xi)\}$, 
\item $\nbQ_\alpha$ is a $\bP_\alpha$--name for a dense suborder of
$\nbQ_\alpha'$. 
\end{enumerate}
Then for each $\alpha\leq\alpha^*$, $\bP_\alpha$ is a dense suborder of
$\bP_\alpha'$ and $\langle\bP_\alpha,\nbQ_\alpha: \alpha<\alpha^*
\rangle$ is a $(<\kappa)$--support iteration. \QED
\end{fact}
We finish our overview of basic facts with the following observation, which
will be used several times later (perhaps even without explicit referring to
it). 
\begin{fact}
\label{seqfromN}
Let $\bQ$ be a forcing notion which does not add new $({<}\theta)$--sequences
of elements of $\lambda$ (i.e.~$\forces_{\bQ}$``$\lambda^{<\theta}=\lambda^{
<\theta}\cap\V$''). Suppose that $N$ is an elementary submodel of $\Hchi$ such
that $\|N\|=\lambda$ and $N^{<\theta}\subseteq N$. Let $G\subseteq\bQ$ be a
generic filter over $\V$. Then 
\[\V[G]\models N[G]^{<\theta}\subseteq N[G].\]
\end{fact}

\Proof Suppose that $\bar{x}=\langle x_i: i<i^*\rangle\in N[G]^{<\theta}$,
$i^*<\theta$. By the definition of $N[G]$, for each $i<i^*$ there is a
$\bQ$--name $\name{\tau}_i\in N$ such that $x_i=\name{\tau}^G_i$. Look at the
sequence $\langle\name{\tau}_i:i<i^*\rangle\in \V[G]$. By the assumptions on
$\bQ$ we know that $\langle\name{\tau}_i:i<i^*\rangle\in\V$ (remember $i^*<
\theta$, $\|N\|=\lambda$) and therefore, as each $\name{\tau}_i$ is in $N$ and
$N^{<\theta}\subseteq N$, we have $\langle\name{\tau}_i: i<i^*\rangle\in
N$. This implies that $\bar{x}\in N[G]$. \QED

\subsection{Notation}
We will define several properties of forcing notions using $\Hchi$. In all
these definitions any ``large enough'' regular cardinal $\chi$ works. 

\begin{definition}
\label{most}
{\em For most $N\prec\Hchi$ with PROPERTY we have$\ldots$} will mean:
\begin{description}
\item[] there is $x\in {\cal H}(\chi)$ such that 

if $x\in N\prec\Hchi$ and $N$ has the {\em PROPERTY} then $\ldots$.
\end{description}
Similarly, {\em for most sequences $\bar{N}=\langle N_i: i<\alpha\rangle$ of
elementary submodels of $\Hchi$ with PROPERTY we have$\ldots$} will mean: 
\begin{description}
\item[] there is $x\in {\cal H}(\chi)$ such that 

if $x\in N_0$, $\bar{N}=\langle N_i: i<\alpha\rangle$, $N_i\prec\Hchi$ and
$\bar{N}$ has the {\em PROPERTY} then $\ldots$. 
\end{description}
In these situations we call the element $x\in{\cal H}(\chi)$ {\em a witness}. 
\end{definition}

\begin{notation}
We will keep the following rules for our notation:
\begin{enumerate}
\item $\alpha,\beta,\gamma,\delta,\xi,\zeta, i,j\ldots$ will denote ordinals,
\item $\kappa,\lambda,\mu,\mu^*\ldots$ will stand for cardinal numbers,
\item a bar above a name indicates that the object is a sequence, usually
$\bar{X}$ will be $\langle X_i: i<\lh(\bar{X})\rangle$, where $\lh(\bar{X})$
denotes the length of $\bar{X}$,
\item a tilde indicates that we are dealing with a name for an object in
forcing extension (like $\name{x}$),
\item $S$, $S_i$, $S^j_i$, $E$, $E_i$, $E^j_i\ldots$ will be used to denote
sets of ordinals,  
\item $\family$, $\family_i$, $\family^j_i$, $\cal E$, ${\cal E}_i$, ${\cal
E}^j_i\ldots$ will stand for families of sets of ordinals of size $<\kappa$,
and finally  
\item $\hat{\family}$, $\hat{\family}_i$, $\hat{\family}^j_i$, $\baza$,
$\baza_i$, $\baza^j_i$ will stand for families of sequences of sets of
ordinals of size $<\kappa$. 
\item In groups we will use the additive convention (so in particular $0_G$
will stand for the neutral element of the group $G$). $G,H,K,L$ will denote
(always Abelian) groups.

\end{enumerate}
\end{notation}

\PART{A}{Case}

\noindent In this part of the paper we are dealing with the {\bf Case A} (see
the introduction), so naturally we assume the following. 

\begin{ourassum}
$\lambda$, $\kappa$, $\mu^*$ are uncountable cardinal numbers such that
\[\lambda^{<\lambda}=\lambda<\lambda^+=2^\lambda=\kappa\leq\mu^*.\]
We will keep these assumptions for some time (unless otherwise stated) and we
may forget to remind them.
\end{ourassum}

\section{Complete forcing notions}
\label{caseA}
In this section we introduce several notions of completeness of forcing
notions and prove basic results about them. We define when a forcing notion
$\bQ$ is: $(\theta,S)$--strategically complete, strongly $\family$--complete,
basically $(\family_0,\hat{\family}_1)$--complete and really $(\family_0,
\hat{\family}_1, D)$--complete. The notions which we will use are strong
$\family_0$--completeness and real
$(\family_0,\hat{\family}_1,D)$--completeness, however the other definitions 
seem to be interesting too. They are, in some sense, successive approximations
to real completeness (which is as weak as the iteration theorem allows) and
they might be of some interest in another contexts. But a reader not
interested in a general theory may concentrate on definitions
\ref{complete}(3), \ref{dkldef}, \ref{good}(3) and \ref{really} only.

\begin{definition}
\label{complete}
Let $\bQ$ be a forcing notion, $\theta$ be an ordinal and let $S\subseteq
\theta$. 
\begin{enumerate}
\item For a condition $r\in\bQ$, let $\Gra(\bQ,r)$ be the following game of two
players, $\com$ (for {\em complete}) and $\inc$ (for {\em incomplete}):
\begin{quotation}
\noindent the game lasts $\theta$ moves and during a play the players
construct a sequence $\langle (p_i,q_i): i<\theta\rangle$ of conditions from
$\bQ$ in such a way that $(\forall j<i<\theta)(r\leq p_j\leq q_j\leq p_i)$ and
at the stage $i<\theta$ of the game

\noindent if $i\in S$ then $\com$ chooses $p_i$ and $\inc$ chooses $q_i$, and

\noindent if $i\notin S$ then $\inc$ chooses $p_i$ and $\com$ chooses $q_i$.
\end{quotation}
The player $\com$ wins if and only if for every $i<\theta$ there are legal
moves for both players.
\item We say that the forcing notion $\bQ$ is {\em $(\theta,S)$--strategically
complete} if the player {\rm COM} has a winning strategy in the game
$\Gra(\bQ,r)$ for each condition $r\in\bQ$. We say that $\bQ$ is {\em
strategically $(<\theta)$--complete} if it is
$(\theta,\emptyset)$--strategically complete.

\item We say that the forcing notion $\bQ$ is {\em $(<\theta)$--complete} if
every increasing sequence $\langle q_i: i<\delta\rangle\subseteq\bQ$ of length
$\delta<\theta$ has an upper bound in $\bQ$.
\end{enumerate}
\end{definition}

\begin{proposition}
\label{implik}
Suppose that $\bQ$ is a forcing notion, $\theta$ is an ordinal and let
$S\subseteq\theta$. 
\begin{enumerate}
\item If $\bQ$ is $(<\theta)$--complete then it is $(\theta,S)$--strategically
complete. 
\item If $S'\subseteq S$ and $\bQ$ is $(\theta,S')$--strategically complete
then it is $(\theta,S)$--strategically complete. 
\item If $\bQ$ is $(\theta,S)$--strategically complete forcing notions then
the forcing with $\bQ$ does not add new sequences of ordinals of length
$<\theta$. 
\end{enumerate}
\end{proposition}

\Proof 1) and 3) should be clear.\\
2)\quad Note that if all members of $S$ are limit ordinals, or at least
$\alpha\in S\quad \Rightarrow\quad\alpha+1\notin S$, then one may easily
translate a winning strategy for $\com$ in ${\cal G}^\theta_{S'}(\bQ,r)$ to
the one in $\Gra(\bQ,r)$. In a general case, however, we have to be slightly
more careful. First note that we may assume that $\theta$ is a limit ordinal
(if $\theta$ is not limit consider the game ${\cal G}^{\theta+\omega}_S(\bQ,
r)$). Now, for a set $S\subseteq\theta$ and a condition $r\in\bQ$ we define a
game ${}^*{\cal G}^\theta_S(\bQ,r)$:
\begin{quotation}
\noindent the game lasts $\theta$ moves and during a play the players
construct a sequence $\langle p_i: i<\theta\rangle$ of conditions from $\bQ$
such that $r\leq p_i\leq p_j$ for each $i<j<\theta$ and

if $i\in S$ then $p_i$ is chosen by $\com$,

if $i\notin S$ then $p_i$ is determined by $\inc$.

\noindent The player $\com$ wins if and only if there are legal moves for each
$i<\theta$.
\end{quotation}
Note that clearly, if $S'\subseteq S\subseteq\theta$ and Player $\com$ has a
winning strategy in ${}^*{\cal G}^\theta_{S'}(\bQ,r)$ then he has one in
${}^*{\cal G}^\theta_S(\bQ,r)$.

For a set $S\subseteq\theta$ let $S^\bot=\{2\alpha: \alpha\in S\}\cup
\{2\alpha+1: \alpha\in\theta\setminus S\}$. (Plainly $S^\bot\subseteq \theta$
as $\theta$ is limit.)

\begin{claim}
\label{cl16}
For each set $S\subseteq\theta$ the games $\Gra(\bQ,r)$ and ${}^*{\cal
G}^\theta_S(\bQ,r)$ are equivalent [i.e.~$\com$/$\inc$ has a winning strategy
in $\Gra(\bQ,r)$ if and only if he has one in ${}^*{\cal G}^\theta_{S^\bot}(
\bQ,r)$]. 
\end{claim}

\noindent{\em Proof of the claim:}\hspace{0.15in} Look at the definitions of
the games and the set $S^\bot$.

\begin{claim}
\label{cl17}
Suppose that $S_0,S_1\subseteq\theta$ are such that for every non--successor
ordinal $\delta<\theta$ we have
\begin{description}
\item[(a)] $\delta\in S_0\quad\equiv\quad \delta\in S_1$,
\item[(b)] $(\exists^\infty n\in\omega)(\delta+n\in S_0)$, $(\exists^\infty
n\in\omega)(\delta+n\notin S_0)$, $(\exists^\infty n\in\omega)(\delta+n\in
S_1)$, and $(\exists^\infty n\in\omega)(\delta+n\notin S_0)$.
\end{description}
Then the games ${}^*{\cal G}^\theta_{S_0}(\bQ,r)$ and ${}^*{\cal
G}^\theta_{S_1}(\bQ,r)$ are equivalent.
\end{claim}

\noindent{\em Proof of the claim:}\hspace{0.15in} Should be clear once you
realize that finitely many successive moves by the same player may be
interpreted as one move.
\medskip

Now we may prove \ref{implik}(2). Let $S'\subseteq S\subseteq\theta$ (and
$\theta$ be limit). Let $S^*=\{\delta\in S^\bot: \delta$ is not a successor$\;
\}\cup\{\delta\in (S')^\bot: \delta$ is a successor$\;\}$. Note that
$(S')^\bot\subseteq S^*$ and the sets $S^*,S^\bot$ satisfy the demands {\bf
(a)}, {\bf (b)} of \ref{cl17}. Consequently, by \ref{cl16} and \ref{cl17}:
\[\begin{array}{l}
\mbox{Player }\com\mbox{ has a winning strategy in }{\cal G}^\theta_{S'}(\bQ,r)
\quad\Rightarrow\\
\mbox{Player }\com\mbox{ has a winning strategy in }{}^*{\cal
G}^\theta_{(S')^\bot}(\bQ,r)\quad\Rightarrow\\ 
\mbox{Player }\com\mbox{ has a winning strategy in }{}^*{\cal G}^\theta_{S^*}
(\bQ,r)\quad\Rightarrow\\
\mbox{Player }\com\mbox{ has a winning strategy in }{}^*{\cal
G}^\theta_{S^\bot}(\bQ,r)\quad\Rightarrow\\
\mbox{Player }\com\mbox{ has a winning strategy in }\Gra(\bQ,r).\qquad
\QED 
  \end{array}\]

\begin{proposition}
\label{firstiter}
Assume $\kappa$ is a regular cardinal and $\theta\leq\kappa$. Suppose that
$\bar{\bQ}=\langle\bP_\alpha,\nbQ_\alpha:\alpha<\gamma\rangle$ is a $(<
\kappa)$--support iteration of $(<\theta)$--complete
($(\theta,S)$--strategically complete, strategically $(<\theta)$--complete,
respectively) forcing notions. Then $\bP_\gamma$ is $(<\theta)$--complete
($(\theta,S)$--strategically complete, strategically $(<\theta)$--complete,
respectively). 
\end{proposition}

\Proof Easy if you remember that union of less than $\kappa$ sets of size
less that $\kappa$ is of size $<\kappa$, and use \ref{implik}(3). \QED

Note that if we pass from a $(<\lambda)$--complete forcing notion $\bQ$ to its 
completion $\hat{\bQ}$ we may loose $(<\lambda)$--completeness. However, a
large amount of the completeness is preserved. 

\begin{proposition}
\label{denssub}
Suppose that $\bQ$ is a dense suborder of $\bQ'$.
\begin{enumerate}
\item If $\bQ$ is $(<\lambda)$--complete (or just $(<\lambda)$--strategically
complete) then $\bQ'$ is $(<\lambda)$--strategically complete.
\item If $\bQ'$ is $(<\lambda)$--strategically complete then so is $\bQ$.
\end{enumerate}
\end{proposition}

\Proof 1)\ \ \ We describe a winning strategy for player $\com$ in the game
${\cal G}^\lambda_\emptyset(\bQ',r)$ ($r\in\bQ'$), such that it says player
$\com$ to choose elements of $\bQ$ only. So 
\begin{quotation}
\noindent at stage $i<\lambda$ of the play, $\com$ chooses the
$<^*_\chi$--first condition $q_i\in\bQ$ stronger than $p_i\in\bQ'$ chosen by
$\inc$ right before.
\end{quotation}
This strategy is the winning one, as at a limit stage $i<\lambda$ of the play,
the sequence $\langle q_j: j<i\rangle$ has an upper bound in $\bQ$ (remember
$\bQ$ is $(<\lambda)$--complete). 

\noindent 2)\ \ \ Even easier. \QED
 
\begin{definition}
\label{dkldef}
\begin{enumerate}
\item By $\dkl$ we will denote the collection of all families
$\family\subseteq [\mu^*]^{<\kappa}$ such that for every large enough regular
cardinal $\chi$, for some $x\in{\cal H}(\chi)$ we have

\quad if $x\in N\prec\Hchi$, $\|N\|<\kappa$, $N^{<\lambda}\subseteq N$ and
$N\cap\kappa$ is an ordinal,

\quad then $N\cap\mu^*\in\family$ (compare with \ref{good}).\\
If $\lambda=\aleph_0$ then we may omit it.

\item By $\dak$ we will denote the collection of all sets $\hat{\family}$ such
that 
\[\begin{array}{r}
\hat{\family}\subseteq\big\{\bar{a}=\langle a_i: i\leq\alpha\rangle: \mbox{
the sequence $\bar{a}$ is increasing continuous and }\ \\
\mbox{ each $a_i$ is from $[\mu^*]^{<\kappa}$}\big\}\end{array}\]
and for every large enough regular cardinal $\chi$, for some $x\in{\cal
H}(\chi)$ we have

\quad if $\bar{N}=\langle N_i: i\leq\alpha\rangle$ is increasing continuous,
$x\in N_0$ and for each $i<j\leq\alpha$: \qquad $N_i\prec N_j\prec\Hchi$,\
$\langle N_\xi:\xi\leq j\rangle\in N_{j+1}$,\ $\|N_j\|<\kappa,\quad N_j\cap
\kappa\in\kappa$ and 
\[j\mbox{ is non-limit }\quad\Rightarrow\quad {N_j}^{<\lambda}\subseteq N_j,\]
\quad then $\langle N_i\cap\mu^*:i\leq\alpha\rangle\in\hat{\family}$.

\item For a family ${\frak D}\subseteq {\cal P}({\cal X})$ (say ${\cal
X}=\bigcup\limits_{X\in {\frak D}} X$) let ${\frak D}^+$ stand for the family
of all $\family\subseteq {\cal X}$ such that  
\[(\forall {\cal C}\in{\frak D})(\family\cap {\cal C}\neq\emptyset).\]
[So ${\frak D}^+$ is the collection of all ${\frak D}$--positive subsets of
${\cal X}$.]

\item For $\family\in(\dkl)^+$ we define $\dak[\family]$ like $\dak$ above,
except that its members $\hat{\family}$ are subsets of
\[\begin{array}{r}
\big\{\bar{a}=\langle a_i: i\leq\alpha\rangle: \bar{a}\ \mbox{ is increasing
continuous and for each }i\leq\alpha,\ \ \ \\ 
a_i\in [\mu^*]^{<\kappa}\ \mbox{ and if $i$ is not limit then }a_i\in\family
\big\},\end{array}\]
and, naturally, we consider only these sequences $\bar{N}=\langle N_i: i\leq
\alpha\rangle$ for which 
\[i\leq\alpha\mbox{ is non-limit }\quad\Rightarrow\quad N_i\cap\mu^*\in
\family.\]
\end{enumerate}
As $\lambda$ is determined by $\kappa$ in our present case we may forget to
mention it.
\end{definition}

\begin{remark}
{\em
\begin{enumerate}
\item These are normal filters in natural senses.
\item Concerning $\dak$, we may not distinguish $\bar{a}_0$, $\bar{a}_1$ which
are similar enough (e.g. see \ref{really} below).
\item Remember: our case is GCH, $\lambda=\cf(\lambda)$, $\kappa=\lambda^+$
and $\alpha=\lambda$. 
\end{enumerate}
}
\end{remark}

\begin{definition}
\label{good}
Assume $\family\subseteq [\mu^*]^{\leq\lambda}$.
\begin{enumerate}
\item Let $\chi$ be a large enough regular cardinal. We say that an elementary
submodel $N$ of $\Hchi$ is {\em $(\lambda,\family)$--good} if 
\[\|N\|=\lambda,\quad N^{<\lambda}\subseteq N,\quad\mbox{ and }\ N\cap\mu^*\in
\family.\]
\item We say that a forcing notion $\bQ$ is {\em strongly $\family$--complete}
if for most (see \ref{most}) $(\lambda,\family)$--good elementary submodels
$N$ of $\Hchi$ such that $\bQ\in N$ and for each $\bQ$--generic over $N$
increasing sequence $\bar{p}=\langle p_i:i<\lambda \rangle \subseteq\bQ\cap N$
there is an upper bound of $\bar{p}$ in $\bQ$.  

[Recall that a sequence $\bar{p}=\langle p_i: i<\lambda\rangle\subseteq\bQ\cap
N$ is {\em $\bQ$--generic over $N$} if for every open dense subset ${\cal I}$
of $\bQ$ from $N$ for some $i<\lambda$, $p_i\in {\cal I}$.]
\item Let $N\prec\Hchi$ be $(\lambda,\family)$--good. For a forcing notion
$\bQ$, a set $S\subseteq\lambda$ and a condition $r\in\bQ\cap N$ we define a
game ${\cal G}_{N,S}(\bQ,r)$ like the game ${\cal G}^{\lambda+1}_S(\bQ,r)$
with an additional requirement that during a play all choices below $\lambda$
have to be done from $N$, i.e. $p_i,q_i\in N\cap\bQ$ for all $i<\lambda$.

If $S=\emptyset$ then we may omit it.
\item Let $\bar{S}:\family\longrightarrow {\cal P}(\lambda)$. We say that a
forcing notion $\bQ$ is {\em $(\family,\bar{S})$--complete} if for most
$(\lambda,\family)$--good models $N$, for every condition $r\in N\cap\bQ$
the player $\com$ has a winning strategy in the game ${\cal G}_{N,\bar{S}(N
\cap\mu^*)}(\bQ,r)$.

If $\bar{S}(a)=\emptyset$ for each $a\in\family$ then we write
{\em $\family$--complete}. (In both cases we may add ``strategically''.)
\end{enumerate}
\end{definition}

\begin{remark}
{\em 
In the use of {\em most} in \ref{good} (and later) we do not mention
explicitly the witness $x$ for it. And in fact, normally it is not
necessary. If $\chi_1,\chi$ are large enough, $2^{<\chi_1}<\chi$ (so ${\cal 
H}(\chi_1)\in {\cal H}(\chi)$), $\family,\bQ,\ldots\in N$ then there is a
witness in ${\cal H}(\chi_1)$ and therefore there is such a witness in
$N$. Consequently we may forget it.
}
\end{remark}

\begin{remark}
{\em 
\begin{enumerate}
\item Most popular choice of $\mu^*$ is $\kappa$; then $\family\in(\dkl)^+$ if
and only if the set $\{\delta<\kappa: \cf(\delta)=\lambda\ \&\ \delta\in
\family\}$ is stationary. So $\family$ ``becomes'' a stationary subset of
$\kappa$. 
\item Also here we have obvious monotonicities and implications. 
\end{enumerate}
}
\end{remark}

\begin{proposition}
\label{nonewseq}
Suppose that $\family\in (\dkl)^+$ and a forcing notion $\bQ$ is
$\family$--complete. Then the forcing with $\bQ$ adds no new
$\lambda$--sequences of ordinals (or, equivalently, of elements of $\V$) and
$\forces_{\bQ}$``$\family\in (\dkl)^+$''. 
\end{proposition}

\Proof Standard; compare the proof of \ref{iterstrong}. \QED

\begin{proposition}
\label{scomimpcom}
Let $\family\subseteq [\mu^*]^{\leq\lambda}$. If a forcing notion $\bQ$ is
strongly $\family$--complete and is $(<\lambda)$--complete, then it is
$\family$--complete. \QED
\end{proposition}

The strong $\family$--completeness is preserved if we pass to the completion
of a forcing notion.

\begin{proposition}
\label{strondense}
Suppose that $\family\subseteq [\mu^*]^{\leq\lambda}$ and $\bQ$ is a dense
suborder of $\bQ'$. Then
\begin{enumerate}
\item $\bQ'$ is strongly $\family$--complete if and only if $\bQ$ is strongly
$\family$--complete, 
\item similarly for $(\family,\bar{S})$--completeness.
\end{enumerate}
\end{proposition}

\Proof 1)\ \ \ Assume $\bQ'$ is strongly $\family$--complete and $x'\in {\cal
H}(\chi)$ be a witness for the ``most'' in the definition of this fact. Let
$x=\langle x',\bQ'\rangle$. Suppose that $N\prec\Hchi$ is
$(\lambda,\family)$--good and $\bQ,x\in N$. Then $\bQ',x'\in N$ too. Now
suppose that $\bar{q}=\langle q_i:i<\lambda\rangle\subseteq\bQ\cap N$ is an
increasing $\bQ$--generic sequence over $N$. Since $\bQ$ is dense in $\bQ'$,
$\bar{q}$ is $\bQ'$--generic over $N$ and thus, as $\bQ'$ is strongly
$\family$--complete, it has an upper bound in $\bQ'$ (and so in $\bQ$).

\noindent Now suppose $\bQ$ is strongly $\family$--complete with a witness
$x\in {\cal H}(\chi)$ and let $x'=\langle x,\bQ\rangle$. Let $N$ be
$(\lambda,\family)$--good and $\bQ',x'\in N$. So $\bQ,x\in N$. Suppose that
$\bar{q}=\langle q_i:i<\lambda\rangle\subseteq\bQ'\cap N$ is increasing and
$\bQ'$--generic over $N$. For each $i<\lambda$ choose a condition
$p_i\in\bQ\cap N$ and an ordinal $\varphi(i)<\lambda$ such that
\[q_i\leq_{\bQ'} p_i\leq_{\bQ'} q_{\varphi(i)}\]
(possible by the genericity of $\bar{q}$; remember that $\bar{q}$ is
increasing). Look at the sequence 
\[\langle p_i: i<\lambda\ \ \&\ \ (\forall j<i)(\varphi(j)<i)\rangle.\]
It is an increasing $\bQ$--generic sequence over $N$, so it has an upper bound
in $\bQ$. But this upper bound is good for $\bar{q}$ in $\bQ'$ as well. 

\noindent 2)\ \ \ Left to the reader. \QED

\begin{proposition}
\label{iterstrong}
Suppose that $\bar{\bQ}=\langle\bP_\alpha,\nbQ_\alpha:\alpha<\gamma\rangle$ is
a $(<\kappa)$--support iteration, $\family\in\V$, $\family\in (\dkl)^+$.
\begin{enumerate}
\item If for each $\alpha<\gamma$ 
\[\forces_{\bP_\alpha}\mbox{`` $\nbQ_\alpha$ is strongly $\family$--complete''
}\]
then the forcing notion $\bP_\gamma$ is strongly $\family$--complete (and even
each quotient $\bP_\beta/\bP_\alpha$ is strongly $\family$--complete for
$\alpha\leq\beta\leq\gamma$). 
\item Similarly for $(\family,\bar{S})$--completeness.
\end{enumerate}
\end{proposition}

\Proof 1)\ \ \ The proof can be presented as an inductive one (on $\gamma$),
so then we assume that each $\bP_\alpha$ ($\alpha<\gamma$) is strongly
$\family$--complete. However, the main use of the inductive hypothesis will be
that it helps to prove that no new sequences of length $\lambda$ are added
(hence $\lambda$ is not collapsed, so in $\V^{\bP_\alpha}$ (for
$\alpha<\gamma$) we may talk about $(\lambda,\family)$--good models without
worrying about the meaning of the definition if $\lambda$ is not a cardinal,
and $N[G_{\bP_\alpha}]$ is $(\lambda,\family)$--good). 

For each $\alpha<\gamma$ and $p\in\bP_\alpha$ fix a $\bP_\alpha$--name
$\name{f}^\alpha_p$ for a function from $\lambda$ to $\V$ such that
\begin{quotation}
\noindent if $p\forces_{\bP_\alpha}$``there is a new function from $\lambda$
to $\V$'' then $p\forces_{\bP_\alpha} \name{f}^\alpha_p\notin\V$, and otherwise
$p\forces_{\bP_\alpha}$``$\name{f}^\alpha_p$ is constantly  0''.
\end{quotation}
Let 
\[\begin{array}{lr}
{\cal I}_\alpha=\{p\in\bP_\alpha:&\mbox{ either $p\forces_{\bP_\alpha}$''there
is no new function from $\lambda$ to $\V$'' }\ \\
\ &\mbox{  or }p\forces_{\bP_\alpha}\name{f}^\alpha_p\notin\V\;\}.
  \end{array}\]
Clearly ${\cal I}_\alpha$ is an open dense subset of $\bP_\alpha$. Let
$\name{x}_\alpha$ (for $\alpha<\gamma$) be a $\bP_\alpha$--name for a witness
to the assumption that $\forces_{\bP_\alpha}$``$\nbQ_\alpha$ is strongly
$\family$--complete''. Let $x=\langle \name{x}_\alpha: \alpha<\gamma 
\rangle\conc\langle\bar{\bQ}\rangle\conc\langle({\cal I}_\alpha,
\name{f}^\alpha_p): \alpha<\gamma\ \&\ p\in\bP_\alpha\rangle$. 

\noindent Suppose that $N$ is $(\lambda,\family)$--good, $\bP_\gamma,x\in N$
and $\bar{p}$ is a $\bP_\gamma$--generic sequence over $N$. Note that
$\bar{\bQ}\in N$. We define a condition $r^*\in\bP_\gamma$: \quad we let
$\dom(r^*)=N\cap\gamma$ and we inductively define $r^*(\alpha)$ for
$\alpha\in\dom(r^*)$ by 
\begin{quotation}
\noindent if there is a $\bP_\alpha$--name $\name{\tau}$ such that 
\[r^*\rest\alpha\forces_{\bP_\alpha}\mbox{``}\name{\tau}\in\nbQ_\alpha\mbox{
is an upper bound to }\langle p_i(\alpha): i<\lambda\rangle\mbox{''}\]
then $r^*(\alpha)$ is the $<^*_\chi$--first such a name;

\noindent if there is no $\name{\tau}$ as above then
$r^*(\alpha)={\bf 0}_{\nbQ_\alpha}$.
\end{quotation}
It should be clear that $r^*\in\bP_\gamma$ (as $\|N\|=\lambda<\kappa$). What we
have to do is to show that $r^*$ is an upper bound to $\bar{p}$ in
$\bP_\gamma$. We do this showing by induction on $\alpha\leq\gamma$ that
\begin{description}
\item[($\otimes_\alpha$)] \quad for each $i<\lambda$,\ \  $p_i\rest\alpha
\leq_{\bP_\alpha} r^*\rest\alpha$.
\end{description}

For $\alpha=0$ there is nothing to do.

For $\alpha$ limit this is immediately by the induction hypothesis.

If $\alpha=\beta+1$ and $\beta\notin N$ then we use the induction hypothesis
and the fact that for each $i<\lambda$, $\dom(p_i)\subseteq\gamma\cap N$
(remember $p_i\in\bP_\gamma\cap N$, $\lambda\subseteq N$ and
$\|\dom(p_i)\|\leq\lambda$). 

So we are left with the case $\alpha=\beta+1$, $\beta\in N$. Suppose that
$G_\beta\subseteq\bP_\beta$ is a generic filter over $\V$ such that
$r^*\rest\beta\in G_\beta$ (so necessarily $p_i\rest\beta\in G_\beta$ for 
each $i<\lambda$). We will break the rest of the proof into several
claims. Each of them has a very standard proof, but we will sketch the proofs
for reader's convenience. Remember that we are in the case $\beta\in N$, so in
particular $\bP_\beta,\bP_{\beta+1},\name{x}_\beta,{\cal I}_\beta\in N$ and
$\langle p_i\rest\beta:i<\lambda\rangle\subseteq N$ is a $\bP_\beta$--generic
sequence over $N$.

\begin{claim}
\label{cl18}
\[r^*\rest\beta\forces_{\bP_\beta}\mbox{`` there is no new function from
$\lambda$ to $\V$ ''.}\]
\end{claim}

\noindent{\em Proof of the claim:}\hspace{0.15in} Since ${\cal I}_\beta\in N$
is an open dense subset of $\bP_\beta$ we know that $p_i\rest\beta\in{\cal
I}_\beta$ for some $i<\lambda$. If the condition $p_i\rest\beta$ forces that
``there is no new function from $\lambda$ to $\V$ then we are done (as
$r^*\rest\beta\geq p_i\rest\beta$). So suppose otherwise. Then $p_i\rest\beta
\forces_{\bP_\beta}$`` $\name{f}^\beta_{p_i\rest\beta}\notin\V$ ''. But, as
$\beta, p_i\rest\beta\in N$, we have $\name{f}^\beta_{p_i\rest\beta}\in N$ and
therefore for each $\zeta<\lambda$ there is $j<\lambda$ such that the
condition $p_j\rest\beta$ decides the value of $\name{f}^\beta_{p_i\rest
\beta}(\zeta)$. Consequently the condition $r^*\rest\beta$ decides all values
of $\name{f}^\beta_{p_i\rest\beta}$, so $r^*\rest\beta\forces_{\bP_\beta}
\name{f}^\beta_{p_i\rest\beta}\in\V$, a contradiction. 

\begin{claim}
\label{cl1}
$N[G_\beta]\cap\V=N$ (so $N[G_\beta]\cap\mu^*\in\family$).
\end{claim}

\noindent{\em Proof of the claim:}\hspace{0.15in} Suppose that $\name{\tau}\in
N$ is a $\bP_\beta$--name for an element of $\V$. As the sequence $\langle
p_i\rest\beta:i<\lambda\rangle$ is $\bP_\beta$--generic over $N$, for some
$i<\lambda$, the condition $p_i\rest\beta$ decides the value of the name
$\name{\tau}$. Since $p_i\rest\beta\in N$ the decision belongs to $N$
(remember the elementarity of $N$) and hence $\name{\tau}^{G_\beta}\in N$.

\begin{claim}
\label{cl2}
\[\|N[G_\beta]\|=\lambda,\quad N[G_\beta]^{<\lambda}\subseteq N[G_\beta]\quad
\mbox{ and }\quad N[G_\beta]\prec\Hchi^{\V[G_\beta]}.\]
Consequently,\quad $\V[G_\beta]\models$``the model $N[G_\beta]$ is $(\lambda,
\family)$--good and $\name{x}_\beta^{G_\beta}\in N[G_\beta]$''.
\end{claim}

\noindent{\em Proof of the claim:}\hspace{0.15in} Names for elements of
$N[G_\beta]$ are from $N$, so clearly $\|N[G_\beta]\|=\lambda=\|N\|$. It
follows from \ref{seqfromN} and \ref{cl18} that $N[G_\beta]^{<\lambda}
\subseteq N[G_\beta]$. To check that $N[G_\beta]$ is an elementary submodel of
$\Hchi$ (in $\V[G_\beta]$) we use the genericity of $\langle p_i\rest\beta:
i<\lambda\rangle$ and the elementarity of $N$: each existential formula of the
language of forcing (with parameters from $N$) is decided by some $p_i\rest
\beta$. If the decision is positive then there is in $N$ a name for a witness
for the formula. So we finish by Tarski--Vaught criterion. 

\begin{claim}
\label{cl3}
\[\V[G_\beta]\models\mbox{``}\langle p_i(\beta)^{G_\beta}\!:i<\lambda\rangle
\mbox{ is an increasing $\nbQ_\beta^{G_\beta}$--generic sequence over
}N[G_\beta]\mbox{''}.\] 
\end{claim}

\noindent{\em Proof of the claim:}\hspace{0.15in} By the induction hypothesis,
the condition $r^*\rest\beta$ is stronger then all $p_i\rest\beta$. Hence, as
$\bar{p}$ is increasing, the sequence $\langle p_i(\beta)^{G_\beta}:i<\lambda
\rangle$ is increasing (in $\nbQ_\beta^{G_\beta}$). Suppose now that
$\name{I}\in N$ is a $\bP_\beta$--name for an open dense subset of
$\nbQ_\beta$. Look at the set $\{p\in\bP_{\beta+1}: p\rest\beta\forces_{
\bP_\beta}p(\beta)\in\name{I}\}$. It is an open dense subset of
$\bP_{\beta+1}$ from $N$. But $\bP_{\beta+1}\lesdot\bP_\gamma$, so for some
$i<\lambda$ we have  
\[p_i\rest\beta\forces_{\bP_\beta} p_i(\beta)\in\name{I},\]
finishing the claim.
\medskip

By \ref{cl2}, \ref{cl3} (and the assumption that $\forces_{\bP_\beta}$
``$\nbQ_\beta$ is strongly $\family$--complete'') we conclude that, in
$\V[G_\beta]$, the sequence $\langle p_i(\beta)^{G_\beta}:i<\lambda\rangle
\subseteq\nbQ_\beta^{G_\beta}$ has an upper bound (in
$\nbQ_\beta^{G_\beta}$). Now, as $G_\beta$ was an arbitrary generic filter
containing $r^*\rest\beta$ we conclude that there is a $\bP_\beta$--name
$\name{\tau}$ such that 
\[r^*\rest\beta\forces_{\bP_\beta}\mbox{``}\name{\tau}\in\nbQ_\beta\mbox{
is an upper bound to }\langle p_i(\beta): i<\lambda\rangle\mbox{''}.\]
Now look at the definition of $r^*(\beta)$. 

\noindent 2)\ \ \ Left to the reader. \QED 

\begin{definition}
\label{basically}
Let (of course, $\kappa=\lambda^+$, and) $\family_0\in(\dkl)^+$ and $\hat{
\family}_1\in(\dlk[\family_0])^+$. Let $\bQ$ be a forcing notion and $\chi$ be
a large enough regular cardinal. 
\begin{enumerate}
\item We say that a sequence $\bar{N}=\langle N_i: i\leq\lambda\rangle$ is
{\em $(\lambda,\kappa,\hat{\family}_1,\bQ)$--considerable} if
\begin{quotation}
\noindent $\bar{N}$ is an increasing continuous sequence of elementary
submodels of $\Hchi$ such that $\lambda\cup\{\lambda,\kappa,\bQ\}\subseteq
N_0$, the sequence $\langle N_i\cap\mu^*: i\leq\lambda\rangle$ is in
$\hat{\family}_1$ and for each $i<\lambda$
\[\|N_i\|<\kappa\quad\mbox{ and }\quad \langle N_j: j\leq i\rangle\in N_{i+1}
\quad\mbox{ and}\]
\[i\mbox{ is non-limit }\quad\Rightarrow\quad (N_i)^{<\lambda}\subseteq N_i.\]
\end{quotation}
\item For a $(\lambda,\kappa,\hat{\family}_1,\bQ)$--considerable sequence
$\bar{N}=\langle N_i: i\leq \lambda\rangle$ and a condition $r\in N_0\cap\bQ$,
let $\innagra(\bQ,r)$ be the following game of two players, $\com$ and $\inc$:
\begin{quotation}
\noindent the game lasts $\lambda$ moves and during a play the players
construct a sequence $\langle (p_i,\bar{q}_i):i<\lambda\rangle$ such that each
$p_i$ is a condition from $\bQ$ and $\bar{q}_i=\langle q_{i,\xi}:
\xi<\lambda\rangle$ is an increasing $\lambda$--sequence of conditions from
$\bQ$ (we may identify it with its least upper bound in the completion
$\hat{\bQ}$) and at the stage $i<\lambda$ of the game: 

\noindent the player $\com$ chooses a condition $p_i\in N_{-1+i+1}\cap\bQ$
such that 
\[r\leq p_i,\qquad(\forall j<i)(\forall\xi<\lambda)(q_{j,\xi}\leq p_i),\]
and the player $\inc$ answers choosing a $\leq_{\bQ}$--increasing
$\bQ$--generic over $N_{-1+i+1}$ sequence $\bar{q}_i=\langle q_{i,\xi}:
\xi<\lambda\rangle\subseteq N_{-1+i+1}\cap \bQ$ such that
\[p_i\leq q_{i,0},\quad\mbox{ and }\quad\bar{q}_i\in N_{-1+i+2}.\]
\end{quotation}
The player $\com$ wins the play of $\innagra(\bQ,r)$ if the sequence $\langle
p_i: i<\lambda\rangle$ constructed by him during the play has an upper bound
in $\bQ$. 

\item We say that the forcing notion $\bQ$ is {\em basically $(\family_0,
\hat{\family}_1)$--complete} if 
\begin{description}
\item[($\alpha$)] $\bQ$ is $(<\lambda)$--complete (see \ref{complete}(3)), and
\item[($\beta$)]  $\bQ$ is strongly $\family_0$--complete (see \ref{good}(3)),
and 
\item[($\gamma$)] for most $(\lambda,\kappa,\hat{\family}_1,
\bQ)$--considerable sequences $\bar{N}=\langle N_i: i\leq\lambda\rangle$ and a
condition $r\in N_0\cap\bQ$ the player $\inc$ DOES NOT have a winning strategy
in the game $\innagra(\bQ,r)$. 
\end{description}
\end{enumerate}
\end{definition}

\begin{remark}
\label{mogagrac}
{\em 
1)\quad Why do we have ``strongly $\family_0$--complete'' in
\ref{basically}(3)($\beta$) and not ``strategically $\family_0$--complete''?
To help proving the preservation theorem. 

\noindent 2)\quad Note that if a forcing notion $\bQ$ is strongly
$\family_0$--complete and $(<\lambda)$--complete, and $\bar{N}$ is
$(\lambda,\kappa,\hat{\family}_1, \bQ)$--considerable then both players have
always legal moves in the game $\innagra(\bQ,r)$. Moreover, if $\bQ$ is a
dense suborder of $\bQ'$ and $\bQ'\in N_0$ and the player $\com$ plays
elements of $\bQ$ only then both players have legal moves in the game
$\innagra(\bQ',r)$.\\ 
\relax [Why? Arriving at stage $i$ of the game, the player $\com$ has to
choose a condition $p_i\in N_{-1+i+1}\cap\bQ$ stronger than all $q_{j,\xi}$
(for $j<i$, $\xi<\lambda$). If $i$ is a limit ordinal, $\com$ looks at the
sequence $\langle p_j: j<i\rangle$ constructed by him so far. Since
$(N_{i+1})^{<\lambda}\subseteq N_{i+1}$ we have that $\langle p_j: j<i\rangle
\in N_{i+1}$ and, as $\bQ$ is $(<\lambda)$--complete, this sequence has an
upper bound in $N_{i+1}$ (remember that $N_{i+1}$ is an elementary submodel of
$\Hchi$). This upper bound is good for $q_{j,\xi}$ ($j<i$, $\xi<\lambda$) too.
If $i=i_0+1$ then the player $\com$ looks at the sequence $\bar{q}_{i_0}\in
N_{-1+i_0+2}$ only. It is $\bQ$-generic over $N_{-1+i_0+1}$, $\bQ$ is strongly
$\family_0$--complete and $N_{-1+i_0+1}$ is $(\lambda,\family_0)$--good. 
Therefore, there is an upper bound to $\bar{q}_{i_0}$, and by the elementarity
there is one in $N_{-1+i_0+2}$. Now, the player $\inc$ may always use the fact
that $\bQ$ is $(<\lambda)$--complete to build above $p_i$ an increasing
sequence $\bar{q}_i\subseteq\bQ\cap N_{-1+i+1}$ which is generic over
$N_{-1+i+1}$. Since $N_{-1+i+1}\in N_{-1+i+2}$, by the elementarity there are
such sequences in $N_{-1+i+2}$.\\
Concerning the ``moreover'' part note that the only difference is when $\com$
is supposed to choose an upper bound to $\bar{q}_{i_0}$. But then he proceeds
like in \ref{strondense} reducing the task to finding an upper bound to a
sequence (generic over $N_{-1+i_0+1}$) of elements of $\bQ$.]
}
\end{remark}

Unfortunately, the amount of completeness demanded in \ref{basically} is too
large to capture the examples we have in mind (see the next section). 
Therefore we slightly weaken the demand \ref{basically}(3)($\gamma$) (or
rather we change the appropriate game a little). In definition \ref{really}
below we formulate the variant of completeness which seems to be the right one
for our case.  

\begin{definition}
\label{really}
Let $\family_0\in(\dkl)^+$ and $\hat{\family}_1\in (\dlk[\family_0])^+$. Let
$D$ be a function such that $\dom(D)=\hat{\family}_1$ and for every $\bar{a}
\in\hat{\family}_1$
\[D(\bar{a})=D_{\bar{a}}\ \mbox{ is a filter on }\lambda.\]
Let $\bQ$ be a forcing notion. 
\begin{enumerate}
\item We say that an increasing continuous sequence $\bar{N}=\langle N_i:i
\leq\lambda\rangle$ of elementary submodels of $\Hchi$ is {\em $(\lambda,
\kappa,\hat{\family}_1,D,\bQ)$--suitable} if\\
$\lambda\cup\{\lambda,\kappa,\bQ\}\subseteq N_0$, $\|N_i\|<\kappa$, $\langle
N_j: j\leq i\rangle\in N_{i+1}$ and there are $\bar{a}\in\hat{\family}_1$ and
$X\in D_{\bar{a}}$ such that for each $i\in X$
\[(N_{i+1})^{<\lambda}\subseteq N_{i+1}\quad\&\quad N_{i+1}\cap\mu^*=a_{i+1}\]
(compare with \ref{basically}(1)); we can add $N_i\cap\mu^*=a_i$ if
$D_{\bar{a}}$ is normal. 

A pair $(\bar{a},X)$ witnessing the last demand on $\bar{N}$ will be called
{\em suitable base for $\bar{N}$}.
\item For a $(\lambda,\kappa,\hat{\family}_1,D,\bQ)$--suitable sequence 
$\bar{N}=\langle N_i: i\leq \lambda\rangle$, a suitable base $(\bar{a},X)$ for
$\bar{N}$ and a condition $r\in N_0$, let $\tagra(\bQ,r)$ be the following
game of two players, $\com$ and $\inc$: 
\begin{quotation}
\noindent The game lasts $\lambda$ moves and during a play the players
construct a sequence $\langle (p_i,\zeta_i,\bar{q}_i):i<\lambda\rangle$ such
that $\zeta_i\in X$, $p_i\in\bQ$ and $\bar{q}_i=\langle q_{i,\xi}:\xi<\lambda
\rangle\subseteq\bQ$ in the following manner.\\
At the stage $i<\lambda$ of the game: 

\noindent player $\com$ chooses $\zeta_i\in X$ above all $\zeta_j$ chosen so
far and then he picks a condition $p_i\in N_{\zeta_i+1}\cap\bQ$ such that 
\[r\leq p_i,\qquad(\forall j<i)(\forall\xi<\lambda)(q_{j,\xi}\leq p_i),\]
after this player $\inc$ answers choosing a $\leq_{\bQ}$--increasing
$\bQ$--generic over $N_{\zeta_i+1}$ sequence $\bar{q}_i=\langle q_{i,\xi}:
\xi<\lambda\rangle\subseteq N_{\zeta_i+1}\cap \bQ$ such that
\[p_i\leq q_{i,0},\quad\mbox{ and }\quad\bar{q}_i\in N_{\zeta_i+2}.\]
\end{quotation}
The player $\com$ wins the play of $\tagra(\bQ,r)$ if $\{\zeta_i:i<\lambda\}
\in D_{\bar{a}}$ and the sequence $\langle p_i: i<\lambda\rangle$ constructed
by him during the play has an upper bound in $\bQ$. 

\item We say that the forcing notion $\bQ$ is {\em really $(\family_0,
\hat{\family}_1,D)$--complete} if
\begin{description}
\item[($\alpha$)] $\bQ$ is $(<\lambda)$--complete (see \ref{complete}(3)), and
\item[($\beta$)]  $\bQ$ is strongly $\family_0$--complete (see \ref{good}(3)),
and 
\item[($\gamma$)] for most $(\lambda,\kappa,\hat{\family}_1,D,
\bQ)$--suitable sequences $\bar{N}=\langle N_i: i\leq\lambda\rangle$, a
suitable basis $(\bar{a},X)$ for $\bar{N}$ and a condition $r\in N_0\cap\bQ$,
the player $\inc$ DOES NOT have a winning strategy in the game
$\tagra(\bQ,r)$.  
\end{description}
\end{enumerate}
\end{definition}

\begin{remark}
{\em 
If a forcing notion $\bQ$ is strongly $\family_0$--complete and
$(<\lambda)$--complete, and $\bar{N}$ is  $(\lambda,\kappa,\hat{\family}_1,D,
\bQ)$--suitable (witnessed by $(\bar{a},X)$) then both players have always
legal moves in the game $\tagra(\bQ,r)$. Moreover, if $\bQ$ is dense in
$\bQ'$, $\bQ'\in N_0$ and $\com$ plays elements of $\bQ$ only then both
players have legal moves in $\tagra(\bQ',r)$\\
\relax [Why? Like in \ref{mogagrac}.]
}
\end{remark}

\begin{remark}
\label{renotation}
{\em
We may equivalently describe the game $\tagra(\bQ,r)$ in the following
manner. Let $\hat{\bQ}$ be the completion of $\bQ$.
\begin{quotation}
\noindent The play lasts $\lambda$ moves during which the players construct a
sequence $\langle p_i,q_i:i<\lambda\rangle$ such that $p_i\in N_{i+1}\cap
(\bQ\cup\{*\})$ (where $*\notin\bQ$ is a fixed element of $N_0$), $q_i\in
N_{i+2}\cap\hat{\bQ}$ .\\
At the stage $i<\lambda$ of the game, $\com$ chooses $p_i$ in such a way that
\[p_i\neq *\quad\Rightarrow\quad i\in X\ \&\ (\forall j<i)(q_j\leq_{\hat{\bQ}}
p_i),\] 
and $\inc$ answers choosing $q_i$ such that

if $p_i=*$ then $q_i$ is the least upper bound of $\langle q_j: j<i\rangle$ in
$\hat{\bQ}$,  

if $p_i\neq *$ then $q_i\in N_{i+2}\cap\hat{\bQ}$ is the least upper bound
of a $\bQ$--generic filter over $N_{i+1}$ containing $p_i$. 
\end{quotation}
The player $\com$ wins if $\{i<\lambda: p_i\neq *\}\in D_{\bar{a}}$ and the
sequence $\langle p_i: p_i\neq *\rangle$ has an upper bound.

There is no real difference between \ref{really}(2) and the description given
above. Here, instead of ``jumping'' player $\com$ puts $*$ (which has the
meaning of {\em I am waiting}) and he uses the existence of the least upper
bounds to replace a generic sequence by its least upper bound.
}
\end{remark}

\begin{proposition}
\label{reallydense}
Suppose that $\family_0\in(\dkl)^+$, $\hat{\family}_1\in (\dlk[\family_0])^+$
and $D_{\bar{a}}$  is a filter on $\lambda$ for $\bar{a}\in \hat{\family}_1$. 
Assume that $\bQ$ is a dense suborder of $\bQ'$, $\bar{N}$ is $(\lambda,
\kappa,\hat{\family}_1,D,\bQ)$--suitable (witnessed by $(\bar{a},X)$),
$\bQ'\in N_0$. Then for each $r\in\bQ$:
\begin{quotation}
\noindent the player $\com$ has a winning strategy in $\tagra(\bQ,r)$ (the
player $\inc$ does not have a winning strategy in $\tagra(\bQ,r)$,
respectively)  

if and only if 

\noindent he has a winning strategy in $\tagra(\bQ',r)$ (the player $\inc$
does not have a winning strategy in $\tagra(\bQ',r)$, resp.). 
\end{quotation}
\end{proposition}

\Proof Suppose that $\com$ has a winning strategy in $\tagra(\bQ,r)$. We 
describe a winning strategy for him in $\tagra(\bQ',r)$ which tells him to
play elements of $\bQ$ only. The strategy is very simple. At each stage
$i<\lambda$, $\com$ replaces the sequence $\bar{q}_i\subseteq\bQ'$ by a
sequence $\bar{q}^*_i\subseteq\bQ$ which has the same upper bounds in $\bQ$ as
$\bar{q}_i$, is increasing and generic over $N_{\zeta_i+1}$. To do this he
applies the procedure from the proof of \ref{strondense} (in $N_{\zeta_i+2}$,
of course). Then he may use his strategy from $\tagra(\bQ,r)$.
The converse implication is easy too: if the winning strategy of $\com$ in
$\tagra(\bQ',r)$ tells him to play $\zeta_i,p_i$ then he puts $\zeta_i$ and any
element $p^*_i$ of $\bQ\cap N_{\zeta_i+1}$ stronger than $p_i$. Note that this
might be interpreted as playing $p_i$ followed by a sequence $p^*_i\conc
\bar{q}_i$. \QED  

\begin{proposition}
Suppose $\family_0\in(\dkl)^+$ and $\hat{\family}_1\in (\dlk[\family_0])^+$
(and as usual in this section, $\kappa=\lambda^+$). Let $D_{\bar{a}}$ be the
club filter of $\lambda$ for each $\bar{a}\in\hat{\family}_1$. Then any
really $(\family_0,\hat{\family}_1,D)$--complete forcing notion preserves
stationarity of $\family_0$, $\hat{\family}_1$ in the respective filters. \QED
\end{proposition}

\section{Examples}
\label{sekcjaexamples}
Before we continue with the general theory let us present a simple example
with the properties we are investigating. Of course it is related to {\em
guessing clubs}; remember that there are ZFC theorems saying that many times
we can guess clubs (see \cite{Sh:g}). 

\begin{hypothesis}
Assume $\lambda^{<\lambda}=\lambda$ and $\lambda^+=\kappa$. Suppose that
$\family_0=S_0\subseteq S^\kappa_\lambda$ is a stationary set such that $S
\stackrel{\rm def}{=}S^\kappa_\lambda\setminus S_0$ is stationary too (but the
definitions below are meaningful also when $S=\emptyset$). Let  
\[\begin{array}{r}
\hat{\family}_1=\big\{\bar{a}=\langle a_i: i\leq\lambda\rangle: \bar{a}\
\mbox{ is increasing continuous and for each }i\leq\lambda,\ \ \\ 
a_i\in \kappa\ \mbox{ and if $i$ is not limit then }a_i\in S_0\big\}.
\end{array}\]
[Check that $\family_0\in ({\frak D}_{<\kappa,<\lambda}(\kappa))^+$ and
$\hat{\family}_1\in {\frak D}_{<\kappa,<\lambda}^\lambda[\family_0]$.]
\end{hypothesis}
Note that (provably in ZFC, see \cite[Ch III, \S 2]{Sh:g}) there is a sequence
$\bar{C}=\langle C_\delta:\delta\in S\rangle$ such that for each $\delta\in S$
\begin{quotation}
\noindent $C_\delta$ is a club of $\delta$ of the order type $\lambda$ and

\noindent if $\alpha\in\nacc(C_\delta)$ then $\cf(\alpha)=\lambda$ 
\end{quotation}
such that $\emptyset\notin {\rm id}^p(\bar{C})$, i.e. for every club $E$ of
$\kappa$ for stationary many $\delta\in S$, $\delta=\sup(E\cap\nacc(C_\delta)
)$, even $\{\alpha<\delta:\min(C_\delta\setminus(\alpha+1))\in E\}$ is a
stationary subset of $\delta$. We can use this to show that some natural
preservation of not adding bounded subsets of $\kappa$ (or just not collapsing
cardinals) necessarily fails, just considering the forcing notion killing the
property of such $\bar{C}$. [Why? As in the result such $\bar{C}$ exists, but
by iterating we could have dealt with all possible $\bar{C}$'s.] We will show
that we cannot demand   
\[\alpha\in\nacc(C_\delta)\quad\Rightarrow\quad\cf(\alpha)<\lambda,\]
that is in some forcing extension preserving GCH there is no such $\bar{C}$. 
So, for $\bar{C}$ as earlier but with the above demand we want to add
generically a club $E$ of $\lambda^+$ such that  
\[(\forall\delta\in S)(E\cap\nacc(C_\delta)\mbox{ is bounded in }\delta).\]
We will want our forcing to be quite complete. To get the consistency of no
guessing clubs we need to iterate, what {\em is} our main theme. 

\begin{definition}
\label{qjeden}
Let $\bar{C}=\langle C_\delta: \delta\in S\rangle$ be a sequence such that
for every $\delta\in S$: 
\begin{quotation}
\noindent $C_\delta$ is a club of $\delta$ of the order type $\lambda$ and

\noindent if $\alpha\in\nacc(C_\delta)$ then $\cf(\alpha)<\lambda$ (or at
least $\alpha\notin S_0$).
\end{quotation}
We define a forcing notion $\bQ_{\bar{C}}^1$ to add a desired club $E\subseteq
\lambda^+$:

\noindent{\bf a condition} in $\bQ_{\bar{C}}^1$ is a closed subset $e$ of
$\lambda^+$ such that $\alpha_e\stackrel{\rm def}{=}\sup(e)<\lambda^+$ and
\[(\forall \delta\in S\cap (\alpha_e+1))(e\cap\nacc(C_\delta)\mbox{ is
bounded in }\delta),\]

\noindent{\bf the order} $\leq_{\bQ_{\bar{C}}^1}$ of $\bQ_{\bar{C}}^1$ is
defined by

$e_0\leq_{\bQ_{\bar{C}}^1} e_1$ \quad if and only if\quad $e_0$ is an initial
segment of $e_1$. 
\end{definition}
It should be clear that $(\bQ_{\bar{C}}^1,\leq_{\bQ_{\bar{C}}^1})$ is a
partial order. We claim that is quite complete.

\begin{proposition}
\label{stepAB}
\begin{enumerate}
\item $\bQ_{\bar{C}}^1$ is $(<\lambda)$--complete.
\item $\bQ_{\bar{C}}^1$ is strongly $\family_0$--complete.
\end{enumerate}
\end{proposition}

\Proof 1)\ \ \ Should be clear.

\noindent 2)\ \ \ Suppose that $N\prec\Hchi$ is $(\lambda,\family_0)$--good
(see \ref{good}) and $\bQ_{\bar{C}}^1\in N$. Further suppose that $\bar{e}
=\langle e_i: i<\lambda\rangle\subseteq\bQ_{\bar{C}}^1\cap N$ is an increasing 
$\bQ_{\bar{C}}^1$--generic sequence over $N$. Let $e\stackrel{\rm def}{=}
\bigcup\limits_{i<\lambda}e_i\cup\{\sup(\bigcup\limits_{i<\lambda}e_i)\}$.

\begin{claim}
\label{cl4}
$e\in\bQ_{\bar{C}}^1$.
\end{claim}

\noindent{\em Proof of the claim:}\hspace{0.15in} First note that as each
$e_i$ is the end extension of all $e_j$ for $j<i$, the set $e$ is
closed. Clearly $\alpha_e\stackrel{\rm def}{=}\sup(e)<\lambda^+$ (as each
$\alpha_{e_i}$ is below $\lambda^+$). So what we have to check is that
\[(\forall \delta\in S\cap (\alpha_e+1))(e\cap\nacc(C_\delta)\mbox{ is
bounded in }\delta).\]
Suppose that $\delta\in S\cap (\alpha_e+1)$. If $\delta<\alpha_e$ then for
some $i<\lambda$ we have $\delta\leq\alpha_{e_i}$ and $e\cap\delta=e_i\cap
\delta$ and therefore $e\cap\nacc(C_\delta)$ is bounded in $\delta$. So a
problem could occur only if $\delta=\alpha_e=\sup\limits_{i<\lambda}
\alpha_{e_i}$, but we claim that it is impossible. Why? Let $\delta^*=N\cap
\lambda^+$, so $\delta^*\in\family_0$ (as $N$ is $(\lambda,\family_0)$--good)
and therefore $\delta^*\neq\delta$ (as $S_0\cap S=\emptyset$). For each $\beta
<\delta^*$ the set   
\[{\cal I}_\beta\stackrel{\rm def}{=}\{q\in\bQ_{\bar{C}}^1: q\setminus\beta
\neq\emptyset\}\]
is open dense in $\bQ_{\bar{C}}^1$ (note that if $q\in\bQ_{\bar{C}}^1$, $q
\setminus\beta=\emptyset$ then $q\leq q\cup\{\alpha_q,\beta+1\}\in
\bQ_{\bar{C}}^1$). Clearly ${\cal I}_\beta\in N$. Consequently, by the
genericity of $\bar{e}$, $e_i\in {\cal I}_\beta$ for some $i<\lambda$ and thus
$\alpha_{e_i}>\beta$. Hence $\sup\limits_{i<\lambda}\alpha_{e_i}\geq
\delta^*$. On the other hand, as each $e_i$ is in $N$ we have $\alpha_{e_i}<
\delta^*$ (for each $i<\lambda$) and hence $\delta^*=\sup\limits_{i<\lambda}
\alpha_{e_i}=\delta$, a contradiction. 

\begin{claim}
\label{cl5}
For each $i<\lambda$, $e_i\leq e$.
\end{claim}

\noindent{\em Proof of the claim:}\hspace{0.15in} Should be clear.
\medskip

Now, by \ref{cl4}+\ref{cl5}, we are done. \QED   

\begin{proposition}
\label{getreally}
For each $\bar{a}\in\hat{\family}_1$, let $D_{\bar{a}}$ be the club filter of
$\lambda$ (or any normal filter on $\lambda$). Then the forcing notion
$\bQ_{\bar{C}}^1$ is really $(\family_0, \hat{\family}_1,D)$--complete.
\end{proposition}

\Proof By \ref{stepAB} we have to check demand \ref{really}(3$\gamma$)
only. So suppose that $\bar{N}=\langle N_i: i\leq\lambda\rangle$ is
$(\lambda,\kappa,\hat{\family}_1,D,\bQ_{\bar{C}}^1)$--suitable and $(\bar{a},
X)$ is a suitable basis for $\bar{N}$ (and we may think that $X$ is a closed
unbounded subset of $\lambda$). Let $r\in N_0$. We are going to describe a
winning strategy for player $\com$ in the game $\tagra(\bQ_{\bar{C}}^1,r)$. 
There are two cases to consider here: $N_\lambda\cap\kappa\in S$ and
$N_\lambda\cap\kappa\notin S$. The winning strategy for $\com$ in $\tagra(
\bQ_{\bar{C}}^1,r)$ is slightly more complicated in the first case, so
let us describe it only then. So we assume $N_\lambda\cap\kappa\in S$.

Arriving at the stage $i<\lambda$ of the game, $\com$ chooses $\zeta_i$
according to the following rules:
\begin{description}
\item[if $i=0$] \qquad\quad then he takes $\zeta_i=\min X$,
\item[if $i=i_0+1$] \quad  then he takes
\[\zeta_i=\min\big\{j\in X: \zeta_{i_0}+1<j\quad\&\quad (N_{\zeta_{i_0}+1}
\cap\kappa, N_j\cap\kappa)\cap C_{N_\lambda\cap\kappa}\neq\emptyset\big\},\]
\item[if $i$ is limit]\ \ \ then he lets $\zeta_i=\sup\limits_{j<i}\zeta_j$.
\end{description}
Note that as $C_{N_\lambda\cap\kappa}$ is unbounded in $N_\lambda\cap\kappa$
and $X$ is a club of $\lambda$, the above definition is correct; i.e. the
respective $\zeta_i$ exists, belongs to $X$ and is necessarily above all
$\zeta_j$ chosen so far. Next $\com$ plays $p_i$ defined as follows. The first 
$p_0$ is just $r$. If $i>0$ then $\com$ takes the first ordinal $\gamma_i$
such that 
\[\sup(N_{\zeta_i+1}\cap C_{N_\lambda\cap\kappa})<\gamma_i<N_{\zeta_i+1}\cap
\kappa\]
and he  puts 
\[p_i=\bigcup_{\scriptstyle j<i\atop \scriptstyle \xi<\lambda} q_{j,\xi} \cup
\{N_{\bigcup\limits_{j<i}\zeta_j+1}\cap\kappa\}\cup\{\gamma_i\}.\]
Note that $\cf(N_{\zeta_i+1}\cap\kappa)=\lambda$ so $C_{N_\lambda\cap\kappa}
\cap N_{\zeta_i+1}\cap\kappa$ is bounded in $N_{\zeta_i+1}\cap\kappa$ and the
$\gamma_i$ above is well defined. Moreover, by arguments similar to that of
\ref{stepAB}, one easily checks that $\bigcup\limits_{\scriptstyle j<i\atop
\scriptstyle \xi<\lambda} q_{j,\xi}\in \bQ_{\bar{C}}^1$ and then easily
$p_i\in\bQ_{\bar{C}}^1$ and it is $\leq_{\bQ_{\bar{C}}^1}$--stronger then all
$q_{j,\xi}$ (for $j<i$, $\xi<\lambda$). Consequently, the procedure described
above produces a legal strategy for $\com$ in $\tagra(\bQ_{\bar{C}}^1,r)$. But
why is this a winning strategy for $\com$? Suppose that $\langle (p_i,\zeta_i,
\bar{q}_i):i<\lambda\rangle$ is the result of a play in which $\com$ follows
our strategy. First note that the sequence $\langle\zeta_i:i<\lambda\rangle$
is increasing continuous so it contains a club of $\lambda$ and thus $\{
\zeta_i:i<\lambda\}\in D_{\bar{a}}$. Now, let $e=\bigcup\limits_{i<\lambda}
p_i\cup\{N\cap\kappa^+\}$. We claim that $e\in\bQ_{\bar{C}}^1$. First note
that it is a closed subset of $\lambda^+$ with $\sup e\stackrel{\rm def}{=}
\alpha_e=N_\lambda\cap\kappa$. So suppose now that $\delta\in S\cap(\alpha_e
+1)$. If $\delta<\alpha_e$ then necessarily $\delta<\alpha_{p_i}$ for some
$i<\lambda$ and therefore $e\cap \nacc(C_\delta)=p_i\cap\nacc(C_\delta)$ is
bounded in $\delta$. The only danger may come from $\delta=N_\lambda\cap
\kappa$. Thus assume that $\beta\in e$ and we ask where does $\beta$ come
from? If it is from $p_0\cup\bigcup\limits_{\xi<\lambda}q_{0,\xi}$ then we
cannot say anything about it (this is the part of $e$ that we do not
control). But in all other instances we may show that $\beta\notin\nacc(
C_{N_\lambda\cap\kappa})$.  Why? If $\beta\in\bigcup\limits_{\xi<\lambda}
q_{i,\xi}\setminus p_i$ for some $0<i<\lambda$, then by the choice of
$\gamma_i$ and $p_i$ and the demand that $\bar{q}_i\subseteq N_{\zeta_i+1}$ we
have that $\beta\notin C_{N_\lambda\cap\kappa}$. Similarly if $\beta=
\gamma_i$. So the only possibility left is that $\beta=N_{\bigcup\limits_{j<
i}\zeta_j+1}\cap\kappa$. If $i$ is not limit then $\cf(N_{\bigcup\limits_{j<
i}\zeta_j+1}\cap\kappa)=\lambda$ so $\beta\notin\nacc(C_{N_\gamma\cap
\kappa})$. If $i$ is limit then, by the choice of the $\zeta_j$'s we have
$N_{\bigcup\limits_{j<i}\zeta_j+1}\cap\kappa\in\acc(C_{N_\gamma\cap\kappa})$
and we are clearly done. 
\medskip

Note that if $N_\lambda\cap\kappa\notin S$ then the winning strategy for $\com$
is much simpler: choose successive elements of $X$ as the $\zeta_j$'s and 
play natural bounds to sequences constructed so far. \QED

\begin{remark}
{\em
1)\ \ \ Note that one cannot prove that the forcing notion $\bQ_{\bar{C}}^1$
is basically $(\family_0,\hat{\family}_1)$--complete. The place in which a try
to repeat the proof of \ref{getreally} fails is the limit case of
$N_i\cap\kappa$. If we do not allow $\com$ to make ``jumps'' (the choices of
$\zeta_i$) then he cannot overcome difficulties coming from the case
exemplified by 
\[C_{N_\lambda\cap\kappa}=\{N_{\omega\cdot i}\cap\kappa: i<\lambda\}.\]
2)\ \ \ The instance $S=S^{\lambda^+}_\lambda$ is not covered here, but we
will deal with it later.
}
\end{remark}

The following forcing notion is used to get ${\rm Pr}_S$ (see
\ref{uniformization}). 

\begin{definition}
\label{defqdwa}
Let $\bar{C}=\langle C_\delta: \delta\in S\rangle$ be as in \ref{qjeden} and
let $\bar{h}=\langle h_\delta:\delta\in S\rangle$ be a sequence such that
$h_\delta:C_\delta\longrightarrow\lambda$ for $\delta\in S$. Further let
$\bar{D}=\langle D_\delta: \delta\in S\rangle$ be such that each $D_\delta$ is
a filter on $C_\delta$. 
\begin{enumerate}
\item We define a forcing notion $\bQ^2_{\bar{C},\bar{h}}$:

\noindent{\bf a condition} in $\bQ_{\bar{C},\bar{h}}^2$ is a function
$f:\alpha_f\longrightarrow\lambda$ such that $\alpha_f<\lambda^+$ and
\[(\forall \delta\in S\cap (\alpha_f+1))(\{\beta\in C_\delta: h_\delta(\beta)
=f(\beta)\}\mbox{ is a co-bounded subset of }C_\delta),\]

\noindent{\bf the order} $\leq_{\bQ_{\bar{C},\bar{h}}^2}$ of $\bQ_{\bar{C},
\bar{h}}^2$ is the inclusion (extension). 
\item The forcing notion $\bQ^{2,\bar{D}}_{\bar{C},\bar{h}}$ is defined
similarly, except that we demand that a condition $f$ satisfies
\[(\forall \delta\in S\cap (\alpha_f+1))(\{\beta\in C_\delta: h_\delta(
\beta)=f(\beta)\}\in D_\delta).\]
\end{enumerate}
\end{definition}

\begin{proposition}
\label{basqcom}
Let $D_{\bar{a}}$ be a club filter of $\lambda$ for $\bar{a}\in
\hat{\family}_1$. Then the forcing notion $\bQ_{\bar{C},\bar{h}}^2$ is really
$(\family_0,\hat{\family}_1,D)$--complete.
\end{proposition}

\Proof This is parallel to \ref{getreally}. It should be clear that
$\bQ^2_{\bar{C},\bar{h}}$ is $(<\lambda)$--complete. The proof that it is
strongly $\family_0$--complete goes like that of \ref{stepAB}(2), so what we
need is the following claim.

\begin{claim}
\label{cl19}
For each $\beta<\lambda^+$ the set 
\[{\cal I}_\beta\stackrel{\rm def}{=}\{f\in\bQ^2_{\bar{C},\bar{h}}:\beta\in
\dom(f)\}\]
is open dense in $\bQ^2_{\bar{C},\bar{h}}$.
\end{claim}

\noindent{\em Proof of the claim:}\hspace{0.15in} Let $f\in\bQ^2_{\bar{C},
\bar{h}}$. We have to show that for each $\delta<\lambda^+$ there is a
condition $f'\in\bQ^2_{\bar{C},\bar{h}}$ such that $f\leq f'$ and $\delta\leq
\alpha_{f'}$. Assume that for some $\delta<\lambda^+$ there is no suitable
$f'\geq f$, and let $\delta$ be the first such ordinal (necessarily $\delta$
is limit). Choose an increasing continuous sequence $\langle\beta_\zeta:\zeta
<\cf(\delta)\rangle$ cofinal in $\delta$ and such that $\beta_0=\alpha_f$ and
$\beta_\zeta\in\delta\setminus S$ for $0<\zeta<\cf(\delta)$. For each $\zeta<
\cf(\delta)$ pick a condition $f_\zeta\geq f$ such that $\alpha_{f_\zeta}=
\beta_\zeta$ and let $f^*=\bigcup\limits_{\zeta<\cf(\beta)}f_{\zeta+1}\rest
[\beta_\zeta,\beta_{\zeta+1})$. If $\delta\notin S$ then easily $f^*\in
\bQ^2_{\bar{C},\bar{h}}$ is a condition stronger than $f$. Otherwise we take
$f':\delta\longrightarrow\lambda$ defined by
\[f'(\xi)=\left\{
\begin{array}{ll}
h_\delta(\xi)&\mbox{if }\xi\in C_\delta\setminus\alpha_f,\\
f^*(\xi)     &\mbox{otherwise.}
\end{array}
\right.\]
Plainly, $f'\in\bQ^2_{\bar{C},\bar{h}}$ and it is stronger than $f$. Thus in
both cases we may construct a condition $f'$ stronger than $f$ and such that
$\delta=\alpha_{f'}$, a contradiction.
\medskip

With \ref{cl19} in hands we may repeat the proof of \ref{stepAB}(2) with no
substantial changes. 

The proof that $\bQ^2_{\bar{C},\bar{h}}$ is really $(\family_0,\hat{\family}_1,
D)$--complete is similar to that of \ref{getreally}. So let $\bar{N}$,
$(\bar{a},x)$ and $r$ be as there and suppose that $N_\lambda\cap\kappa\in
S$. The winning strategy for $\com$ says him to choose $\xi_i$ as in
\ref{getreally} and play $p_i$ defined as follows. The first $p_0$ is $r$. If
$i>0$ then $\com$ lets $p_i^\prime=\bigcup\limits_{\scriptstyle j<i\atop
\scriptstyle \xi<\lambda} q_{j,\xi}$ (which clearly is a condition in
$\bQ^2_{\bar{C},\bar{h}}$) and chooses $p_i\in\bQ^2_{\bar{C},\bar{h}}\cap
N_{\zeta_i+1}$ such that 
\[\begin{array}{l}
p^\prime_i\leq p_i,\quad C_{N_\lambda\cap\kappa}\cap N_{\zeta_i+1} \subseteq
\dom(p_i)\quad\mbox{ and}\\
(\forall \beta\in C_{N_\lambda\cap\kappa}\cap N_{\zeta_i+1})(\alpha_{
p^\prime_i}<\beta\ \Rightarrow\ p_i(\beta)=h_{N_\lambda\cap\kappa}(\beta)).
  \end{array}\]
Clearly this is a winning strategy for $\com$. \QED

\begin{remark}
{\em
\begin{enumerate}
\item The proof of \ref{basqcom} shows that $\bQ^2_{\bar{C},\bar{h}}$ is
actually basically $(\family_0,\hat{\family}_1)$--complete. The same applies
to \ref{A210}.
\item In \ref{defqdwa}, \ref{basqcom} we may consider $\bar{h}$ such that for
some $h^*:\kappa\longrightarrow\kappa$, for each $\delta\in S$ we have 
\[(\forall\alpha\in C_\delta)(h_\delta(\alpha)<h^*(\alpha)),\]
no real change here.
\item Why we need above $h^*$ at all? If we allow e.g.\ $h_\delta$ to be
constantly $\delta$ then clearly there is no function $f$ with domain $\kappa$
and such that $(\forall \delta\in S)(\delta>\sup\{\alpha\in C_\delta:
f(\alpha)\neq h_\delta(\alpha)\})$ (by Fodor lemma). We may still ask if we
could just demand $h_\delta:C_\delta\longrightarrow\delta$? Even this
necessarily fails, as we may let $h_\delta(\alpha)=\min(C_\delta\setminus
(\alpha+1))$. Then, if $f$ is as above, the set $E=\{\delta<\kappa: \delta$ is
a limit ordinal and $(\forall\alpha<\delta)(f(\alpha)<\delta)\}$ is a club of
$\kappa$. Hence for some $\delta\in S$ we have:\ \ \
$\lambda^2<\delta=\sup(E\cap \delta)=\otp(E\cap\delta)$ and we get an easy
contradiction. 
\end{enumerate}
}
\end{remark}

Another example of forcing notions which we have in mind when developing the
general theory is related to the following problem. Let $K$ be a
$\lambda$--free Abelian group of cardinality $\kappa$. We want to make it
a Whitehead group.

\begin{definition}
\label{group}
Suppose that
\begin{description}
\item[(a)] $K_1$ is a strongly $\kappa$--free Abelian group of cardinality
$\kappa$, $\langle K_{1,\alpha}:\alpha<\kappa\rangle$ is a filtration of $K_1$
(i.e. it is an increasing continuous sequence of subgroups of $K_1$ such that
$K_1= \bigcup\limits_{\alpha<\kappa} K_{1,\alpha}$ and each $K_{1,\alpha}$ is
of size $<\kappa$), 
\[\Gamma=\{\alpha<\kappa: K_1/K_{1,\alpha}\mbox{ is not
$\lambda$--free}\;\},\]
\item[(b)] $K_2$ is an Abelian group extending ${\Bbb Z}$,
$h:K_2\stackrel{\rm onto}{\longrightarrow} K_1$ is a homomorphism with kernel
${\Bbb Z}$. 
\end{description}
We define a forcing notion $\bQ^3_{K_2,h}$:

\noindent{\bf a condition} in $\bQ_{K_2,h}^3$ is a homomorphism $g:
K_{1,\alpha}\longrightarrow K_2$ such that $\alpha\in\kappa\setminus\Gamma$
and $h\comp g={\rm id}_{K_1,\alpha}$,

\noindent{\bf the order} $\leq_{\bQ_{K_2,h}^3}$ of $\bQ_{K_2,h}^3$ is the
inclusion (extension).  
\end{definition}

\begin{proposition}
\label{A210}
Let $D_{\bar{a}}$ be a club filter of $\lambda$ for $\bar{a}\in
\hat{\family}_1$. Assume $K_1,K_{1,\alpha},K_2$ and $\Gamma$ are as in
assumptions of \ref{group} and $\Gamma\subseteq S$.\\
Then the forcing notion $\bQ_{K_2,h}^3$ is really
$(\family_0,\hat{\family}_1,D)$--complete. 
\end{proposition}

\Proof Similar. \QED

\section{The iteration theorem}
In this section we will prove the preservation theorem needed for {\bf Case
A}. Let us start with some explanations which (hopefully) will help the reader
to understand what and why we do to get our result. 

We would like to prove that if $\bar{\bQ}=\langle\bP_i,\nbQ_i: i<\gamma
\rangle$ is a $(<\kappa)$--support iteration, $(\family_0,\hat{\family}_1,D)$
are as in \ref{really} then:
\begin{description}
\item[\quad if] $\bar{N}=\langle N_i:i\leq\lambda\rangle$ is an increasing
continuous sequence of elementary submodels of $\Hchi$, $\|N_i\|=\lambda$,
$\lambda+1\subseteq N_i$, $(N_i)^{<\lambda}\subseteq N_i$ for non-limit $i$,
and for some $\bar{a}\in\hat{\family}_1$ and $X\in D_{\bar{a}}$
\[(\forall i\in X)(N_i\cap\mu^*=a_i\quad\&\quad N_{i+1}\cap\mu^*=a_{i+1})\]
and $p\in\bP_\gamma\cap N_0$
\item[then]  there is a condition $q\in\bP_\gamma$ stronger than $p$ and
$(N_\lambda,\bP_\gamma)$--generic. 
\end{description}
For each $\nbQ_i$ we may get respective $q$, but the problem is with the
iteration. We can start with increasing successively $p$ to $p_i\in N_\lambda$
($i<\lambda$) and we can keep meeting dense sets due to
$(<\lambda)$--completeness. But the main question is: {\em why} is there a
limit? For each $\alpha\in\gamma\cap N_\lambda$ we have to make sure that the
sequence $\langle p_i(\alpha): i<\lambda\rangle$ has an upper bound in
$\nbQ_\alpha$, but for this we need information which is a $\bP_\gamma$--name
which does not belong to $N_\lambda$, e.g.\ if $\nbQ_i$ is $\bQ^2_{\name{
\bar{C}},\name{\bar{h}}}$ we need to know $\name{C}_{N\cap\kappa,
\name{h}_{N\cap\kappa}}$. {\em But} for each $i$, the size of the information
needed is $<\lambda$. 

As the life in our context is harder than for proper forcing iterations, we
have to go back to pre--proper tools and methods and we will use trees of
names (see \cite{Sh:64}). A tree of conditions is essentially a
non-deterministic condition; in the limit we will show that {\em some} choice
of a branch through the tree does the job.

\noindent [Note that one of difficulties one meets here is that we cannot
diagonalize objects of type $\lambda\times\omega$ when $\lambda>\aleph_0$.]
     
\begin{definition}
\begin{enumerate}
\item A tree $(T,<)$ is {\em normal} if for each $t_0,t_1\in T$,

if $\{s\in T:s<t_0\}=\{s\in T:s<t_1\}$ has no last element

then $t_0=t_1$.

\item For an ordinal $\gamma$, $\Tr(\gamma)$ stands for the family of all
triples 
\[\tree=(T^\tree,<^\tree,\rk^\tree)\]
such that $(T^\tree,<^\tree)$ is a normal tree and $\rk^\tree:T^\tree
\longrightarrow\gamma+1$ is an increasing continuous function.  

We will keep the convention that $\tree^x_y=(T^x_y,<^x_y,\rk^x_y)$. Sometimes
we may write $t\in\tree$ instead $t\in T^\tree$ (or $t\in T$).
\end{enumerate}
\end{definition}
The main case and examples we have in mind are triples $(T,<,\rk)$ such that
for some $w\subseteq\gamma$ (where $\gamma$ is the length of our iteration),
$T$ is a family partial functions such that:
\[(\forall t\in T)(\dom(t)\mbox{ is an initial segment of }w\quad\mbox{ and }
\quad (\forall\alpha\in w)(t\rest \alpha\in T));\]
the order is the inclusion and the function $\rk$ is given by
\[\rk(t)=\min\{\alpha\in w\cup\{\gamma\}: \dom(t)=\alpha\cap w\},\]
(see \ref{standleg}). Here we can let $N_\lambda\cap\gamma=\{\alpha_\xi:
\xi<\lambda\}$. Defining $p_i$ we are thinking of why $\langle
p_j(\alpha_\xi):j<\lambda\rangle$ will have an upper bound. Now
$\lambda\times\lambda$ {\em has} a diagonal. 

\noindent Note: {\em starting to take care of $\alpha_\xi$ sometime} is a
reasonable strategy, so in stage $i<\lambda$ we care about $\{\alpha_\xi:\xi
<i\}$ only. 

\noindent But what does this mean to do it? We have to guess the relevant
information which is a $\bP_{\alpha_\xi}$--name and is not present. 

\noindent What do we do? We cover {\em all}. So the tree $\tree_\zeta$ will
consist of objects $t$ which are guesses on what is $\langle${\em information
for }$\alpha_\varepsilon$ {\em up to $\zeta^{\rm th}$ stage}$:\varepsilon<
\zeta\rangle$. Of course we should not inflate, e.g. $\langle p_t^\zeta: t\in
\tree_\xi\rangle\in N_\lambda$.

\noindent It is very nice to have {\em an open option} so that in stage
$\lambda$ we can choose the most convenient branch. But we need to go into all
dense sets and then we have to pay an extra price for having an extra
luggage. We need to put {\em all} the $p_i$'s into a dense set (which is
trivial for a single condition). What will help us in this task is the strong
$\family_0$--completeness. Without this {\em big brother to pay our bills},
our schema would have to fail: we do have some ZFC theorems which put
restrictions on the iteration.

\begin{definition}
Let $\bar{\bQ}=\langle\bP_i,\nbQ_i: i<\gamma\rangle$ be a $(<\kappa)$--support
iteration. 
\begin{enumerate}
\item We define
\[\hspace{-0.5cm}
\begin{array}{r}
\ftr(\bar{\bQ})\stackrel{\rm def}{=}\big\{\bar{p}=\langle p_t: t\in T^\tree
\rangle:\tree\in\Tr(\gamma),\ (\forall t\in T^\tree)(p_t\in\bP_{\rk(t)})\
\mbox{ and }\ \\
(\forall s,t\in T^\tree))(s<t\ \Rightarrow\ p_s=p_t\rest\rk(s))\big\},
  \end{array}
\]
and
\[\hspace{-0.5cm}
\begin{array}{r}
\ftrw(\bar{\bQ})\stackrel{\rm def}{=}\big\{\bar{p}=\langle p_t: t\in T^\tree
\rangle:\tree\in\Tr(\gamma),\ (\forall t\in T^\tree)(p_t\in\bP_{\rk(t)})\
\mbox{ and }\ \\
(\forall s,t\in T^\tree))(s<t\ \Rightarrow\ p_s\geq p_t\rest\rk(s))\big\}.
  \end{array}
\]
We may write $\langle p_t: t\in\tree\rangle$. Abusing notation, we mean
$\bar{p}\in \ftrw(\bar{\bQ})$ (and $\bar{p}\in \ftr(\bar{\bQ})$) determines
$\tree$ and we call it $\tree^{\bar{p}}$ (or we may forget and write
$\dom(\bar{p})$). 

Adding primes to $\ftr$, $\ftrw$ means that we allow $p_t(\beta)$ be (a
$\bP_\beta$--name for) an element of the completion $\hat{\nbQ}_\beta$ of
$\nbQ_\beta$. Then $p_t$ is an element of $\bP_{\rk(t)}'$ --- the
$(<\kappa)$--support iteration of the completions $\hat{\nbQ}_\beta$ (see
\ref{uzupelnienie}).   
 
\item If $\tree\in \Tr(\gamma)$, $\bar{p},\bar{q}\in\ftrw'(\bar{\bQ})$,
$\dom(\bar{p})=\dom(\bar{q})=T^{\tree}$ then we let 
\[\bar{p}\leq\bar{q}\quad\mbox{ if and only if }\quad (\forall t\in T^\tree)
(p_t\leq q_t).\]
\item Let $\tree_1,\tree_2\in\Tr(\gamma)$. We say that a surjection $f:T_2
\stackrel{\rm onto}{\longrightarrow} T_1$ is {\em a projection} if for each
$s,t\in T_2$  
\begin{description}
\item[($\alpha$)]\quad $s\leq_2 t\quad\Rightarrow\quad f(s)\leq_1 f(t)$ and
\item[($\beta$)]\quad  $\rk_2(t)\leq\rk_1(f(t))$.
\end{description}
\item Let $\bar{p}^0,\bar{p}^1\in\ftrw'(\bar{\bQ})$, $\dom(\bar{p}^\ell)=
\tree_\ell$ ($\ell<2$) and $f:T_1\longrightarrow T_0$ be a projection. Then we
will write $\bar{p}^0\leq_f\bar{p}^1$ whenever for all $t\in T_1$
\begin{description}
\item[$(\alpha)$]\quad $p^0_{f(t)}\rest\rk_1(t)\leq_{\bP_{\rk_1(t)}'} p^1_t$
\qquad\quad and
\item[$(\beta)$]\quad if $i<\rk_1(t)$ then
\[p_t^1\rest i\forces_{\bP_i}\mbox{``}p^0_{f(t)}(i)\neq p_t^1(i)\quad
\Rightarrow\quad (\exists q\in\nbQ_i)(p^0_{f(t)}(i)\leq_{\hat{\nbQ}_i} q
\leq_{\hat{\nbQ}_i} p_t^1(i))\mbox{''}.\] 
\end{description}
\end{enumerate}
\end{definition}
The projections play the key role in the iteration lemma. Therefore, to
make the presentation more clear we will restrict ourselves to the case we
actually need. 

You may think of $\gamma$ as the length of the iteration, and let $\{
\beta_\xi: \xi<\lambda\}$ list $N\cap\gamma$, $w=\{\beta_\xi:\xi<\alpha\}$. We
are trying to build a generic condition for $(\bP_\gamma,N)$ by approximating
it by a sequence of trees of conditions. In the present tree we are at stage
$\alpha$. Now, for $t\in\tree$, $t(i)$ is a guess on the information needed to
construct a generic for $(N[\name{G}_{\bP_i}],\nbQ_i[\name{G}_{\bP_i}])$, more
exactly the $\alpha$--initial segment of it.

\begin{definition}
\label{standleg}
Let $\gamma$ be an ordinal.
\begin{enumerate}
\item Suppose that $w\subseteq\gamma$ and $\alpha$ is an ordinal. We say that
$\tree\in\Tr(\gamma)$ is {\em a standard $(w,\alpha)^\gamma$--tree} if
\begin{description}
\item[($\alpha$)] $(\forall t\in T^{\tree})(\rk^\tree(t)\in w\cup\{\gamma\})$,
\item[($\beta$)]  if $t\in T^\tree$, $\rk^\tree(t)=\varepsilon$ then $t$ is a
sequence $\langle t_i: i\in w\cap\varepsilon\rangle$, where each $t_i$ is a
sequence of length $\alpha$,
\item[($\gamma$)] $<^\tree$ is the extension (inclusion) relation.
\end{description}
[In $(\beta)$ above we may demand that each $t_i$ is a function with domain
$[i,\alpha)$, $i<\alpha$, but we can use a default value $*$ below $i$. Note
that $T^\tree$ determines $\tree$ in this case; $\langle\rangle$ is the root
of $T$.] 
\item Suppose that $w_0\subseteq w_1\subseteq\gamma$, $\alpha_0<\alpha_1$ and
$\tree=(T,<,\rk)$ is a standard $(w_1,\alpha_1)^\gamma$--tree. We define {\em
the projection $\proj^{(w_1,\alpha_1)}_{(w_0,\alpha_0)}(\tree)$ of $\tree$
onto $(w_0,\alpha_0)$} as $(T^*,<^*,\rk^*)$ such that:

$T^*=\{\langle t_i\rest\alpha_0: i\in w_0\cap\rk(t)\rangle: t=\langle t_i:
i\in w_1\cap\rk(t)\rangle\in T\}$,

$<^*$ is the extension relation,

$\rk^*(\langle t_i\rest\alpha_0: i\in w_0\cap\rk(t)\rangle)=\min(w_0\cup
\{\gamma\}\setminus\rk(t))$ for $t\in T$.

\noindent [Note that $\proj^{(w_1,\alpha_1)}_{(w_0,\alpha_0)}(\tree)$ is a
standard $(w_0,\alpha_0)^\gamma$--tree.]
\item If $w_0\subseteq w_1\subseteq\gamma$, $\alpha_0<\alpha_1$, $\tree_1=
(T_1,<_1,\rk_1)$ is a standard $(w_1,\alpha_1)^\gamma$--tree and $\tree_0=
(T_0,<_0,\rk_0)=\proj^{(w_1,\alpha_1)}_{(w_0,\alpha_0)}(\tree_1)$ then the
mapping 
\[T_1\ni\langle t_i: i\in w_1\cap\rk_1(t)\rangle\ \longmapsto\ \langle t_i
\rest\alpha_0: i\in w_0\cap\rk_1(t)\rangle\in T_0\]
is denoted by $\proj^{\tree_1}_{\tree_0}$ (or $\proj^{(w_1,\alpha_1)}_{(w_0,
\alpha_0)}$).

\noindent [Note that $\proj^{\tree_1}_{\tree_0}$ is a projection from
$\tree_1$ onto $\tree_0$.]
\item We say that $\bar{\tree}=\langle\tree_\alpha: \alpha<\alpha^*\rangle$ is
{\em a legal sequence of standard $\gamma$--trees} if for some
$\bar{w}=\langle w_\alpha:\alpha<\alpha^*\rangle$ we have
\begin{description}
\item[($\alpha$)] $\bar{w}$ is an increasing continuous sequence of subsets of
$\gamma$,
\item[($\beta$)]  for each $\alpha<\alpha^*$, $\tree_\alpha$ is a standard
$(w_\alpha,\alpha)^\gamma$--tree,
\item[($\gamma$)] if $\alpha<\beta<\alpha^*$ then $\tree_\alpha=
\proj^{(w_\alpha,\alpha)}_{(w_\beta,\beta)}(\tree_\beta)$.
\end{description}
\item For a legal sequence $\bar{\tree}=\langle\tree_\alpha:\alpha<\alpha^*
\rangle$ of standard $\gamma$--trees, $\alpha^*$ a limit ordinal, we define
the inverse limit $\inver(\bar{\tree})$ of $\bar{\tree}$ as a triple 
\[(T^{\inver(\bar{\tree})},<^{\inver(\bar{\tree})},\rk^{\inver(\bar{\tree})})\]
such that 
\begin{description}
\item[(a)] $T^{\inver(\bar{\tree})}$ consists of all sequences $t$ such that
\begin{description}
\item[(i)]    $\dom(t)$ is an initial segment of $w\stackrel{\rm def}{=}
\bigcup\limits_{\alpha<\alpha^*} w_\alpha$ (not necessarily proper),
\item[(ii)]   if $i\in\dom(t)$ then $t_i$ is a sequence of length $\alpha^*$,
\item[(iii)]  for each $\alpha<\alpha^*$, $\langle t_i\rest\alpha: i\in
w_\alpha\cap\dom(t)\rangle\in T_\alpha$,
\end{description}
\item[(b)] $<^{\inver(\bar{\tree})}$ is the extension relation,
\item[(c)] $\rk^{\inver(\bar{\tree})}(t)=\min(w\cup\{\gamma\}\setminus
\dom(t))$ for $t\in T^{\inver(\bar{\tree})}$.
\end{description}
[Note that it may happen that $T^{\inver(\bar{\tree})}=\{\langle\rangle\}$,
however not if $\bar{\tree}$ is continuous, see below.]
\item A legal sequence of standard $\gamma$--trees $\bar{\tree}=\langle
\tree_\alpha: \alpha<\alpha^*\rangle$ is {\em continuous} if $\tree_\alpha=
\inver(\tree_\beta:\beta<\alpha)$ for each limit $\alpha<\alpha^*$. 
\end{enumerate}
\end{definition}

\begin{proposition}
Suppose that $\alpha_0,\alpha_1,\alpha_2,\gamma$ are ordinals such that
$\alpha_0\leq\alpha_1\leq\alpha_2$. Let $w_{\alpha_0}\subseteq w_{\alpha_1}
\subseteq w_{\alpha_2}\subseteq\gamma$. Assume that for $\ell<3$, $\tree_\ell$
are standard $(w_\ell,\alpha_\ell)^\gamma$--trees such that
\[\tree_0=\proj^{(w_1,\alpha_1)}_{(w_0,\alpha_0)}(\tree_1)\quad\mbox{ and
}\quad\tree_1=\proj^{(w_2,\alpha_2)}_{(w_1,\alpha_1)}(\tree_2).\]
Then $\tree_0=\proj^{(w_2,\alpha_2)}_{(w_0,\alpha_0)}(\tree_2)$ and
$\proj^{\tree_2}_{\tree_0}=\proj^{\tree_1}_{\tree_0}\comp\proj^{\tree_2}_{
\tree_1}$.\\
Moreover, if $\bar{p}^\ell=\langle p^\ell_t: t\in T_\ell\rangle\in
\ftr'(\bar{\bQ})$ (for $\ell<3$) are such that $\bar{p}^0\leq_{\proj^{
\tree_1}_{\tree_0}}\bar{p}^1$ and $\bar{p}^1\leq_{\proj^{\tree_2}_{\tree_1}}
\bar{p}^2$ then $\bar{p}^0\leq_{\proj^{\tree_2}_{\tree_0}}\bar{p}^2$. \QED
\end{proposition}

\begin{proposition}
\label{liminver}
Let $\gamma,\alpha^*$ be ordinals, $\alpha^*$ limit, and let $\bar{\tree}=
\langle \tree_\alpha: \alpha<\alpha^*\rangle$ be a continuous legal sequence of
standard $\gamma$--trees. 
\begin{enumerate}
\item The inverse limit $\inver(\bar{\tree})$ is a standard
$(\bigcup\limits_{\alpha<\alpha^*} w_\alpha,\alpha^*)^\gamma$--tree and each
$\tree_\alpha$ is a projection of $\inver(\bar{\tree})$ onto
$(w_\alpha,\alpha)$ and the respective projections commute. (Here,
$w_\alpha\subseteq\gamma$ is such that $\tree_\alpha$ is a standard
$(w_\alpha,\alpha)^\gamma$--tree) 

\noindent [So we do not cheat: $\inver(\bar{\tree})$ is really the inverse
limit of $\bar{\tree}$.]
\item If $\lambda^{<\lambda}=\lambda$, $\alpha^*<\lambda$ and $\|T_\alpha\|\leq
\lambda$ for each $\alpha<\alpha^*$ then $\|T^{\inver(\bar{\tree})}\|\leq
\lambda$.
\item If $\alpha^*<\lambda$, $\kappa=\lambda^+$, $\bar{\bQ}=\langle\bP_\xi,
\nbQ_\xi:\xi<\gamma \rangle$ is a $(<\kappa)$--support iteration of
$(<\lambda)$--complete forcing notions and $\bar{p}^\alpha=\langle p^\alpha_t:
t\in\tree_\alpha\rangle\in\ftr'(\bar{\bQ})$ (for each $\alpha<\alpha^*$) are
such that $|T_\alpha|\leq\lambda$ for $\alpha<\alpha^*$ and
\[\beta<\alpha<\alpha^*\quad\Rightarrow\quad \bar{p}^\beta
\leq_{\proj^{\tree_\alpha}_{\tree_\beta}}\bar{p}^\alpha\]
then there is $\bar{p}^{\alpha^*}=\langle p^{\alpha^*}_t: t\in
\inver(\bar{\tree})\rangle\in\ftr'(\bar{\bQ})$ such that  
\[(\forall\alpha<\alpha^*)(\bar{p}^\alpha\leq_{\proj^{\inver(
\bar{\tree})}_{\tree_\alpha}}\bar{p}^{\alpha^*}).\]
\end{enumerate}
\end{proposition}

\Proof 1)\ \ \ Should be clear: just read the definitions.

\noindent 2)\ \ \ It follows from the following inequalities:
\[\|T^{\inver(\bar{\tree})}\|\leq\prod_{\alpha<\alpha^*}\|T_\alpha\|\leq
\lambda^{<\lambda}=\lambda.\]

\noindent 3)\ \ \ For each $t\in\inver(\bar{\tree})$ we define a condition
$p^{\alpha^*}_t\in\bP'_\gamma$ as follows. Let $t^\alpha=\proj^{\inver(
\bar{\tree})}_{\tree_\alpha}(t)$ (for $\alpha<\alpha^*$). We know that the
sequence $\langle p^\alpha_{t^\alpha}\rest\rk^{\inver(\bar{\tree})}(t):
\alpha<\alpha^*\rangle$ is increasing (remember $\rk^{\inver(\bar{\tree})}(t)
\leq\rk_\alpha(t_\alpha)$ for each $\alpha<\alpha^*$) and $p^{\alpha^*}_t$ is
supposed to be an upper bound to it (and $p^{\alpha^*}_t\in\bP'_{\rk^{\inver(
\bar{\tree})}(t)}$). We define $p^{\alpha^*}_t$ quite straightforward. We let
\[\dom(p^{\alpha^*}_t)=\bigcup\{\dom(p^\alpha_{t^\alpha})\cap\rk^{\inver(
\bar{\tree})}(t): \alpha<\alpha^*\}\] 
and next we inductively define $p^{\alpha^*}_t(i)$ for $i\in\dom(p^{\alpha^*
}_t)$. Assume we have defined $p^{\alpha^*}_t\rest i$ such that 
\[(\forall\alpha<\alpha^*)(p^\alpha_{t^\alpha}\rest i\leq_{\bP'_i}
p^{\alpha^*}_t\rest i).\]
Then (remembering our convention that if $i\notin\dom(p)$ then $p(i)={\bf
0}_{\nbQ_i}$) 
\[\begin{array}{r}
p^{\alpha^*}_t\rest i\forces\mbox{`` the sequence }\langle
p^\alpha_{t^\alpha}(i): \alpha<\alpha^*\rangle\subseteq\hat{\nbQ}_i\mbox{ is
$\leq_{\hat{\nbQ}_i}$--increasing and}\\
\alpha<\beta<\alpha^*\ \&\ p^\alpha_{t^\alpha}(i)\neq p^\beta_{t^\beta}(i)
\quad\Rightarrow\quad (\exists q\in\nbQ_i)(p^\alpha_{t^\alpha}(i)
\leq_{\hat{\nbQ}_i} q\leq_{\hat{\nbQ}_i} p^\beta_{t^\beta}(i))\\
\mbox{and $\nbQ_i$ is $(<\lambda)$--complete and $\alpha^*<\lambda$''}.
  \end{array}
\] 
Hence we find a $\bP'_i$--name $p^{\alpha^*}_t(i)$ (and we take the
$<^*_\chi$--first such a name) such that 
\[p^{\alpha^*}_t\rest i\forces\mbox{``}p^{\alpha^*}_t(i)\in\hat{\nbQ}_i\mbox{
is the least upper bound of }\langle p^\alpha_{t^\alpha}(i):\alpha<\alpha^*
\rangle\mbox{ in }\hat{\nbQ}_i\mbox{''}.\]
Now one easily checks that $p^{\alpha^*}_t\in\bP_{\rk^{\inver(\bar{\tree})}
(t)}'$. Consequently the condition $p^{\alpha^*}_t$ is as required. But why
does $\langle p^{\alpha^*}_t: t\in\inver(\bar{\tree})\rangle\in\ftr'(
\bar{\bQ})$? We still have to argue that
\[(\forall s,t\in\inver(\bar{\tree}))(s<t\quad\Rightarrow\quad p^{\alpha^*}_s
=p^{\alpha^*}_t\rest\rk^{\inver(\bar{\tree})}(s)).\]
For this note that if $s<t$ are in $\inver(\bar{\tree})$ and
$s_\alpha,t_\alpha$ are their projections to $\tree_\alpha$ then
$s_\alpha\leq_\alpha t_\alpha$ and $p^\alpha_{s_\alpha}= p^\alpha_{t_\alpha}
\rest\rk_\alpha(s_\alpha)$ and $\rk^{\inver(\bar{\tree})}(s)\leq\rk_\alpha
(s_\alpha)$. Thus clearly $\dom(p^{\alpha^*}_s)=\dom(p^{\alpha^*}_t)\cap
\rk^{\inver(\bar{\tree})}(s)$. Next, by induction on $i\in
\dom(p^{\alpha^*}_t)\cap\rk^{\inver(\bar{\tree})}(s)$ we show that
$p^{\alpha^*}_s(i)=p^{\alpha^*}_t(i)$. Assume we have proved that
$p^{\alpha^*}_t\rest i=p^{\alpha^*}_s\rest i$ and look at the way we defined
the respective values at $i$. We looked there at the sequences $\langle
p^\alpha_{t_\alpha}(i): \alpha<\alpha^*\rangle$, $\langle p^\alpha_{s_\alpha}
(i): \alpha<\alpha^*\rangle$ and we chosen the $<^*_\chi$--first names for the
least upper bounds to them. But $i<\rk_\alpha (s_\alpha)$ for all $\alpha<
\alpha^*$, so the two sequences are equal and the choice was the same. \QED  
 
\begin{proposition}
\label{onesteptree}
Assume that $\family_0\subseteq [\mu^*]^{\leq\lambda}$ and $\bar{\bQ}=\langle
\bP_\alpha,\nbQ_\alpha:\alpha<\gamma\rangle$ is a $(<\kappa)$--support
iteration of $(<\lambda)$--complete strongly $\family_0$--complete forcing
notions, and $\name{x}_\alpha$ (for $\alpha<\gamma$) are $\bP_\alpha$--names
such that
\[\forces_{\bP_\alpha}\mbox{``$\name{x}_\alpha$ witnesses the {\em most} in
\ref{good}(2) for $\nbQ_\alpha$''.}\]
Further suppose that
 
\begin{description}
\item[($\alpha$)] $N\prec\Hchi$ is $(\lambda,\family_0)$--good (see
\ref{good}), $\langle\name{x}_\alpha:\alpha<\gamma\rangle, \alpha_0,\bar{\bQ},
\ldots\in N$, 
\item[($\beta$)]  $0\in w_0\subseteq w_1\in N\cap [\gamma]^{<\lambda}$,
$\alpha_0<\lambda$ is an ordinal, $\alpha_1=\alpha_0+1$,
\item[($\gamma$)] $\tree_0=(T_0,<_0,\rk_0)\in N$ is a standard $(w_0,
\alpha_0)^\gamma$--tree, $\|T_0\|\leq\lambda$,
\item[($\delta$)] $\bar{p}=\langle p_t: t\in T_0\rangle\in\ftr'(\bar{\bQ})\cap
N$, 
\item[($\varepsilon$)] $\tree_1=(T_1,<_1,\rk_1)$ is such that 

$T_1$ consists of all sequences $t=\langle t_i: i\in\dom(t)\rangle$ such that
$\dom(t)$ is an initial segment of $w_1$, and
\begin{itemize}
\item each $t_i$ is a sequence of length $\alpha_1$,
\item $t'\stackrel{\rm def}{=}\langle t_i\rest\alpha_0: i\in\dom(t)\cap
w_0\rangle\in T_0$,
\item if $i\in\dom(t)\setminus w_0$, $\alpha<\alpha_0$ then $t_i(\alpha)=*$, 
\item for some $j(t)\in\dom(t)\cup\{\gamma\}$, 

$t_i(\alpha_0)$ is $*$ for every $i\in\dom(t)\setminus j(t)$, and for each
$i\in\dom(t)\cap j(t)$

$t_i(\alpha_0)\in N$ is a $\bP_i$--name for an element of $\nbQ_i$, 
\end{itemize}
$\rk_1(t)=\min(w_1\cup\{\gamma\}\setminus\dom(t))$ and $<_1$ is the extension
relation. 
\end{description}
Then 
\begin{description}
\item[(a)] $\tree_1$ is a standard $(w_1,\alpha_1)^\gamma$--tree,
$\|T_1\|=\lambda$,
\item[(b)] $\tree_0$ is the projection of $\tree_1$ onto $(w_0,\alpha_0)$,
\item[(c)] there is $\bar{q}=\langle q_t: t\in T_1\rangle\in \ftr'(\bar{\bQ})$
such that
\begin{description}
\item[(i)]   $\bar{p}\leq_{\proj^{\tree_1}_{\tree_0}}\bar{q}$,
\item[(ii)]  if $t\in T_1\setminus\{\langle\rangle\}$ and $(\forall i\in
\dom(t))(t_i(\alpha_0)\neq *)$ then the condition $q_t\in\bP_{\rk_1(t)}'$ 
is an upper bound in $\bP'_{\rk_1(t)}$ of a $\bP_{\rk_1(t)}$--generic sequence
over $N$, and for every $\beta\in\dom(q_t)=N\cap\rk_1(t)$, $q_t(\beta)$ is (a
name for) the least upper bound in $\hat{\nbQ}_\beta$ of the family of all
$r(\beta)$ for $r$ from the generic set (over $N$) generated by $q_t$, 
\item[(iii)] if $t\in T_1$, $t'=\proj^{\tree_1}_{\tree_0}(t)\in T_0$, $i\in
\dom(t)$ and $t_i(\alpha_0)\neq *$ then
\[q_t\rest i\forces_{\bP_i}\mbox{``}p_{t'}(i)\leq_{\hat{\nbQ}_i}t_i(\alpha_0)\
\Rightarrow\ t_i(\alpha_0)\leq_{\hat{\nbQ}_i} q_t(i) \mbox{''},\] 
\item[(iv)]  $q_{\langle\rangle}=p_{\langle\rangle}$ and 

if $t\in T_1\setminus\{\langle\rangle\}$ and $j(t)<\gamma$ then $q_t=q_{t\rest
j(t)}\cup p_{t'}\rest[j(t),\rk_1(t))$, where $t'=\proj^{\tree_1}_{\tree_0}(t)
\in T_0$.   
\end{description}
\end{description}
\end{proposition}

\Proof Clauses (a) and (b) should be clear.\\
(c)\ \ \ Let $\langle t_\zeta:\zeta<\lambda\rangle$ list with
$\lambda$-repetitions all elements $t$ of $T_1\setminus\{\langle\rangle\}$
such that $(\forall i\in\dom(t))(t_i(\alpha_0)\neq *)$. For $\alpha\in w_1\cup
\{\gamma\}$ let $\langle {\cal I}^\alpha_\zeta:\zeta<\lambda\rangle$ enumerate
all open dense subsets of $\bP_\alpha$ from $N$. By induction on
$\zeta<\lambda$ choose $r_\zeta$ such that
\begin{itemize}
\item $r_\zeta\in \bP_{\rk_1(t_\zeta)}\cap N$,
\item if $t'=\proj^{\tree_1}_{\tree_0}(t_\zeta)$ then $p_{t'}\rest \rk_1(
t_\zeta)\leq_{\bP'_{\rk_0(t_\zeta)}} r_\zeta$ and for $i\in\dom(t_\zeta)$
\[r_\zeta\rest i\forces_{\bP_i}\mbox{``}p_{t'}(i)\leq_{\hat{\nbQ}_i}t_i(
\alpha_0)\ \Rightarrow\ t_i(\alpha_0)\leq_{\hat{\nbQ}_i}r_\zeta(i)\mbox{''},\] 
\item $r_\zeta\in {\cal I}^{\rk_1(t_\zeta)}_\xi$ for all $\xi\leq\zeta$,
\item if $t\in T_1$, $\xi<\zeta$, $t\leq_1 t_\xi$, $t\leq_1 t_\zeta$ (e.g.\
$t=t_\xi=t_\zeta$) then $r_\xi\rest\rk_1(t)\leq_{\bP_{\rk_1(t)}} r_\zeta\rest
\rk_1(t)$. 
\end{itemize}
Since we have assumed that all the $\nbQ_\alpha$'s are (names for)
$(<\lambda)$--complete forcing notions there are no difficulties in carrying
out the above construction. [First, working in $N$, choose $r^*_\zeta\in
\bP_{\rk_1(t_\zeta)}\cap N$ satisfying the second and the fourth demand. How?
Declare 
\[\dom(r^*_\zeta)=[w_1\cup\bigcup\{\dom(r_\xi): \xi<\zeta\}\cup\dom(p_{
\proj^{\tree_1}_{\tree_0}(t_\zeta)})]\cap\rk_1(t_\zeta)\]
and by induction on $i$ define $r^*_\zeta(i)$ using $(<\lambda)$--completeness
of $\nbQ_i$ and taking care of the respective demands (similar to the choice
of $q_t$ done in details below). Next use the $(<\lambda)$--completeness (see
\ref{firstiter}) to enter all ${\cal I}^{\rk_1(t_\zeta)}_\xi$ for $\xi\leq
\zeta$. Note that the sequence $\langle {\cal I}^{\rk_1(t_\zeta)}_\xi: \xi
\leq\zeta\rangle$ is in $N$, so we may choose the respective $r_\zeta\geq
r^*_\zeta$ in $N$.]\\ 
Now we may define $\bar{q}=\langle q_t:t\in T_1\rangle\in\ftr'(\bar{\bQ})$. If
$t\in T_1$ is such that $j(t)<\rk_1(t)$ then $q_t$ is defined from $q_{t\rest
j(t)}$ and $\bar{p}$ by demand {\bf (c)(iv)}. So we have to define $q_t$ for
these $t\in T_1$ that $(\forall i\in\dom(t))(t(\alpha_0)\neq *)$ (and
$t\neq\langle\rangle$) only. So let $t\in T_1\setminus\{\langle\rangle\}$ be
of this type. Let 
\[\dom(q_t)=\bigcup\{\dom(r_\zeta): \zeta<\lambda\ \&\ t\leq_1 t_\zeta\}\cap
\rk_1(t)\subseteq N\]
and by induction on $i\in\dom(q_t)$ we define $q_t(i)$ (a $\bP_i$--name for a
member of $\hat{\nbQ}_i$). So suppose that $i\in\dom(q_t)$ and we have defined
$q_t\rest i\in\bP_i'$ in such a way that
\[(\forall\zeta<\lambda)(t\leq_1t_\zeta\quad\Rightarrow\quad r_\zeta\rest i
\leq_{\bP'_i} q_t\rest i).\]
Note that this demand implies that $q_t\rest i\in\bP'_i$ is an upper bound of
a generic sequence in $\bP_i$ over $N$ (remember the choice of the
$r_\zeta$'s, and that $i\in N$ and there are many $\zeta<\lambda$ such that
$\rk_1(t_\zeta)=\gamma$, and all open dense subsets of $\bP_i$ from $N$ appear
in the list $\langle {\cal I}^\gamma_\zeta:\zeta<\lambda\rangle$) and
therefore  
\[q_t\rest i\forces_{\bP'_i}\mbox{``the model $N[\name{G}_{\bP'_i}]$ is
$(\lambda,\family_0)$--good''}\]
(remember \ref{seqfromN}). Look at the sequence $\langle r_\zeta(i): t\leq_1
t_\zeta\ \&\ i\in\dom(r_\zeta)\cap\rk_1(t)\rangle$. By the last two demands of
the choice of the $r_\zeta$'s we have
\[\begin{array}{r}
q_t\rest i\forces_{\bP_i}\mbox{``}\langle r_\zeta(i): t\leq_1 t_\zeta\ \&\
i\in\dom(r_\zeta)\cap\rk_1(t)\rangle\mbox{ is an increasing }\ \ \\
\nbQ_i\mbox{--generic sequence over }N[\name{G}_{\bP_i}]\mbox{''}.
  \end{array}\]
Consequently we may use the fact that $\nbQ_i$ is (a name for) strongly
$\family_0$--complete forcing notion and $\name{x}_i\in N$, and we take
$q_t(i)$ to be the $<^*_\chi$--first name for the least upper bound of this
sequence in $\hat{\nbQ}_i$.  

This completes the definition of $\bar{q}$. Checking that it is as required is
straightforward. \QED

\begin{theorem}
\label{iterweak}
Assume $\lambda^{<\lambda}=\lambda$, $\kappa=\lambda^+=2^\lambda\leq\mu^*$.\\
Suppose that $\family_0\in(\dkl)^+$, $\hat{\family}_1\in (\dlk[\family_0])^+$
and $D$ is a function such that $\dom(D)=\hat{\family}_1$ and for every
$\bar{a}\in\hat{\family}_1$ 
\[D(\bar{a})=D_{\bar{a}}\ \mbox{ is a normal filter on }\lambda.\]
Further suppose that $\bar{\bQ}=\langle\bP_i,\nbQ_i: i<\gamma\rangle$ is a
$(<\kappa)$--support iteration such that for each $i<\gamma$
\[\begin{array}{ll}
\forces_{\bP_i}&\mbox{`` }\nbQ_i\mbox{ is really $(\family_0,\hat{\family}_1,
D)$--complete with witness $\name{x}_i$}\\
\ &\ \mbox{ for {\em most} in \ref{really}(3)($\gamma$)''.}
  \end{array}\]
Then:
\begin{description}
\item[(a)] the forcing notion $\bP_\gamma$ is $(<\lambda)$--complete and
strongly $\family_0$--complete,
\item[(b)] if a sequence $\bar{N}=\langle N_i: i\leq\lambda\rangle$ is
$(\lambda,\kappa,\hat{\family}_1,D,\bP_\gamma)$--suitable (see
\ref{really}(1)) and $p\in\bP_\gamma\cap N_0$, $\langle \name{x}_i:i<\gamma
\rangle,\ \langle\family_0,\hat{\family}_1,D\rangle\in N_0$

then there is an $(N_\lambda,\bP_\gamma)$--generic condition $q\in\bP_\gamma$
stronger than $p$, 
\item[(c)] the forcing notion $\bP_\gamma$ is $(\family_0,\hat{\family}_1,
D)$--complete.
\end{description}
\end{theorem}

\Proof {\bf (a)} It is the consequence of \ref{firstiter} and \ref{iterstrong}.

\noindent {\bf (b)} Let $(X,\bar{a})$ be a suitable basis for $\bar{N}$, so 
$\bar{a}\in\hat{\family}_1$, $X\in D_{\bar{a}}$ and 
\[(\forall i\in X)((N_{i+1})^{<\lambda}\subseteq N_{i+1}\ \ \&\ \ N_{i+1}\cap
\mu^*=a_{i+1}).\]
We may assume that all members of $X$ are limit ordinals. Let $w_\lambda=
N_\lambda\cap\gamma$ (so $\|w_\lambda\|=\lambda$). Choose an increasing
continuous sequence $\langle w_\alpha: \alpha<\lambda\rangle$ such that
$\bigcup\limits_{\alpha<\lambda} w_\alpha=w_\lambda$ and for each $\alpha<
\lambda$ 
\[\|w_\alpha\|<\lambda,\quad w_\alpha\subseteq N_\alpha\cap\gamma,\quad 0\in
w_\alpha,\quad \mbox{ and }\quad \mbox{if $\alpha$ is limit then $w_\alpha=
w_{\alpha+1}$}\] 
(so then $\langle w_\beta: \beta\leq \alpha\rangle\in N_{\alpha+1}$). 

Now, by induction on $\alpha\leq\lambda$ we define a legal continuous sequence
of standard $\gamma$--trees $\langle\tree_\alpha: \alpha\leq\lambda\rangle$
and a sequence $\langle\bar{p}^\alpha:\alpha<\lambda\rangle$ such that
$\bar{p}^\beta=\langle p^\beta_t:t\in T_\beta\rangle\in\ftr'(\bar{\bQ})$ and
$\bar{p}^\beta\leq_{\proj^{\tree_\alpha}_{\tree_\beta}} \bar{p}^\alpha$ for
each $\beta<\alpha<\lambda$. 

{\em At stage $\alpha=0$ of the construction:}\\
$T_0$ consists of all sequences $t=\langle t_i:i\in\dom(t)\rangle$ such that
$\dom(t)$ is an initial segment of $w_0$ (not necessarily proper) and for each
$i\in\dom(t)$, $t_i$ is a sequence of length $0$ (i.e. $\langle\rangle$),

\noindent $\rk_0(t)=\min(w_0\cup\{\gamma\}\setminus\dom(t))$ and $<_0$ is the
extension relation; 

\noindent for each $t\in T_0$ we let $p^0_t=p\rest\rk_0(t)$ and finally
$\bar{p}^0=\langle p^0_t:t\in T_0\rangle\in\ftr'(\bar{\bQ})$.

\noindent [Note that $\tree_0=(T_0,<_0,\rk_0)\in N_0$ is a standard
$(w_0,0)^\gamma$--tree, $\bar{p}^0\in\ftr'(\bar{\bQ})\cap N_0$.]

{\em At stage $\alpha=\alpha_0+1$ of the construction:}\\
We have defined a standard $(w_{\alpha_0},\alpha_0)^\gamma$--tree
$\tree_{\alpha_0}\in N_{\alpha_0+1}$ and $\bar{p}^{\alpha_0}=\langle
p^{\alpha_0}_t: t\in\tree_{\alpha_0}\rangle\in\ftr'(\bar{\bQ})\cap N_{\alpha_0
+1}$. Now we consider two cases.

\noindent {\bf If} $\alpha_0\in X$ (so $N_{\alpha_0+1}$ is
$(\lambda,\family_0)$--good) then we apply the procedure of \ref{onesteptree}
inside $N_{\alpha_0+2}$ to $\tree_{\alpha_0}$, $\bar{p}^{\alpha_0}$,
$(w_{\alpha_0+1},\alpha_0+1)$ and $N_{\alpha_0+1}$ (in place of $\tree_0$,
$\bar{p}$, $(w_1,\alpha_1)$ and $N$ there) and we get a standard $(w_{\alpha_0
+1},\alpha_0+1)^\gamma$--tree $\tree_{\alpha_0}\in N_{\alpha_0+2}$ and 
$\bar{p}^{\alpha_0+1}=\langle p^{\alpha_0+1}_t: t\in T_{\alpha_0+1}\rangle
\in\ftr'(\bar{\bQ})\cap N_{\alpha_0+2}$ satisfying the demands
\ref{onesteptree}($\varepsilon$) and \ref{onesteptree}(a)--(c). 

\noindent {\bf If} $\alpha_0\notin X$ then we define
$\tree_{\alpha_0+1}$ as above but we cannot put any new genericity
requirements on $\bar{p}^{\alpha_0+1}$, so we just let $p^{\alpha_0+1}_t=
p^{\alpha_0}_{t'}\rest\rk_{\alpha_0+1}(t)$ where $t'=\proj^{\tree_{\alpha_0+
1}}_{\tree_{\alpha_0}}(t)$. 

\noindent [Note that in both cases $\tree_{\alpha_0}\in N_{\alpha_0+2}$ is a
standard $(w_{\alpha_0+1},\alpha_0+1)^\gamma$--tree, projection of
$\tree_{\alpha_0+1}$ onto $(w_{\alpha_0},\alpha_0)$ is $\tree_{\alpha_0}$ ,
$\bar{p}^{\alpha_0+1}\in\ftr'(\bar{\bQ})\cap N_{\alpha_0+2}$ and
$\bar{p}^{\alpha_0}\leq_{\proj^{\tree_{\alpha_0+1}}_{\tree_{\alpha_0}}} 
\bar{p}^{\alpha_0}$.] 

{\em At limit stage $\alpha$ of the construction:}\\
We let $\tree_\alpha=\inver(\langle\tree_\beta: \beta<\alpha\rangle)\in
N_{\alpha+1}$ and we choose $\bar{p}^\alpha=\langle p^\alpha_t: t\in T_\alpha
\rangle\in\ftr'(\bar{\bQ})\cap N_{\alpha+1}$ applying \ref{liminver} in
$N_{\alpha+1}$.  

\noindent [Note that the corresponding inductive assumptions hold true.]
\medskip

After the construction is carried out we may let $\tree_\lambda=\inver(\langle
\tree_\alpha: \alpha<\lambda\rangle)$. Then $\tree_\lambda$ is a standard
$(w_\lambda,\lambda)^\gamma$--tree, but no longer we have $\|T_\lambda\|\leq
\lambda$.  

Now, by induction on $\alpha\in w_\lambda\cup\{\gamma\}$ we choose conditions
$q_\alpha$ and $\bP_\alpha$--names $\name{X}_\alpha$, $\name{Y}_\alpha$ and
$\name{t}_\alpha$ such that 
\begin{description}
\item[(a)] $\forces_{\bP_\alpha}$``$\name{t}_\alpha\in T_\lambda\ \&\
\rk_\lambda(\name{t}_\alpha)=\alpha$'',
\item[(b)] $\forces_{\bP_\alpha}$``$\name{t}_\beta=\name{t}_\alpha\rest
\beta$'' for $\beta<\alpha$,
\item[(c)] $q_\alpha\in\bP_\alpha$, $\dom(q_\alpha)=w_\lambda\cap\alpha$,
\item[(d)] if $\beta<\alpha$ then $q_\beta=q_\alpha\rest \beta$,
\item[(e)] $q_\alpha\forces_{\bP_\alpha}$``$p^i_{\proj^{\tree_\lambda}_{
\tree_i}(\name{t}_\alpha)}\rest\alpha \in \name{G}_{\bP_\alpha}$'' for each
$i<\lambda$, 
\item[(f)] for each $\beta<\alpha$
\[\begin{array}{ll}
q_\alpha\forces_{\bP_\alpha}&\mbox{``$\name{X}_\beta=\{i<\lambda:
(\name{t}_{\beta+1})_{\beta}(i)\neq *\}\in D_{\bar{a}}$ and the sequence}\\
\ &\langle\langle(i,(\name{t}_{\beta+1})_\beta(i)), p^{i+1}_{\proj^{
\tree_\lambda}_{\tree_{i+1}}(\name{t}_{\beta+1})}(\beta)\rangle:\ i<\lambda
\ \&\ i\in\name{X}_\beta\rangle\\
\ &\mbox{is a result of a play of the game}\\
\ &{\cal G}^\heartsuit_{\langle N_i[\name{G}_\beta]: i\leq
\lambda\rangle,D,\name{Y}_\beta,\bar{a}}(\nbQ_\beta,p^{i_0}_{\proj^{ 
\tree_\lambda}_{\tree_{i_0}}(\name{t}_{\beta+1})}(\beta)),\\
\ &\mbox{[where $i_0<\lambda$ is the first such that $\beta\in
w_{i_0}$]},\\ 
\ &\mbox{ won by player $\com$'',}
\end{array} \]
\item[(g)] the condition $q_\alpha$ forces (in $\bP_\alpha$) that

``the sequence $\langle N_i[\name{G}_{\bP_\alpha}]: i\leq\lambda\rangle$ is
$(\lambda,\kappa,\hat{\family}_1,D,\name{\bQ}_\alpha)$--suitable and
$\name{Y}_\alpha\in D_{\bar{a}}$ is such that $\name{Y}_\alpha\subseteq X$ and
for every $i\in\name{Y}_\alpha$ we have  
\[(N_{i+1}[\name{G}_{\bP_\alpha}])^{<\lambda}\subseteq N_{i+1}[\name{G}_{
\bP_\alpha}]\quad\mbox{ and }\quad \ N_{i+1}[\name{G}_{\bP_\alpha}]\cap\V=
N_{i+1}\] 
and $i\in\name{X}_\xi$ for all $\xi\in\alpha\cap w_i$ (hence
$N_\lambda[G_{\bP_\alpha}]\cap\V= N_\lambda$)'' 
\end{description}
\medskip

\noindent{\sc Case 1:}\qquad $\alpha=0$.\\
We do not have much choice here: we let $q_0=\emptyset$, $\name{t}_0=\langle
\rangle\in\tree_\lambda$ and $\name{Y}_0=X$. Note that clauses {\bf (a)}--{\bf
(e)} and {\bf (g)} are trivially satisfied (for {\bf (g)} remember that
$(\bar{a},X)$ is a suitable basis for $\bar{N}$) and clause {\bf (f)} is not
relevant. 
\medskip

\noindent{\sc Case 2:}\qquad $\alpha=\beta+1$.\\
Arriving at this stage we have defined $q_\beta,\name{t}_\beta,\name{Y}_\beta$
and $\name{X}_\xi$ for $\xi<\beta$, and we want to choose $q_{\beta+1},
\name{t}_{\beta+1}, \name{Y}_{\beta+1}$ and $\name{X}_\beta$.\\
Suppose that $G_\beta\subseteq\bP_\beta$ is a generic filter over $\V$ such
that $q_\beta\in G_\beta$. Then (by clause {\bf (g)} at stage $\beta$) we have 
\[\begin{array}{ll}
\V[G_\beta]\models&\mbox{``the sequence }\langle N_i[G_\beta]:i\leq\lambda
\rangle \mbox{ is $(\lambda,\kappa,\hat{\family}_1,D,\nbQ^{G_\beta}_\beta
)$--suitable}\\
\ &\ \mbox{ and $(\bar{a},\name{Y}_\beta^{G_\beta})$ is a suitable base for
it}\\ 
\ &\ \mbox{ and $(\forall i\in\name{Y}^{G_\beta}_\beta)(\forall\xi\in\beta\cap
w_i)(i\in \name{X}^{G_\beta}_{\xi})$''} 
  \end{array}\]
Let $i_0=\min\{j<\lambda: \beta\in w_j\}$. We know that the player $\inc$
does not have any winning strategy in the game 
\[{\cal G}^\heartsuit_{\langle N_i[G_\beta]:i\leq\lambda\rangle,D,
\name{Y}^{G_\beta}_\beta\setminus (i_0+1),\bar{a}}(\nbQ^{G_\beta}_\beta,
p^{i_0}_{\proj^{\tree_\lambda}_{\tree_{i_0}}(\name{t}^{G_\beta}_\beta)}(
\beta)).\]
Now, using the interpretation of the game presented in \ref{renotation}, we
describe a strategy for player $\inc$ in this game. 
\begin{description}
\item[The strategy says:] {\em during a play $\com$ constructs a sequence\\
$\bar{s}=\langle s(i):i<\lambda\rangle$ of elements of $\nbQ^{G_\beta}_\beta
\cup\{*\}$; let 
\[r_i\stackrel{\rm def}{=}\proj^{\tree_\lambda}_{\tree_i}(\name{t}^{G_\beta}_{
\beta})\conc \bar{s}\rest i\in T_i,\]
and at the stage $i<\lambda$ of the game $\inc$ answers with
$p^{i+1}_{r_{i+1}}(\beta)$.  
}
\end{description}
We have to argue that the strategy described above is a legal one, i.e. that
it always says $\inc$ to play legal moves (assuming that $\com$ plays
according to the rules of the game). For this we show by induction on
$i<\lambda$ that really $r_i\in T_i$ and that if $s(i)\neq *$ then $p^{i+1}_{
r_{i+1}}(\beta)\in N_{i+2}[G_\beta]\cap\hat{\nbQ}^{G_\beta}_\beta$ is the
least upper bound of a $\nbQ^{G_\beta}_\beta$--generic filter over
$N_{i+1}[G_\beta]$ (to which $s(i)$ belongs) and if $s(i)=*$ then $p^{i+1}_{
r_{i+1}}$ is the least upper bound of conditions played by $\inc$ so far.\\
First note that $s(i)=*$ for all $i\leq i_0$ and therefore $r_{i_0+1}\in
T_{i_0+1}$ (just look at the successor stage of the construction of the
$\tree_\alpha$'s; remember that $\dom(\name{t}^{G_\beta}_\beta)=\beta$, so
adding $*$'s at level $\beta$ is allowed by \ref{onesteptree}($\varepsilon$)).
Note that $p^{i_0+1}_{r_{i_0+1}}(\beta)=p^{i_0}_{\proj^{\tree_\lambda}_{
\tree_{i_0}}(\name{t}^{G_\beta}_\beta)}(\beta)$ (remember
\ref{onesteptree}(c)(iv)).\\
If $i<\lambda$ is a limit ordinal above $i_0$, and we know already that
$r_j\in T_j$ for each $j<i$ then $r_i\in\tree_i=\inver(\langle\tree_j: j<i
\rangle)$ as clearly
\[j_0<j_1<i\quad\Rightarrow\quad \proj^{\tree_{j_1}}_{\tree_{j_0}}(r_{j_1})
=r_{j_0}.\]
Note that, by the limit stage of the construction of the $\tree_\alpha$'s and
\ref{liminver}(3) (actually by the construction there), the condition
$p^{i}_{r_i}(\beta)$ is the least upper bound of $\langle p^j_{r_j}(\beta):
j<i\rangle$ in $\hat{\nbQ}^{G_\beta}_\beta$.\\
Suppose now that we have $r_i\in T_i$, $i_0<i<\lambda$ and the player $\com$
plays $s(i)$. If $s(i)=*$ then easily $r_{i+1}\in T_{i+1}$ as adding stars at
``top levels'' does not make any problems (compare the case of
$i_0$). Moreover, as there, we have then 
\[p^{i+1}_{r_{i+1}}(\beta)=p^{i}_{\proj^{\tree_{i+1}}_{\tree_i}(r_{i+1})}(
\beta)=p^i_{r_i}(\beta).\]
If $s(i)\neq *$ then $s(i)\in N_{i+1}[G_\beta]\cap \nbQ^{G_\beta}_\beta$ is a
condition stronger than all conditions played by $\inc$ so far, and thus it is
stronger than $p^i_{r_i}(\beta)$. Moreover, in this case we necessarily have
$i\in \name{Y}_\beta^{G_\beta}$, so $i$ is limit and therefore $w_i=w_{i+1}$.
Hence $(\forall \xi\in w_{i+1}\cap\beta)(i\in \name{X}_\xi^{G_\beta})$. By
clause {\bf (f)} for $q_\beta$ we conclude that $(\forall\xi\in w_{i+1})((
\name{t}_\beta^{G_\beta})_\xi(i)\neq *)$. Therefore, if we look at the way
$\tree_{i+1}$ was constructed, we see that there is no collision in adding
$s(i)$ at the top (i.e. it is allowed by \ref{onesteptree}($\varepsilon$)).
Thus $r_{i+1}\in T_{i+1}$ and by \ref{onesteptree}(c)(ii) we know that
$p^{i+1}_{r_{i+1}}(\beta)\in N_{i+2}[G_\beta]\cap\hat{\nbQ}^{G_\beta}_\beta$
is the least upper bound of a $\nbQ^{G_\beta}_\beta$--generic sequence over
$N_{i+1}[G_\beta]$ to which $s(i)$ belongs (the last is due to
\ref{onesteptree}(c)(iii)).  

Thus we have proved that the strategy presented above is a legal strategy for
$\inc$. It cannot be the winning one, so there is a play $\bar{s}=\langle
s(i):i<\lambda\rangle$ (we give the moves of $\com$ only) in which $\com$
wins. Let $t_\alpha=\name{t}_\beta^{G_\beta}\conc \bar{s}$. We have actually
proved that $t_\alpha\in \tree_\lambda=\inver(\langle\tree_i:i<\lambda
\rangle)$. It should be clear that $\rk_\lambda(t_\alpha)=\alpha$ and
$\name{t}_\beta^{G_\beta}=t_\alpha\rest\beta$. Further let $q_\alpha(\beta)\in
\nbQ^{G_\beta}_\beta$ be any upper bound of $\bar{s}$ in $\nbQ^{G_\beta
}_\beta$ (there is one as $\com$ wins) and $X_\beta$ be the set $\{i<\lambda:
s(i)\neq *\}\in D_{\bar{a}}$. Note that then $q_\alpha(\beta)$ is stronger
than all $p^i_{\proj^{\tree_\lambda}_{\tree_i}(t_\alpha)}(\beta)$ (as these
are answers of the player $\inc$; see above). Lastly, if we let $Y_\alpha=
X_\beta$ then we have
\[\begin{array}{ll}
q_\alpha(\beta)\forces_{\nbQ^{G_\beta}_\beta}&\mbox{``the sequence }\langle
N_i[G_\beta][\name{G}_{\nbQ^{G_\beta}_\beta}]:i\leq\lambda\rangle\mbox{ is 
suitable and $(\bar{a},Y_\alpha)$}\\
\ &\ \mbox{ is a suitable base for it and $(\forall i\in Y_\alpha)(\forall\xi
\in\alpha\cap w_i)(i\in \name{X}_\xi^{G_\beta})$''}  
  \end{array}
\]
(compare the arguments in the proof of \ref{iterstrong}). This is everything
we need: as $G_\beta$ was any generic filter containing $q_\beta$ we may take
names $\name{t}_\alpha$, $\name{X}_\beta$, $\name{Y}_\alpha$ for the objects
defined above and the name for $q_\alpha(\beta)$ and conclude that
$q_\beta\conc q_\alpha(\beta)$ forces that they are as required. 
\medskip

\noindent{\sc Case 3:}\qquad $\alpha$ is a limit ordinal.\\
Arriving at this stage we have defined $q_\beta,\name{t}_\beta,\name{Y}_\beta$
and $\name{X}_\beta$ for $\beta\in\alpha\cap w_\lambda$ and we are going to
define $q_\alpha$, $\name{t}_\alpha$ and $\name{Y}_\alpha$. The first two
objects to be defined are determined by clauses {\bf (a)}--{\bf (d)}. The only
possible problem that may appear here is that we want $\name{t}_\alpha$ to be
(a name for) an element of $T_\lambda$ and thus of $\V$. But by
\ref{iterweak}(a) and \ref{nonewseq} + \ref{scomimpcom} we know that the
forcing with $\bP_\alpha$ adds no new sequences of length $<\kappa$ of
elements of $\V$ (remember $\kappa=\lambda^+$). Therefore the sequence
$\langle\name{t}_\beta: \beta\in w_\lambda\cap\alpha\rangle$ is a
$\bP_\alpha$--name for a sequence {\em from $\V$} and its limit
$\name{t}_\alpha$ is forced to be in $\tree_\lambda$. Now we immediately get
that $q_\alpha$, $\name{t}_\alpha$ satisfy demands {\bf (a)}--{\bf (f)} (for
{\bf (e)} note that $\dom(p^i_{\proj^{\tree_\lambda}_{\tree_i}(\name{t}_\alpha
)})\subseteq w_\lambda$ and
\begin{description}
\item[$(\boxtimes)$] for each $\beta\in\alpha\cap w_\lambda$ and $i<\lambda$
we have $\dom(p^i_{\proj^{\tree_\lambda}_{\tree_i}(\name{t}_\beta)})\subseteq
w_\lambda$ and  
\[p^i_{\proj^{\tree_\lambda}_{\tree_i}(\name{t}_\alpha)}\rest\rk_i(
\proj^{\tree_\lambda}_{\tree_i}(\name{t}_\beta))=p^i_{\proj^{\tree_\lambda}_{
\tree_i}(\name{t}_\beta)}\quad\mbox{ and }\quad\rk_i(\proj^{\tree_\lambda}_{
\tree_i}(\name{t}_\beta))\geq\beta,\]
\end{description}
hence we may use the clause {\bf (e)} from stages $\beta<\alpha$). Finally we
let 
\[\name{Y}_\alpha\stackrel{\rm def}{=}\{i<\lambda: i\mbox{ is limit and }
(\forall\xi\in w_i\cap\alpha)(i\in\name{X}_\xi)\mbox{ and }\alpha\in w_i\}.\]
We have to check that the demand {\bf (g)} is satisfied. Suppose that
$G_\alpha\subseteq\bP_\alpha$ is a generic filter over $\V$ containing
$q_\alpha$. The sequence $\langle \name{X}_\xi^{G_\alpha}: \xi\in
w_\lambda\cap \alpha\rangle$ is a sequence of length $<\kappa$ of elements of
$\V$, and the forcing with $\bP_\alpha$ adds no new such sequences. 
Consequently  
\[\langle \name{X}_\xi^{G_\alpha}: \xi\in w_\lambda\cap \alpha\rangle\in \V.\]
If for $j<\lambda$ we let $Z_j=\bigcap\limits_{\xi\in w_j\cap\alpha}
\name{X}_\xi^{G_\alpha}$ we will have 
\[\langle Z_j:j<\lambda\rangle\in\V,\quad\mbox{ and }\quad (\forall j<\lambda)
(Z_j\in D_{\bar{a}})\]
(as the filter $D_{\bar{a}}$ is $\lambda$--complete) and therefore (by the
normality of $D_{\bar{a}}$)
\[\name{Y}^{G_\alpha}_\alpha\supseteq\mathop{\bigtriangleup}\limits_{j<\lambda}
Z_j=\{i<\lambda:i\mbox{ is limit and }(\forall j<i)(i\in Z_j)\}\in
D_{\bar{a}}.\] 
Next note that $(\bar{a},\name{Y}_\alpha^{G_\alpha})$ is a suitable basis for
the sequence $\langle N_i[G_\alpha]: i\leq\lambda\rangle$. Why? Suppose that
$i\in\name{Y}_\alpha^{G_\alpha}$ and let $t=\proj^{\tree_\lambda}_{
\tree_{i+1}}(\name{t}_\alpha^{G_\alpha})$. By the choice of the $w_i$'s we
know that $w_{i+1}=w_i$ (remember $i$ is limit). Since $\alpha\in w_i$ we have
$\rk_{i+1}(t)=\alpha$ and since $i\in\bigcap\limits_{\xi\in w_i\cap\alpha}
\name{X}_\xi^{G_\alpha}$ we have $t_\xi(i)\neq *$ for each $\xi\in w_i\cap
\alpha=w_{i+1}\cap\alpha$. So look now at the way we defined $\bar{p}^{i+1}$:
we were at the case when $\bar{p}^{i+1}_t$ was given by
\ref{onesteptree}(c)(ii). In particular, the condition $p^{i+1}_t\in
\bP_{\rk_{i+1}(t)}'\cap N_{i+2}$ generates a $\bP_{\rk_{i+1}(t)}$--generic
filter over $N_{i+1}$. We know already that $q_\alpha,\name{t}_\alpha$ satisfy
{\bf (e)} (or use just $(\boxtimes)$) and therefore $p^{i+1}_t\in G_\alpha$. 
This is enough to conclude that
\[N_{i+1}[G_\alpha]\prec\Hchi,\quad (N_{i+1}[G_\alpha])^{<\lambda}\subseteq
N_{i+1}[G_\alpha],\quad N_{i+1}[G_\alpha]\cap\V=N_{i+1},\] 
(like in \ref{iterstrong}) and therefore to finish the construction.
\medskip

To finish the proof of this case of the theorem note that our demands on
conditions $q_\alpha$ imply that each of them is
$(N_\lambda,\bP_\alpha)$--generic, so in particular $q_\gamma$ is as
required.  \QED

\begin{theorem}
Assume $\lambda^{<\lambda}=\lambda$, $\kappa=\lambda^+=2^\lambda\leq\mu^*$.
Suppose that $\family_0\in(\dkl)^+$, $\hat{\family}_1\in (\dlk[\family_0])^+$.
Let $\bar{\bQ}=\langle\bP_i,\nbQ_i: i<\gamma\rangle$ be a $(<\kappa)$--support
iteration such that for each $i<\gamma$ 
\[\forces_{\bP_i}\mbox{`` }\nbQ_i\mbox{ is basically $(\family_0,
\hat{\family}_1)$--complete''}.\]
Then the forcing notion $\bP_\gamma$ is basically $(\family_0,
\hat{\family}_1)$--complete.
\end{theorem}

\Proof  Similar to the proof of \ref{iterweak} (but easier) and not used, so
we do not give details. \QED

\section{The Axiom}

\begin{definition}
\label{firstaxiom}
Suppose that $\lambda^{<\lambda}=\lambda$, $\kappa=\lambda^+=2^\lambda\leq
\mu^*$ and $\theta$ is a regular cardinal. Let $\family_0\in(\dkl)^+$,
$\hat{\family}_1\in (\dlk[\family_0])^+$ and let $D$ be a function from
$\hat{\family}_1$ such that each $D_{\bar{a}}$ is a normal filter on
$\lambda$. Let $\Axkt(\family_0,\hat{\family}_1)$, {\em the forcing axiom for
$(\family_0,\hat{\family}_1)$ and $\theta$}, be the following sentence:
\begin{quotation}
\noindent If $\bQ$ is a really $(\family_0,\hat{\family}_1,D)$--complete
forcing notion of size $\leq\kappa$ and $\langle {\cal I}_i: i<i^*<\theta
\rangle$ is a sequence of dense subsets of $\bQ$,

\noindent then there exist a directed set $H\subseteq\bQ$ such that
\[(\forall i<i^*)(H\cap {\cal I}_i\neq\emptyset).\]
\end{quotation}
\end{definition}

\begin{theorem}
\label{noproofaxiom}
Assume that $\lambda,\kappa=\mu^*,\theta$ and $(\family_0,\hat{\family}_1,D)$
are as in \ref{firstaxiom} and 
\[\kappa<\theta=\cf(\theta)\leq\mu=\mu^\kappa\]
(e.g.\ $(\circledast)$\quad $S_0\subseteq S^\lambda_\kappa$, $S_1=
S^\lambda_\kappa\setminus S_0$ are stationary subsets of $\kappa$,
$\family_0=S_0$, $\hat{\family}_1=\{\bar{a}: \bar{a}$ is an increasing
continuous sequence of ordinals, $a_0\in S_0$, $a_{i+1}\in S_0$, $a_\lambda\in
S_1\ \}$). 

\noindent Then there is a forcing notion $\bP$ of cardinality $\mu$ such that 
\begin{description}
\item[($\alpha$)] $\bP$ satisfies the $\kappa^+$--cc,
\item[($\beta$)]  $\forces_{\bP}\family_0\in(\dkl)^+\quad\&\quad \hat{
\family}_1\in (\dlk[\family_0])^+$ and even more:
\item[($\beta^+$)] if $\hat{\family}_1^*\subseteq\hat{\family}_1$ is
such that $\hat{\family}_1^*\in (\dlk[\family_0])^+$ 

then $\forces_{\bP}\hat{\family}_1^*\in (\dlk[\family_0])^+$,
\item[($\gamma$)] $\forces_{\bP}\Axkt(\baza_0,\baza_1)$,
\item[($\delta$)] if $(\circledast)$ then all stationary subsets of $\kappa$
are preserved. 
\end{description}
\end{theorem}

\Proof It is parallel to \ref{getaxiom} which is later done elaborately. \QED

\PART{B}{Case}

\noindent While {\bf Case D} (see the introduction; $\kappa$ inaccessible, $S$
has stationary many inaccessible members) may be treated similarly to {\bf
Case A}, we need to refine our machinery to deal with {\bf Case B}. Our
prototype here is {\em $\kappa$ is the first strongly inaccessible cardinal},
however the tools developed in this part will be applicable to cases {\bf A},
{\bf C}, {\bf D} too (and other strong inaccessibles in {\bf Case B}, of
course).   

\begin{ourassum}
$\kappa$ is a strongly inaccessible cardinal and $\mu^*\geq\kappa$ is a
regular cardinal.\\
These assumptions will be kept in the present part (unless otherwise stated)
and we may forget to remind them.
\end{ourassum}

There are two main difficulties which one meets when dealing with the present
case. First problem, a more general one, is that $(<\mu)$--completeness is
not reasonable even for $\mu=\aleph_1$. Why? As we would like to force the
Uniformization Property for $\langle S_\delta: \delta\in S\rangle$, where
$S\subseteq\{\delta<\kappa:\cf(\delta)=\aleph_0\}$ is stationary not
reflecting. The second problem is related to closure properties of models we
consider. In {\bf Case A}, when $\kappa=\lambda^+$, the demand $N^{<\lambda}
\subseteq N$ was reasonable. If $\kappa$ is Mahlo, $\|N\|=N\cap\kappa$ is an
inaccessible cardinal $<\kappa$ then the demand $N^{<N\cap\kappa}\subseteq N$
is reasonable too. However, if $\kappa$ is the first inaccessible this does
not work. (Note that these models are parallel of {\em countable
$N\prec\Hchi$} of the case $\kappa=\aleph_1$.)  To handle these problems we
will use exclusively sequences $\bar{N}=\langle N_i: i\leq\alpha\rangle$ of
models and all action will take place at limit stages only. For example, we
will have completeness for $\bar{N}=\langle N_i: i\leq\omega\rangle$ by
looking at $N_\omega$, BUT the equivalence class $\bar{N}/\approx$ will be
important too, where for two sequences $\bar{N}$, $\bar{N}'$ of length
$\omega$ we write $\bar{N}\approx\bar{N}'$ if
\[(\forall n\in\omega)(\exists m\in\omega)(N_n\subseteq N_m')\quad\mbox{ and
}\quad (\forall n\in\omega)(\exists m\in\omega)(N_n'\subseteq N_m)\]

\section{More on complete forcing notions}
In this section we introduce more notions of completeness of forcing
notions. In some sense we will generalize and develop the notions introduced
in section \ref{caseA}.

\begin{definition}
\label{ckibaza}
\begin{enumerate}
\item Let $\bar{N}=\langle N_i: i\leq\alpha\rangle$ be a sequence of models and
$\bar{a}=\langle a_i: i\leq\alpha\rangle$ be a sequence of elements of
$[\mu^*]^{<\kappa}$. We say that {\em $\bar{N}$ obeys $\bar{a}$ with an error
$n\in\omega$} if
\[(\forall i<\alpha)(a_i\subseteq N_i\cap\mu^*\subseteq a_{i+n}).\]
When we say {\em $\bar{N}$ obeys $\bar{a}$} we mean  {\em with some error
$n\in\omega$}. 
\item By $\ck$ we will denote the collection of all sets $\baza$ such that 
\[\begin{array}{r}
\baza\subseteq\big\{\bar{a}=\langle a_i: i\leq\alpha\rangle: \mbox{
the sequence $\bar{a}$ is increasing continuous, }\ \\
\alpha<\kappa\ \mbox{ and }\ (\forall i\leq\alpha)(a_i\in [\mu^*]^{<\kappa}
\ \&\ a_i\cap\kappa\in\kappa)\big\},
\end{array}\]
and for every regular large enough cardinal $\chi$, for every $x\in {\cal H}(
\chi)$ and a regular cardinal $\theta<\kappa$ there are $\bar{N}=\langle N_i:
i\leq\theta\rangle$ and $\bar{a}=\langle a_i: i\leq\theta\rangle$ such that
\begin{description}
\item[(a)] $\bar{N}$ is an increasing continuous sequence of elementary
submodels of $\Hchi$ such that $(\forall i<\theta)(\bar{N}\rest (i+1)\in
N_{i+1}\quad\&\quad \|N_i\|<\kappa)$ and  $x\in N_0$, 
\item[(b)] $\bar{a}\in \baza$,
\item[(c)] $\bar{N}$ obeys $\bar{a}$.
\end{description}
\item If $\bar{a}\in\baza$, $\bar{N}$ is an increasing continuous sequence of
elementary submodels of $\Hchi$ such that $(\forall i+1<\lh(\bar{N}))(\bar{N}
\rest (i+1)\in N_{i+1}\ \&\ \|N_i\|<\kappa)$ and $\bar{N}$ obeys $\bar{a}$
(with error $n$, respectively) then we say that $(\bar{N},\bar{a})$ is {\em an
$\baza$--complementary pair} ({\em an $(\baza,n)$--complementary pair},
respectively).  
\item We say that a family $\baza\in\ck$ is {\em closed} if for every sequence 
$\bar{a}=\langle a_i: i\leq\alpha\rangle\in\baza$ and ordinals $\beta,\gamma$
such that $\beta+\gamma\leq\alpha$ we have
\[\langle a_{\beta+i}:i\leq\gamma\rangle\in\baza\]
(or, in other words, $\baza$ is closed under both initial and end segments).
\end{enumerate}
\end{definition}

\begin{remark}
{\em
\begin{enumerate}
\item Definition \ref{ckibaza} is from \cite[\S 1]{Sh:186}.
\item The exact value of the error $n$ in \ref{ckibaza}(2) is not important at
all, we may consider here several other variants as well.
\item Note that $N_i,\|N_i\|\in N_{i+1}$. Sometimes we may add to
\ref{ckibaza}(1) a requirement that $2^{\|N_i\|}\subseteq a_{i+n}$ (saying
that then {\em $\bar{N}$ strongly obeys $\bar{a}$}). Note that this naturally
occurs to be true for strongly inaccessible $\kappa$, as we demand that \quad
$\bar{a}\in\baza\quad\Rightarrow\quad a_i\cap\kappa\in\kappa$.\quad So then
$2^{\|N_i\|}\in a_{i+n}$, but $a_{i+n}\cap\kappa\in\kappa$ so we have
$2^{\|N_i\|}\subseteq a_{i+n}$. 

In this situation, {\em if\/} $\chi_1<\chi$ are large enough, $\chi_1\in N_0$
and for non-limit $i$, $N_i'$ is the closure of $N_i\cap{\cal H}(\chi_1)$
under Skolem functions and sequences of length $\leq\|N_i\|$, and for limit
$i$, $N_i'=N_i\cap{\cal H}(\chi_1)$ {\em then\/} the sequence $\langle N_i':
i\leq\alpha\rangle$ will have closure properties and will obey $\bar{a}$ (as
$N_i\in N_{i+1}$, ${\cal H}(\chi_1)\in N_{i+1}$ imply $N_i'\in N_{i+1}$ and so
$N_i'\subseteq N_{i+n}$).

\item The presence of ``regular $\theta<\kappa$'' in \ref{ckibaza}(2) is not
accidental; it will be of special interest when $\kappa$ is a successor of a
singular strong limit cardinal, as then $\theta=\cf(\theta)<\kappa=\mu^+$
implies $\theta<\mu$. 
\end{enumerate}
}
\end{remark}

\begin{definition}
\label{Ecomplete}
Let $\baza\in\ck$ and let $\bQ$ be a forcing notion.
\begin{enumerate}
\item Let $\bar{N}=\langle N_i: i\leq \delta\rangle$ be an increasing
continuous sequence of elementary submodels of $\Hchi$, $\bQ\in N_0$ and
$\bar{p}=\langle p_i: i<\delta\rangle$ be an increasing sequence of conditions
from $\bQ\cap N_\delta$, $n\in\omega$. We say that {\em $\bar{p}$ is
$(\bar{N},\bQ)^n$--generic} if for each $i<\delta$  
\[\bar{p}\rest (i+1)\in N_{i+1}\ \mbox{ and }\ p_{i+n}\in\bigcap\{{\cal I}\in
N_i: {\cal I}\mbox{ is an open dense subset of }\bQ\}.\]
When we say that {\em $\bar{p}$ is $(\bar{N},\bQ)^*$--generic} we mean that it
is $(\bar{N},\bQ)^n$--generic for some $n\in\omega$. We may say then that
$\bar{p}$ is {$(\bar{N},\bQ)^*$--generic with an error $n$}.
\item We say that {\em $\bQ$ is complete for $\baza$} if for large enough
$\chi$, for some $x\in {\cal H}(\chi)$ the following condition is satisfied:
\begin{description}
\item[$(\circledast)^{\baza}_x$] {\em if} 
\begin{description}
\item[(a)] $(\bar{N},\bar{a})$ is an $\baza$--complementary pair (see
\ref{ckibaza}(3)), $\bar{a}\in\baza$, $\bar{N}=\langle N_i: i\leq\delta
\rangle$, $\bQ,x\in N_0$, and
\item[(b)] $\bar{p}$ is an increasing $(\bar{N},\bQ)^*$--generic sequence
\end{description}
{\em then} $\bar{p}$ has an upper bound in $\bQ$.
\end{description}
\item We say that a forcing notion $\bQ$ is {\em strongly complete for
$\baza$} if it is complete for $\baza$ and does not add sequences of ordinals
of length $<\kappa$. 
\end{enumerate}
\end{definition}

\begin{remark}
{\em
\begin{enumerate}
\item The $x$ in definition \ref{Ecomplete}(2) is the way to say ``for most'',
compare \ref{most}. 
\item In the present applications, we will have $\mu^*=\kappa$ and a
stationary set $S\subseteq\kappa$ such that 
\[\begin{array}{r}
\baza^c_S\stackrel{\rm def}{=}\big\{\bar{a}: \bar{a}\mbox{ an increasing
sequence of ordinals from }\kappa\setminus S\ \ \\
\mbox{of length $<\kappa$ with the last element from }S\big\}
  \end{array}\]
will be in $\ck$. The forcing notions will be complete for $\baza^c_S$, so the
iteration will add no sequences of length $<\kappa$ (see \ref{nextiter}
below). On $S$ the behavior will be more interesting doing the
uniformization. Thus the pair $(\baza^c_S,S)$ corresponds to the pair
$(\family_0,\hat{\family}_1)$ from the previous part (on Case A). 

For example, if $C_\delta\subseteq\delta=\sup(C_\delta)$, $\otp(C_\delta)=
\cf(\delta)$, $(\forall \delta\in S)(\cf(\delta)<\delta)$ and $h_\delta:
C_\delta\longrightarrow 2$ then 
\[\begin{array}{ll}
\bQ=\{g:&\mbox{for some }\alpha<\kappa,\ \ g:\alpha\longrightarrow 2\ \ \mbox{
and}\\
\ 	&(\forall\delta\in (\alpha+1)\cap S)(\forall\gamma\in C_\delta\mbox{
large enough})(g(\gamma)=h_\delta(\gamma))\big\}
  \end{array}\]
is such a forcing (but we need that $S$ is not reflecting or $\langle
C_\delta: \delta\in S\rangle$ is somewhat free, so that for each
$\alpha<\kappa$ there are $g\in\bQ$ with $\dom(g)=\alpha$). 
\item If we want to have $S$ reflecting in a stationary set though still
``thin'' then things are somewhat more complicated, but manageable, see later.
\end{enumerate}
}
\end{remark}

\begin{proposition}
\label{bazapres}
Suppose that $\baza\in\ck$ is closed and $\bQ$ is a forcing notion.
\begin{enumerate}
\item Assume $(\bar{a},\bar{N})$ is an $(\baza,n_1)$--complementary pair,
$\bar{a}\in\baza$, $\bar{N}=\langle N_i: i\leq\delta\rangle$, $\bQ\in N_0$. If 
$\bar{p}\subseteq\bQ\cap N_\delta$ is $(\bar{N},\bQ)^{n_2}$--generic (see
\ref{Ecomplete}(1)) and $q\in\bQ$ is an upper bound of $\bar{p}$ in $\bQ$ then
\[q\forces_{\bQ}\mbox{``}(\langle N_i[\name{G}_{\bQ}]:i\leq\delta\rangle,
\bar{a})\mbox{ is an $(\baza,n_1+n_2+1)$--complementary pair''}.\]
\item If $\bQ$ is strongly complete for $\baza$ then
$\forces_{\bQ}\baza\in\ck$.  
\end{enumerate}
\end{proposition}

\Proof 1)\ \ \ Since $\bar{p}$ is $(\bar{N},\bQ)^{n_2}$--generic, for each
$i<\delta$ and every $\bQ$--name $\name{\tau}\in N_i$ for an element of $\V$,
the condition $p_{i+n_2}$ decides the value of $\name{\tau}$ and the decision
belongs to $N_{i+n_2+1}$ (remember $p_{i+n_2}\in N_{i+n_2+1}$). Now, by
standard arguments (like in the proofs of \ref{cl1} and \ref{cl2}) we conclude
that for each $i<\delta$ 
\[\begin{array}{ll}
p_{i+n_2+1}\forces_{\bQ}&\mbox{``}N_i[\name{G}_{\bQ}]\cap\V\subseteq N_{i+n_2
+1}\quad\mbox{ and }\quad N_i[\name{G}_{\bQ}]\prec\Hchi\quad\mbox{ and}\\
\ &\ \langle N_j[\name{G}_{\bQ}]\!:j\leq i\rangle\in N_{i+1}[\name{G}_{\bQ}]
\mbox{''.}
  \end{array}\]
Since $a_{i+n_2+1}\subseteq N_{i+n_2+1}\subseteq a_{i+n_2+1+n_1}$ (for
$i<\delta$) we get
\[q\forces_{\bQ}\mbox{``}(\langle N_i[\name{G}_{\bQ}]:i\leq\delta\rangle,
\bar{a})\mbox{ is an $(\baza,n_1+n_2+1)$--complementary pair''}.\]

\noindent 2)\ \ \ Suppose that $p\forces_{\bQ}\name{x}\in {\cal H}(\chi)$ and
let $\theta<\kappa$ be a regular cardinal. Since $\baza\in\ck$ we find an
$(\baza,n_1)$--complementary pair $(\bar{N},\bar{a})$ such that $\lh(\bar{N})=
\lh(\bar{a})=\theta+1$ and $(p,\name{x},\bQ,\baza)\in N_0$. Now, by induction
on $i<\theta$, we define an $(\bar{N},\bQ)^1$--generic sequence $\bar{p}=
\langle p_i: i<\theta\rangle$:
\begin{quotation}
\noindent $p_i\in N_{i+1}\cap\bQ$ is the $<^*_\chi$--first element $q$ of
$\bQ$ such that 
\begin{description}
\item[(i)$_i$] \ $p\leq q$ and $(\forall j<i)(p_j\leq q)$,
\item[(ii)$_i$]  $q\in\bigcap\{{\cal I}\in N_i: {\cal I}\subseteq\bQ$ is open
dense$\}$. 
\end{description}
\end{quotation}
To show that this definition is correct we have to prove that, for each
$i<\theta$, there is a condition $q\in\bQ$ satisfying {\bf (i)$_i$}+{\bf
(ii)$_i$} and $\bar{p}\rest i\in N_{i+1}$. Note that once we know this, we
are sure that the $<^*_\chi$--first condition with these properties is in
$N_{i+1}$ and therefore $\bar{p}\rest (i+1)\in N_{i+1}$ too.\\
There are no problems for $i=0$, so suppose that $i=i_0+1$ and we have already
$\bar{p}\rest i_0\in N_{i_0+1}$ and $p_{i_0}\in N_{i_0+1}$ and we defined
hence $\bar{p}\rest (i_0+1)\in N_{i_0+1}\prec N_{i_0+2}$. The forcing notion
$\bQ$ does not add new sequences of ordinals of length $<\kappa$ and $\|N_{
i_0+1}\|<\kappa$. Therefore we find a condition $q\in\bQ$ stronger than
$p_{i_0}$ and such that $q$ decides all $\bQ$--names for ordinals from
$N_{i_0+1}$ (i.e. $q\in\bigcap\{{\cal I}\in N_i: {\cal I}\subseteq\bQ$ is open
dense$\}$).\\ 
Suppose now that we arrive to a limit stage $i$ and we have defined $\bar{p}
\rest i$. Since $\langle N_j: j\leq i\rangle\in N_{i+1}$ we know that $\bar{p}
\rest i\in N_{i+1}$ (as all the parameters needed for the definition of
$\bar{p}\rest i$ are in $N_{i+1}$ and we have left no freedom). Note that
$\bar{a}\rest (i+1)\in\baza$ (as $\baza$ is closed), $(\bar{a}\rest (i+1),
\bar{N}\rest (i+1))$ is an $(\baza,n_1)$--complementary pair and the sequence
$\bar{p}\rest i$ is $(\bar{N}\rest (i+1),\bQ)^1$--generic. Since $\bQ$ is
strongly complete for $\baza$ we conclude that there is an upper bound to
$\bar{p}\rest i$ in $\bQ$. Now it should be clear that such an upper bound
satisfies {\bf (i)$_i$}+{\bf (ii)$_i$} (remember that $\bar{N}$ is increasing
continuous).  

Now look at the sequence $\bar{p}=\langle p_i: i<\theta\rangle$. Immediately
by its definition we see that $\bar{p}$ is $(\bar{N}\rest
(i+1),\bQ)^1$--generic. Since $\bQ$ is strongly complete for $\baza$ we find
an upper bound $q\in\bQ$ of $\bar{p}$. Now, by the first part of the
proposition, we conclude that  
\[q\forces_{\bQ}\mbox{``}(\langle N_i[\name{G}_{\bQ}]:i\leq\delta\rangle,
\bar{a})\mbox{ is an $(\baza,n_1+2)$--complementary pair''},\]
what finishes the proof. \QED

\begin{theorem}
\label{nextiter}
Suppose that $\baza\in\ck$ is closed and $\langle\bP_i,\nbQ_i:i<\gamma\rangle$
is a $(<\kappa)$--support iteration such that for each $i<\gamma$ 
\[\forces_{\bP_i}\mbox{``the forcing notion $\nbQ_i$ is strongly complete for
$\baza$''.}\]
Then  $\bP_\gamma$ is strongly complete for $\baza$.
\end{theorem}

\Proof We prove the theorem by induction on $\gamma$.

\noindent{\sc Case 1:}\ \ \ $\gamma=0$.\\
There is nothing to do in this case.

\noindent{\sc Case 2:}\ \ \ $\gamma=\beta+1$.\\
By the induction hypothesis we know that $\bP_\beta$ is strongly complete for
$\baza$ and therefore, by \ref{bazapres}, $\forces_{\bP_\beta}\baza\in\ck$.\\
Clearly the composition of two forcing notions not adding new sequences of
length $<\kappa$ of ordinals does not add such sequences. Thus what we have to
prove is that $\bP_{\beta+1}=\bP_\beta*\nbQ_\beta$ is complete for $\baza$
(i.e. \ref{Ecomplete}(2)).\\
Let $y\in {\cal H}(\chi)$ be the witness for ``$\bP_\beta$ is complete for
$\baza$'' and let $\name{x}$ be a $\bP_\beta$--name for the witness for
``$\nbQ_\beta$ is complete for $\baza$''. We are going to show that the
composition $\bP_{\beta+1}=\bP_\beta*\nbQ_\beta$ satisfies the condition
$(\circledast)^{\baza}_{\langle y,\name{x},\baza,\bP_{\beta+1}\rangle}$ of
\ref{Ecomplete}(2). So suppose that
\begin{description}
\item[(a)] $(\bar{N},\bar{a})$ is an $\baza$--complementary pair (with an
error, say, $n_1$), $y,\name{x},\baza,\bP_{\beta+1}\in N_0$, $\lh(\bar{N})=
\lh(\bar{a})=\delta+1$,
\item[(b)] $\bar{p}=\langle p_i: i<\delta\rangle$ is an increasing $(\bar{N},
\bP_{\beta+1})^{n_2}$--generic sequence.
\end{description}
It should be clear that the sequence $\langle p_i\rest\beta: i<\delta\rangle$
is $(\bar{N},\bP_\beta)^{n_2}$--generic. Therefore, as $\bP_\beta$ is complete
for $\baza$ and $y\in N_0$, we find a condition $q^*\in\bP_\beta$ stronger
than all $p_i\rest\beta$ (for $i<\delta$). By \ref{bazapres}(1) we know that 
\[q^*\forces_{\bP_\beta}\mbox{``}(\langle N_i[\name{G}_{\bP_\beta}]:i\leq
\delta\rangle,\bar{a})\mbox{ is an $(\baza,n_1+n_2+1)$--complementary 
pair''}.\]
Moreover
\[q^*\forces_{\bP_\beta}\mbox{``}\langle p_i(\beta)\!:\beta\!<\!\delta\rangle
\mbox{ is an increasing }(\langle N_i[\name{G}_{\bP_\beta}]\!:i\!\leq\!\delta
\rangle,\nbQ_\beta)^{n_2}\mbox{-generic sequence''.}\]
[Why? Like in \ref{cl3}, if $\name{\cal I}\in N_i$ is a $\bP_\beta$--name for
an open dense subset of $\nbQ_\beta$ then the set
\[\{p\in\bP_{\beta+1}: p\rest\beta\forces_{\bP_\beta} p(\beta)\in\name{\cal
I}\} \in N_i\]
is open dense in $\bP_{\beta+1}$; now use the choice of $q^*$.] Consequently,
we find a $\bP_\beta$--name $\name{\tau}$ for an element of $\nbQ_\beta$ such
that 
\[q^*\forces_{\bP_\beta}\mbox{``}(\forall i<\delta)(p_i(\beta)
\leq_{\nbQ_\beta}\name{\tau})\mbox{''.}\]
Let $q=q^*\cup\{(\beta,\name{\tau})\}$. Clearly $q\in\bP_{\beta+1}$ is an
upper bound of $\bar{p}$.

\noindent{\sc Case 3:}\ \ \ $\gamma$ is a limit ordinal.\\
Let $\name{x}_\beta$ (for $\beta<\gamma$) be a $\bP_\beta$--name for the
witness for $\forces_{\bP_\beta}$``$\nbQ_\beta$ is complete for $\baza$''. Let
$x=\langle\langle\name{x}_\beta:\beta<\gamma\rangle, \langle\bP_\beta,
\nbQ_\beta: \beta<\gamma\rangle\rangle$. 

\begin{claim}
\label{cl6}
Suppose that $(\bar{N},\bar{a})$ is an $\baza$--complementary pair,
$\lh(\bar{N})=\lh(\bar{a})=\delta+1$, $\delta$ is a limit ordinal and $x\in
N_0$. Further assume that $\bar{p}=\langle p_i:i<\delta\rangle\subseteq
\bP_\gamma$ is an increasing sequence of conditions from $\bP_\gamma$ such that
\begin{description}
\item[(a)] $(\forall i<\delta)(\bar{p}\rest (i+1)\in N_{i+1})$, \quad and
\item[(b)] for every $\beta\in\gamma\cap N_\delta$ there are $n<\omega$ and
$i_0<\delta$ such that
\[(\forall i\in [i_0,\delta))(p_{i+n}\rest\beta\in\bigcap\{{\cal I}\in N_i:
{\cal I}\mbox{ is an open dense subset of }\bP_\beta\}).\]
\end{description}
Then the sequence $\bar{p}$ has an upper bound in $\bP_\gamma$.\\
\relax [Note: we do not put any requirements on meeting dense subsets of
$\bP_\gamma$!] 
\end{claim}

\noindent{\em Proof of the claim:}\hspace{0.15in} We define a condition $q\in
\bP_\gamma$. First we declare that $\dom(q)=N_\delta\cap\gamma$ and next we
choose $q(\beta)$ by induction on $\beta\in N_\delta\cap\gamma$ in such a way
that $(\forall i<\delta)(p_i\rest\beta\leq_{\bP_\beta}q\rest\beta)$. So
suppose that we have defined $q\rest\beta\in\bP_\beta$, $\beta\in\gamma\cap
N_\delta$. Let $n\in\omega$ and $i_0<\delta$ be given by the assumption {\bf
(b)} of the claim for $\beta+1$. We may additionally demand that $\beta\in
N_{i_0}$. (Note that $n,i_0$ are good for $\beta$ too, remember $\bP_\beta
\lesdot\bP_{\beta+1}$.) Since $\baza$ is closed we know that $(\bar{N}\rest
[i_0,\delta],\bar{a}\rest [i_0,\delta])$ is an $\baza$--complementary pair and
the sequence $\langle p_i\rest\beta: i_0\leq i<\delta\rangle$ is $(\bar{N}
\rest [i_0,\delta],\bP_\beta)^n$--generic. Consequently, by \ref{bazapres}(1),
we get 
\[q\rest\beta\forces_{\bP_\beta}\mbox{``}(\bar{N}[\name{G}_{\bP_\beta}]\rest
[i_0,\delta],\bar{a}\rest [i_0,\delta])\mbox{ is an $\baza$--complementary
pair''.}\]
Moreover, like in the previous case, the condition $q\rest\beta$ forces (in
$\bP_\beta$) that 
\[\mbox{``}\langle p_i(\beta)\!:i_0\leq i<\delta\rangle\mbox{ is an increasing
$(\bar{N}[\name{G}_{\bP_\beta}]\rest [i_0,\delta],\nbQ_\beta)^n$-generic
sequence''.}\] 
Thus, as $\name{x}_\beta\in N_{i_0}$ and $\nbQ_\beta$ is a name for a forcing
notion which is complete for $\baza$ with the witness $\name{x}_\beta$, we
find a $\bP_\beta$--name $q(\beta)$ such that
\[q\rest\beta\forces_{\bP_\beta}\mbox{``}(\forall i<\delta)(p_i(\beta)
\leq_{\nbQ_\beta} q(\beta))\mbox{''.}\]
Now we finish the proof of the claim notifying that if $\beta\in\gamma\cap
N_\delta$ is limit and for each $\alpha\in\beta\cap N_\delta$, $q\rest\alpha$
is an upper bound to $\langle p_i\rest\alpha:i<\delta\rangle$ then $q\rest
\beta$ is an upper bound of $\langle p_i\rest\beta:i<\delta\rangle$ (remember
$\dom(p_i)\subseteq N_\delta$ for each $i<\delta$). 

\begin{claim}
\label{cl7}
Suppose that $M\prec\Hchi$, $\|M\|<\kappa$ and $p\in\bP_\gamma$. Then there is
a condition $q\in\bP_\gamma$ stronger than $p$ and such that
\[(\forall\beta\in M\cap\gamma)(q\rest\beta\in\bigcap\{{\cal I}\in M: {\cal I}
\mbox{ is an open dense subset of }\bP_\beta\}).\]
\end{claim}

\noindent{\em Proof of the claim:}\hspace{0.15in} Let $\theta=\cf(\otp(M\cap
\gamma))$ and let $\langle\gamma_i:i\leq\theta\rangle$ be an increasing
continuous sequence such that $\gamma_0=0$, $\gamma_\theta=\sup(M\cap\gamma)$ 
and $\gamma_i\in M\cap\gamma$ (for non-limit $i<\theta$). As $\baza\in\ck$, we
find $\bar{N}=\langle N_i:i\leq\theta\rangle$ and $\bar{a}=\langle a_i:i\leq
\theta\rangle\in\baza$ such that $\langle\gamma_i: i\leq\theta\rangle,x,p\in
N_0$ and $(\bar{N},\bar{a})$ is an $\baza$--complementary pair and $M\subseteq
N_0$. The last demand may seem to be too strong, but we use the fact that
$\baza$ is closed and
\[M\in N'\prec N''\prec \Hchi\ \&\ \sup(N'\cap\kappa)\subseteq N''\quad
\Rightarrow\quad M\subseteq N'.\]
(Alternatively, first we take an $\baza$--complementary pair $(\bar{N}^*,
\bar{a}^*)$ such that $\lh(\bar{N})=\lh(\bar{a})=\|M\|^++1$ and $\langle
\gamma_i: i\leq\theta\rangle,x,p,M\in N_0^*$. Next look at the model $N^*_{
\|M\|+1}$ -- it contains all ordinals below $\|M\|$, $M$ and $\|M\|$. Hence
$M\subseteq N^*_{\|M\|+1}$. Take $\bar{N}=\bar{N}^*\rest [\|M\|+1,\|M\|+
\theta]$ and $\bar{a}=\bar{a}^*\rest [\|M\|+1,\|M\|+\theta]$.)

Next, by induction on $i\leq\theta$, we define a sequence $\langle p_i:i\leq
\theta\rangle\subseteq\bP_\gamma$: 
\begin{quotation}
\noindent $p_i\in\bP_\gamma$ is the $<^*_\chi$--first element $q$ of
$\bP_\gamma$ such that  
\begin{description}
\item[(i)$_i$] \ \ $p\rest\gamma_i\leq_{\bP_{\gamma_i}} q\rest\gamma_i$ and
$(\forall j<i)(p_j\rest\gamma_i\leq_{\bP_{\gamma_i}} q\rest\gamma_i)$, 
\item[(ii)$_i$]  \ $q\rest\gamma_i\in\bigcap\{{\cal I}\in N_i: {\cal I}
\subseteq\bP_{\gamma_i}$ is open dense$\}$,
\item[(iii)$_i$]   $q\rest [\gamma_i,\gamma)=p\rest [\gamma_i,\gamma)$. 
\end{description}
\end{quotation}
We have to show that this definition is correct and for this we prove by
induction on $i\leq\theta$ that there is a condition $q\in\bP_{\gamma_i}$
satisfying {\bf (i)$_i$}--{\bf (iii)$_i$} and $\bar{p}\rest i\in N_{i+1}$. By
the way $p_i$'s are defined we will have that then $\bar{p}\rest (i+1)\in
N_{i+1}$ for $i<\theta$.\\ 
If $i$ is not limit (and we have $p_j$ for $j<i$) then there is no problem in
finding the respective condition $q$ once one realizes that, by the inductive
hypothesis of the theorem, the forcing notion $\bP_{\gamma_i}$ does not add
new sequences of length $<\kappa$ of ordinals and $\|N_i\|<\kappa$. So we just
pick up a condition in $\bP_{\gamma_i}$ stronger than the (respective
restriction of the) previous condition (if there is any) and which decides all
names for ordinals from $N_i$. This takes care of {\bf (i)$_i$} and {\bf
(ii)$_i$}. Next we extend our condition to a condition in $\bP_\gamma$ as the
requirement {\bf (iii)$_i$} says. Arriving to a limit stage $i$ we use Claim
\ref{cl6}.  So we have defined $\bar{p}\rest i$ and by the way it was defined
we know that $\bar{p}\rest i\in N_{i+1}$ (as all parameters are there). Since
$\baza$ is closed we know that $(\bar{N}\rest (i+1),\bar{a}\rest (i+1))$ is an
$\baza$--complementary pair. Now apply \ref{cl6} to $\gamma_i$,
$\bP_{\gamma_i}$, $\bar{p}\rest i$, $\bar{N}\rest (i+1)$ and $\bar{a}\rest
(i+1)$ in place of $\gamma$, $\bP_\gamma$, $\bar{p}$, $\bar{N}$, and $\bar{a}$
there. Note that the assumptions are satisfied: for {\bf (b)} use the fact that
$i$ is limit, so if $\beta<\gamma_i$ then for some $j<i$ we have $\beta<
\gamma_j$ and now this $j$ works as $i_0$ there with $n=1$. Consequently the
sequence $\bar{p}\rest i$ has an upper bound in $\bP_{\gamma_i}$. Now,
similarly as in the non-limit case, we may find a condition $q\in\bP_\gamma$
(stronger than this upper bound) satisfying {\bf (i)$_i$}--{\bf (iii)$_i$}.

Now look at the condition $p_\theta\in\bP_\gamma$. If $\beta\in M\cap\gamma$
and $i<\theta$ is such that $\beta<\gamma_i$ then $p_i\rest\gamma_i$ decides
all $\bP_{\gamma_i}$--names from $N_i$ for ordinals. But $M\subseteq N_0$,
$p_i\rest\gamma_i\leq_{\bP_{\gamma_i}} p_\theta\rest\gamma_i$ and $\bP_\beta
\lesdot\bP_{\gamma_i}$. Hence $p_\theta\rest\beta\in\bigcap\{{\cal I}\in M:
{\cal I}\subseteq\bP_\beta$ is open dense$\}$. As $p_\theta$ is stronger than
$p$, this finishes the proof of the claim.

\begin{claim}
\label{cl8}
$\bP_\gamma$ is complete for $\baza$
\end{claim}

\noindent{\em Proof of the claim:}\hspace{0.15in} We are going to show that
$\bP_\gamma$ satisfies the condition $(\circledast)^{\baza}_{\langle x,\baza
\rangle}$ of \ref{Ecomplete}(2). So suppose that $(\bar{N},\bar{a})$ is an
$\baza$--complementary pair, $\bar{N}=\langle N_i: i\leq\delta\rangle$,
$x,\baza\in N_0$ and $\bar{p}=\langle p_i:i<\delta\rangle$ is an increasing
$(\bar{N},\bP_\gamma)^{n_1}$--generic sequence. For $i<\delta$ let
\[{\cal I}^*_i\stackrel{\rm def}{=}\big\{q\!\in\!\bP_\gamma\!:
(\forall\beta\!\in \! N_i\cap\gamma)(q\rest\beta\in\bigcap\{{\cal I}\in N_i\!:
{\cal I}\subseteq \bP_\beta\mbox{ is open dense in }\bP_\beta\})\big\}.\]
Note that Claim \ref{cl7} says that each ${\cal I}^*_i$ is an open dense subset
of $\bP_\gamma$. Clearly ${\cal I}^*_i$ is in $N_{i+1}$, as it is defined from
$N_i$. Hence, for each $i<\delta$, $p_{i+1+n_1}\in {\cal I}^*_i$. Now look at
the assumptions of Claim \ref{cl6}: both {\bf (a)} and {\bf (b)} there are
satisfied (for the second note that if $\beta\in N_\delta\cap\gamma$ then we
may take $i_0<\delta$ so large that $\beta\in N_{i_0}$ and let $n=n_1+1$). 
Thus we may conclude that $\bar{p}$ has an upper bound in $\bP_\gamma$. 

\begin{claim}
\label{cl9}
Forcing with $\bP_\gamma$ does not add new sequences of length $<\kappa$ of
ordinals. 
\end{claim}

\noindent{\em Proof of the claim:}\hspace{0.15in} First note that for a
forcing notion $\bP$, ``not adding new sequences of length $\theta$ of
ordinals'' is equivalent to ``not adding new sequences of length $\theta$ of
elements of $\V$''. Next note that, for a forcing notion $\bP$, if $\theta$ is
the first ordinal such that for some $\bP$--name $\name{\tau}$ and a condition
$p\in\bP$ we have 
\[p\forces_{\bP}\mbox{``}\name{\tau}:\theta\longrightarrow\V\quad\mbox{and}\quad
\name{\tau}\notin\V\mbox{''}\]
then $\cf(\theta)=\theta$. [Why? Clearly such a $\theta$ has to be limit; if
$\cf(\theta)<\theta$ then take an increasing cofinal in $\theta$ sequence
$\langle\zeta_i: i<\cf(\theta)\rangle$ and look at $\langle\name{\tau}\rest
\zeta_i: i<\cf(\theta)\rangle$. Each $\name{\tau}\rest\zeta_i$ is forced to be
in $\V$, so the sequence of them is in $\V$ too -- a contradiction.]
Consequently it is enough to prove that for every regular cardinal
$\theta<\kappa$, forcing with $\bP_\gamma$ does not add new sequences of
length $\theta$ of elements of $\V$. So suppose that, for $i<\theta$,
$\name{\tau}_i$ is a $\bP_\gamma$--name for an element of $\V$, and 
$p\in\bP_\gamma$. Take an $\baza$--complementary pair $(\bar{N},\bar{a})$ such
that $\bar{N}=\langle N_i: i\leq\theta\rangle$ and $x,p,\langle\name{\tau}_i:
i<\theta\rangle\in N_0$ (exists as $\baza\in\ck$). Now, by induction on
$i\leq\theta$, define a sequence $\langle p_i: i\leq\theta\rangle\subseteq
\bP_\gamma$:
\begin{quotation}
\noindent $p_i\in\bP_\gamma$ is the $<^*_\chi$--first element $q$ of
$\bP_\gamma$ such that  
\begin{description}
\item[(i)$_i$] \ \ $p\leq_{\bP_\gamma} q$ and $(\forall j<i)(p_j\leq_{
\bP_\gamma} q)$, 
\item[(ii)$_i$]  \ if $\beta\in N_i\cap\gamma$ then $q\rest\beta\in\bigcap
\{{\cal I}\in N_i: {\cal I}\subseteq\bP_\beta$ is open dense$\}$,
\item[(iii)$_i$]   $q$ decides the value of $\name{\tau}_i$ (when $i<\theta$). 
\end{description}
\end{quotation}
Checking that this definition is correct is straightforward (compare the proof
of \ref{cl7}). At successor stages $i<\theta$ we use \ref{cl7} to show that
there is a condition $q'\in\bP_\gamma$ satisfying {\bf (i)$_i$}+{\bf (ii)$_i$}
and next we extend it to a condition $q$ deciding the value of
$\name{\tau}_i$. At limit stages $i\leq\theta$ we know, by the definition of
$\bar{p}\rest i$, that for each $j\leq i$, $\bar{p}\rest j\in N_{j+1}$. 
Moreover, we may apply \ref{cl6} to $\bar{N}\rest (i+1)$, $\bar{a}\rest (i+1)$
and $\bar{p}\rest i$ to conclude that $\bar{p}\rest i$ has an upper bound
$q'\in \bP_\gamma$. Now take $q\geq q'$ which decides the value of
$\name{\tau}_i$ (if $i<\theta$) -- it satisfies the demands {\bf (i)$_i$}--{
\bf (iii)$_i$}.\\
Finally look at the condition $p_\theta\in \bP_\gamma$: it forces values to
all $\name{\tau}_i$ (for $i<\theta$) and so $p_\theta\forces_{\bP_\gamma}
\langle\name{\tau}_i: i<\theta\rangle\in\V$, finishing the proof of the claim
and thus that of the theorem. \QED 

\begin{definition}
\label{ruler}
\begin{enumerate}
\item Let $\ckm$ be the family of all subsets of 
\[\begin{array}{r}
\big\{\bar{a}=\langle a_i: i\leq\alpha\rangle: \mbox{ the sequence $\bar{a}$
is increasing continuous, }\ \\ 
\alpha<\kappa\ \mbox{ and }(\forall i\leq\alpha)(a_i\in [\mu^*]^{<\kappa}\ \&\
a_i\cap\kappa\in\kappa)\big\}. 
\end{array}\]
\item Let $\bar{M}=\langle M_i: i\leq\alpha\rangle$ be an increasing
continuous sequence of elementary submodels of $\Hchi$, $\baza_0,\baza_1\in
\ckm$. We say that {\em $\bar{M}$ is ruled by $(\baza_0,\baza_1)$} if 
\begin{description}
\item[(a)] $\bar{M}\rest (i+1)\in M_{i+1}$ and $2^{\|M_i\|}+1\subseteq
M_{i+1}$ for all $i<\alpha$,
\item[(b)] $\langle M_i\cap\mu^*: i\leq\alpha\rangle\in\baza_1$,
\item[(c)] for each $i<\alpha$ (and we allow $i=-1$) there is an
$\baza_0$--complementary pair $(\bar{N}^i,\bar{a}^i)$ such that
\begin{description}
\item[($\alpha$)] $\lh(\bar{N}^i)=\lh(\bar{a}^i)=\delta_i+1$, $\cf(\delta_i)>
2^{\|M_i\|}$ and, for simplicity, $\delta_i$ is additively indecomposable,
\item[($\beta$)]  $\bar{M}\rest (i+1)\in N^i_0$, $N^i_{\delta_i}=M_{i+1}$ and
\item[($\gamma$)] $\|N^i_\varepsilon\|^{2^{\|M_i\|}}+1\subseteq
N^i_{\varepsilon+1}$. 
\end{description}
\end{description}
The sequence $\langle \bar{N}^i:i<\alpha\rangle$ given by the clause
{\bf (c)} above will be called {\em an $\baza_0$--approximation to $\bar{M}$}.
\item $\ckp$ is the family of all pairs $(\baza_0,\baza_1)$ such that
$\baza_0,\baza_1\in\ckm$, $\baza_0$ is closed and for every large enough
regular cardinal $\chi$, for every $x\in {\cal H}(\chi)$ there is a sequence
$\bar{M}$ ruled by $(\baza_0,\baza_1)$ and such that $x\in M_0$ and every end
segment of $\bar{M}$ is ruled by $(\baza_0,\baza_1)$ (i.e.\ $\baza_1$ is
closed under end segments).
\end{enumerate}
\end{definition}

\begin{remark}
\label{remB48}
{\em
\begin{enumerate}
\item Condition \ref{ruler}(2)(c) is the replacement for 
\[\|N_{i+1}\|=\lambda\quad\mbox{ and }\quad (N_{i+1})^{<\lambda}\subseteq
N_{i+1}\]
in {\bf Case A}. Here, there are no natural closed candidates for $M_{i+1}$,
as in that case. So we use a relative candidate.
\item In \ref{ruler}(2)(c)($\gamma$) we may put stronger demands (if required
in applications). For example one may consider a demand that
$\|N^i_\varepsilon\|^{h^*(\|M_i\|)}+1\subseteq N^i_{\varepsilon+1}$, for some
function $h^*:\kappa\longrightarrow\kappa$.
\item Note that if $(\baza_0,\baza_1)\in\ckp$ then necessarily
$\baza_0\in\ck$.

[Why? If $\theta=\cf(\theta)<\kappa$, $x=\langle\theta,y\rangle$ then
$\lh(\bar{N}^i)>\theta$.]
\item Note that in examples there is no need to assume that $\baza_1$ is
closed under end segments as ``complete for $(\baza_0,\baza_1)$'' (see
\ref{grapik}) is preserved, as this just restricts the choice of the ``bad
guy'' $\inc$ of $i_0$ (and so $p$) to those in the end segment.
\end{enumerate}
}
\end{remark} 

\begin{definition}
\label{grapik}
Let $(\baza_0,\baza_1)\in\ckp$ and let $\bQ$ be a forcing notion. 
\begin{enumerate}
\item For a sequence $\bar{M}=\langle M_i:i\leq\delta\rangle$ ruled by
$(\baza_0,\baza_1)$ with an $\baza_0$--approximation $\langle\bar{N}^i:
i<\delta\rangle$ and a condition $r\in\bQ$ we define a game $\grapik(\bQ,r)$
between two players $\com$ and $\inc$.
\begin{quotation}
\noindent The play lasts $\delta$ moves during which the players construct a
sequence $\langle i_0,p,\langle p_i,\bar{q}_i:i_0-1\leq i<\delta\rangle
\rangle$ such that $i_0<\delta$ is non-limit, $p\in M_{i_0}\cap\bQ$, $p_i
\in M_{i+1}\cap\bQ$,
$\bar{q}_i=\langle q_{i,\varepsilon}:\varepsilon<\delta_i\rangle\subseteq\bQ$
(where $\delta_i+1= \lh(\bar{N}^i)$).\\ 
The player $\inc$ first decides what is $i_0<\delta$ and then he chooses a
condition $p\in\bQ\cap M_{i_0}$ stronger than $r$. Next, at the stage $i\in
[i_0-1,\delta)$ of the game, $\com$ chooses $p_i\in\bQ\cap M_{i+1}$ such that  
\[p\leq_{\bQ}p_i\quad\mbox{ and }\quad (\forall j<i)(\forall \varepsilon<
\delta_j)(q_{j,\varepsilon}\leq_{\bQ} p_i),\]
and $\inc$ answers choosing an increasing sequence $\bar{q}_i=\langle q_{i,
\varepsilon}:\varepsilon<\delta_i\rangle$ such that $p_i\leq_{\bQ} q_{i,0}$
and $\bar{q}_i$ is $(\bar{N}^i\rest [\alpha,\delta_i],\bQ)^*$--generic for
some $\alpha<\delta$.    
\end{quotation}
The player $\com$ wins if he has always legal moves and the sequence $\langle
p_i: i<\delta\rangle$ has an upper bound. 
\item We say that the forcing notion $\bQ$ is {\em complete for
$(\baza_0,\baza_1)$} if
\begin{description}
\item[(a)] $\bQ$ is strongly complete for $\baza_0$ and
\item[(b)] for a large enough regular $\chi$, for some $x\in{\cal H}(\chi)$,
for every sequence $\bar{M}$ ruled by $(\baza_0,\baza_1)$ with an
$\baza_0$--approximation $\langle\bar{N}^i:i<\delta\rangle$ and such that
$x\in M_0$ and for any condition $r\in\bQ\cap M_0$, the player $\inc$ DOES NOT
have a winning strategy in the game $\grapik(\bQ,r)$.
\end{description}
\end{enumerate}
\end{definition}

\begin{proposition}
\label{B410}
Assume 
\begin{description}
\item[(a)] $(\baza_0,\baza_1)\in\ckp$,
\item[(b)] $\bQ$ is a forcing notion complete for $(\baza_0,\baza_1)$.
\end{description}
Then $\forces_{\bQ}$`` $(\baza_0,\baza_1)\in\ckp$ ''.
\end{proposition}

\Proof Straight (and not used). \QED

\section{Examples for inaccessible}
Let us look at a variant of the examples presented in section
\ref{sekcjaexamples}  relevant for our present case. (Remember 
\ref{remB48}(4).) 

\begin{hypothesis}
\label{hipoteza}
Assume that $\kappa$ is a strongly inaccessible cardinal, $S\subseteq\kappa$
is a stationary set and $\bar{C}=\langle C_\delta:\delta\in S\rangle$ is such
that for each $\delta\in S$:
\begin{quotation}
\noindent $C_\delta$ is a club of $\delta$ such that $\otp(C_\delta)<\delta$,
moreover for simplicity $\otp(C_\delta)<\min(C_\delta)$, $\nacc(C_\delta)
\subseteq\kappa\setminus S$ and

\noindent if $\alpha\in\nacc(C_\delta)$ then $\cf(\alpha)>2^{\max(\alpha\cap
C_\delta)}$ and $S\cap\alpha$ is not stationary,

\noindent if $\alpha\in\acc(C_\delta)\cap S$ then $C_\alpha=C_\delta\cap
\alpha$

\noindent [note that if $S$ does not reflect then we can ask that the last
demand does not occur].
\end{quotation}
Further we assume that $\bar{C}$ guesses clubs, i.e.
\begin{quotation}
\noindent if $E\subseteq\kappa$ is a club 

\noindent then the set $\{\delta\in S: C_\delta\subseteq E\}$ is stationary. 
\end{quotation}
Moreover we demand that for every club $E\subseteq\kappa$, the set
$\kappa\setminus S$ contains arbitrarily long (but $<\kappa$) increasing
continuous sequences from $E$.  
\end{hypothesis}

\begin{definition}
\label{forcing}
Let $\kappa,S,\bar{C}$ be as in Hypothesis \ref{hipoteza} and let
$\mu^*=\kappa$. 
\begin{enumerate}
\item Define
\[\begin{array}{ll}
\baza^S_0=\big\{\bar{\alpha}=\langle\alpha_i:i\leq\gamma\rangle:&\bar{\alpha}
\mbox{ is an increasing continuous sequence}\\
\ &\mbox{of ordinals from }\kappa\setminus S,\quad\gamma<\kappa\big\}
  \end{array}\]
\[\begin{array}{ll}
\baza^{S,\bar{C}}_1=\big\{\bar{\beta}':&\bar{\beta}' \mbox{ is an end segment
(not necessarily proper) of }\bar{\beta}\conc\langle\delta\rangle,\\
\ &\delta\in S \mbox{ and }\bar{\beta}\mbox{ is the increasing enumeration of
}C_\delta\big\}. 
  \end{array}\]
\item Suppose that $\bar{A}=\langle A_\delta: \delta\in S\rangle$, $\bar{h}=
\langle h_\delta:\delta\in S\rangle$ and $\cf(\theta)=\theta<\kappa$ are such
that for each $\delta\in S$:
\[A_\delta\subseteq\delta,\quad \|A_\delta\|<\theta,\quad h_\delta:A_\delta
\longrightarrow\theta,\quad\mbox{ and }\sup(A_\delta)=\delta\]
(we may omit the last demand as only $\bar{A}\rest S'$, for $S'=\{\delta\in S:
\delta=\sup A_\delta\}$, affects the forcing). We define a forcing notion
$\bQ^{S,\theta}_{\bar{A},\bar{h}}$: 

\noindent{\bf a condition} in $\bQ^{S,\theta}_{\bar{A},\bar{h}}$ is a function
$g:\beta\longrightarrow\theta$ (for some $\beta<\kappa$) such that
\[(\forall\delta\in S\cap (\beta+1))(\{\xi\in A_\delta: h_\delta(\xi)\neq
g(\xi)\}\mbox{ is bounded in }\delta),\]

\noindent{\bf the order} $\leq_{\bQ^{S,\theta}_{\bar{A},\bar{h}}}$ of
$\bQ^{S,\theta}_{\bar{A},\bar{h}}$ is the inclusion (extension). 
\item For $\bar{A}$, $\bar{h}$ and $\theta$ as above and $\alpha<\kappa$ we
let 
\[{\cal I}^{\bar{A},\bar{h},\theta}_{\alpha}\stackrel{\rm def}{=}\{g\in
\bQ^{S,\theta}_{\bar{A},\bar{h}}: \alpha\in\dom(g)\}.\]
\end{enumerate}
\end{definition}

\begin{remark}
{\em 
One of the difficulties in handling the forcing notion $\bQ^{S,\theta}_{
\bar{A},\bar{h}}$ is that the sets ${\cal I}^{\bar{A},\bar{h},\theta}_{
\alpha}$ do not have to be dense in $\bQ^{S,\theta}_{\bar{A},\bar{h}}$. Of
course, if this happens then the generic object is not what we expect it to
be. However, if the set $S$ is not reflecting then each ${\cal I}^{\bar{A},
\bar{h},\theta}_{\alpha}$ is dense in $\bQ^{S,\theta}_{\bar{A},\bar{h}}$ and
even weaker conditions are enough for this. One of them is the following:
\begin{description}
\item[$(*)$] $\langle A_\delta:\delta\in S\rangle$ is $\kappa$--free, i.e.\
for every $\alpha<\kappa$ there is a function $g$ such that $\dom(g)=S\cap
\alpha$ and $g(\delta)<\delta$ and the sets $\langle A_\delta\setminus
g(\delta): \delta\in S\cap\alpha\rangle$ are pairwise disjoint.
\end{description}
We can of course weaken it further demanding that $\langle A_\delta:\delta\in
S\cap\alpha\rangle$ has uniformization. (So if we force inductively on all
$\kappa$'s this may be reasonable, or we may ask uniformization just for our
$h_\delta$'s.)
}
\end{remark}

\begin{proposition}
$(\baza^S_0,\baza^{S,\bar{C}}_1)\in\ckp$
\end{proposition}

\Proof Immediately from its definition we get that $\baza^S_0$ is closed. 
Suppose now that $\chi$ is a sufficiently large regular cardinal and $x\in
{\cal H}(\chi)$. First construct an increasing continuous sequence $\bar{W}=
\langle W_j: j<\kappa\rangle$ of elementary submodels of $\Hchi$ such that
$x\in W_0$ and for each $j<\kappa$ 
\[\|W_j\|<\kappa,\quad\mbox{and}\quad W_j\cap\kappa=\|W_j\|,\quad\mbox{and}
\quad \bar{W}\rest (j+1)\in W_{j+1}.\]
Note that then, for each $j<\kappa$, we have $2^{\|W_j\|}+1\subseteq W_{j+1}$.
Clearly the set $E=\{W_j\cap\kappa: j<\kappa \mbox{ is limit\/}\}$ is a club
of $\kappa$ and so $\acc(E)$ is a club of $\kappa$ as well. Thus, by our
assumptions on $\bar{C}$ (see \ref{hipoteza}), we find $\delta\in S$ such that
$C_\delta\subseteq\acc(E)$ (then, of course, $\delta\in \acc(E)$ too). Let
$\bar{M}=\langle M_i:i\leq\otp(C_\delta)\rangle$ be the increasing enumeration
of 
\[\big\{W_j: j<\kappa\ \&\ W_j\cap\kappa\in C_{\delta}\cup\{\delta\}\big\}.\]
Fix $i<\otp(C_\delta)$. Let $j'<j<\kappa$ be such that $W_{j'}=M_i$ and
$W_j=M_{i+1}$, and let $\alpha= M_{i+1}\cap\kappa=W_j\cap\kappa$. Then $\alpha
\in\nacc(C_\delta)\cap\acc(E)$ and, by \ref{hipoteza}, $\alpha\notin S$ and
the set $S$ does not reflect at $\alpha$. Consequently we find a club $C^i$ of
$\alpha$ disjoint from $S\cap\alpha$. Let $\bar{N}^i=\langle N^i_\varepsilon:
\varepsilon\leq\delta_i\rangle$ be the increasing enumeration of
\[\big\{W_\xi: j'<\xi\leq j\ \&\ W_\xi\cap\kappa\in C^i\cup\{\alpha\}\big\}.\]
(Note that the set above is non-empty as $\alpha\in\acc(E)$; passing to a
cofinal subsequence we may demand that $\delta_i$ is additively
indecomposable.) We claim that the sequence $\bar{M}$ is ruled by
$(\baza^S_0,\baza^{S,\bar{C}}_1)$ and $\langle \bar{N}^i: i<\otp(C_\delta)
\rangle$ is an $\baza^S_0$--approximation to $\bar{M}$. For this we have to
check the demands of \ref{ruler}(2). By the choice of the $W_j$'s we have that
the clause {\bf (a)} there is satisfied. As $\langle M_i\cap\kappa:i\leq
\otp(C_\delta)\rangle$ enumerates $C_\delta\cup\{\delta\}$ we get the demand
{\bf (b)} there. For the clause {\bf (c)} fix $i<\otp(C_\delta)$ and look at
the way we defined $\bar{N}^i=\langle N^i_\varepsilon:\varepsilon\leq\delta_i
\rangle$. For each $\varepsilon\leq \delta_i$, $N^i_\varepsilon\cap\kappa\in
C^i\cup\{\alpha\}\subseteq\kappa\setminus S$. Hence $(\bar{N}^i,\langle
N^i_\varepsilon\cap\kappa: \varepsilon\leq\delta_i\rangle)$ is an
$\baza^S_0$--complementary pair. Moreover, 
\[\cf(\delta_i)=\cf(\alpha)>2^{\max(\alpha\cap C_\delta)}=2^{M_i\cap\kappa} =
2^{\|M_i\|}\]
(by \ref{hipoteza}) and $\delta_i$ is additively indecomposable. This verifies
{\bf (c)}($\alpha$). The clauses {\bf (c)}($\beta$) and {\bf (c)}($\gamma$)
should be clear by the choice of the $W_j$'s and that of $\bar{N}^i$. \QED

\begin{proposition}
Suppose that $\bar{A}$, $\bar{h}$, $\theta$ are as in \ref{forcing}(2) and for
each $\alpha<\kappa$ the set ${\cal I}^{\bar{A},\bar{h},\theta}_\alpha$ (see
\ref{forcing}(3)) is dense in $\bQ^{S,\theta}_{\bar{A},\bar{h}}$ (e.g.~$S$
does not reflect).\\
Then the forcing notion $\bQ^{S,\theta}_{\bar{A},\bar{h}}$ is complete for
$(\baza^S_0,\baza^{S,\bar{C}}_1)$.
\end{proposition}

\Proof We break the proof to three steps checking the requirements of
\ref{grapik}(2).

\begin{claim}
\label{cl10}
$\bQ^{S,\theta}_{\bar{A},\bar{h}}$ is complete for $\baza^S_0$.
\end{claim}

\noindent{\em Proof of the claim:}\hspace{0.15in} Suppose that
$(\bar{N},\bar{a})$ is an $\baza^S_0$--complementary pair, $\bar{N}=\langle
N_i: i\leq\delta\rangle$ and $\bar{p}=\langle p_i: i\leq\delta\rangle\subseteq
\bQ^{S,\theta}_{\bar{A},\bar{h}}$ is and increasing $(\bar{N},\bQ^{S,\theta}_{
\bar{A},\bar{h}})^n$--generic sequence. Let $p\stackrel{\rm
def}{=}\bigcup\limits_{i<\delta} p_i$. Note that $p$ is a function from
$\dom(p)=\bigcup\limits_{i<\delta}\dom(p_i)$ to $\theta$. Moreover, as the
sets ${\cal I}^{\bar{A},\bar{h},\theta}_\alpha$ are dense in $\bQ^{S,\theta}_{
\bar{A},\bar{h}}$ (and ${\cal I}^{\bar{A},\bar{h},\theta}_\alpha\in N_i$ if
$\alpha\in N_i\cap\kappa$), we have $N_i\cap\kappa\subseteq\dom(p_{i+n})
\subseteq N_{i+n+1}$. Hence 
\[\dom(p)=\bigcup_{i<\delta}N_i\cap\kappa=N_\delta\cap\kappa\in\kappa.\]
Note that $N_\delta\cap\kappa\notin S$ (by the definition of $\baza^S_0$). 
Suppose that $\alpha\in S\cap (\dom(p)+1)$, so $\alpha\in\dom(p)$. Then for
some $i<\delta$ we have $\alpha\in\dom(p_i)$ and, as $p_i\in\bQ^{S,\theta}_{
\bar{A},\bar{h}}$, the set $\{\xi\in A_\alpha: h_\alpha(\xi)\neq p(\xi)=
p_i(\xi)\}$ is bounded in $\alpha$. This shows that $p\in \bQ^{S,\theta}_{
\bar{A},\bar{h}}$ and clearly it is an upper bound of $\bar{p}$.

\begin{claim}
\label{cl11}
Forcing with $\bQ^{S,\theta}_{\bar{A},\bar{h}}$ does not add new sequences of
length $<\kappa$ of ordinals.
\end{claim}

\noindent{\em Proof of the claim:}\hspace{0.15in} Suppose that $\zeta<\kappa$
and $\name{\tau}$ is a $\bQ^{S,\theta}_{\bar{A},\bar{h}}$--name for a function
from $\zeta$ to ordinals, $p\in \bQ^{S,\theta}_{\bar{A},\bar{h}}$. Take an
increasing continuous sequence $\bar{W}=\langle W_j: j<\kappa\rangle$ of
elementary submodels of $\Hchi$ such that $\bQ^{S,\theta}_{\bar{A},\bar{h}},p,
\name{\tau}\in W_0$, $\zeta+1\subseteq W_0$ and for each $j<\kappa$
\[\|W_j\|<\kappa,\quad\mbox{and}\quad W_j\cap\kappa=\|W_j\|,\quad\mbox{and}
\quad \bar{W}\rest (j+1)\in W_{j+1}.\]
Look at the club $E=\{W_j\cap\kappa: j<\kappa\}$. By the last assumption of
\ref{hipoteza} we find an increasing continuous sequence $\langle j_\xi:\xi
\leq\zeta\rangle$ such that $\{W_{j_\xi}\cap\kappa:\xi\leq\zeta\}\cap S=
\emptyset$. Now we build inductively an increasing sequence $\langle p_\xi:
\xi\leq\zeta\rangle$ of conditions from $\bQ^{S,\theta}_{\bar{A},\bar{h}}$
such that $p\leq_{\bQ^{S,\theta}_{\bar{A},\bar{h}}} p_0$ and for each
$\xi<\zeta$: 
\begin{enumerate}
\item $p_\xi\in W_{j_{\xi+1}}$,
\item $p_\xi$ forces a value to $\name{\tau}(\xi)$, and
\item $W_{j_\xi}\cap\kappa\subseteq\dom(p_\xi)$.
\end{enumerate}
There are no problems with carrying out the construction. At a non-limit stage
$\xi$, we may easily choose a condition $p_\xi$ in $W_{j_{\xi+1}}$ stronger
than the condition chosen before (if any) and such that $W_{j_\xi}\cap\kappa
\subseteq\dom(p_\xi)$ (remember that ${\cal I}^{\bar{A},\bar{h},\theta}_{
W_{j_\xi}\cap\kappa}\in W_{j_{\xi+1}}$ is a dense subset of $\bQ^{S,\theta}_{
\bar{A},\bar{h}}$) and $p_\xi$ decides the value of
$\name{\tau}(\xi)$. Arriving at a limit stage $\xi\leq\zeta$ we take the union
of conditions chosen so far and we note that it is a condition in
$\bQ^{S,\theta}_{\bar{A},\bar{h}}$  as
\[\dom(\bigcup_{i<\xi}p_i)=\bigcup_{i<\xi}\dom(p_i)=\bigcup_{i<\xi} W_{j_i}
\cap\kappa= W_{j_\xi}\cap\kappa\notin S.\]
Now proceed as in the successor case. Finally look at the condition $p_\zeta$
-- it decides the value of $\name{\tau}$ (and is stronger than $p$).

\begin{claim}
\label{cl12}
Assume that $\bar{M}=\langle M_i:i\leq\delta\rangle$ is an increasing
continuous sequence of elementary submodels of $\Hchi$  ruled by
$(\baza_0,\baza_1)$ with an $\baza_0$--approximation $\langle\bar{N}^i:
i<\delta\rangle$ and such that $S,\baza^S_0,\baza^{S,\bar{C}}_1,\bar{A},
\bar{h},\theta,\bQ^{S,\theta}_{\bar{A},\bar{h}}\in M_0$. Let $r\in\bQ^{S,
\theta}_{\bar{A},\bar{h}}\cap M_0$. Then the player $\com$ has a winning
strategy in the game $\grapik(\bQ^{S,\theta}_{\bar{A},\bar{h}},r)$. 
\end{claim}

\noindent{\em Proof of the claim:}\hspace{0.15in} Let us describe a strategy
for player $\com$ in the game $\grapik(\bQ^{S,\theta}_{\bar{A},\bar{h}},r)$.\\
Since $\langle M_i: i\leq\delta\rangle\in\baza^{S,\bar{C}}_1$ and for each
$\alpha\in S$, $\otp(C_\alpha)<\alpha$ (see \ref{hipoteza}) we know that
$\delta=\otp(C_{M_\delta\cap\kappa})<M_\delta\cap\kappa\in S$. Recall
$\otp(C_\delta)<\min(C_\delta)$. Let 
\[Z\stackrel{\rm def}{=}\bigcup\{A_{M_i\cap\kappa}:\ i\leq\delta\ \ \&\ \
M_i\cap\kappa\in S\}.\]
Note that $\|Z\|\leq\delta\cdot\theta<\|M_{i_0}\|$. By induction on $i\leq
\theta^+$ choose an increasing continuous sequence $\langle Z_i:i\leq\theta^+
\rangle$ of subsets of $\kappa$ such that $Z_0=Z$ and $Z_{i+1}=Z_i\cup \bigcup
\{A_\alpha: \alpha\in S\ \&\ \alpha=\sup(Z_i\cap \alpha)\}$. Clearly $\|Z_i\|
\leq\delta\cdot\theta\cdot\|i\|$ for each $i\leq\theta^+$ and if $\alpha=
\sup(\alpha\cap Z_{\theta^+})$ then $A_\alpha\subseteq Z_{\theta^+}$. So as
$\|A_\alpha\|\leq\theta$ we have
\[\alpha\in S\ \&\ \alpha=\sup(Z_{\theta^+}\cap\alpha)\quad\Rightarrow\quad
A_\alpha\subseteq Z_{\theta^+}.\]
Now, in his first move, player $\inc$ chooses non-limit $i_0<\delta$ and $p
\in\bQ^{S,\theta}_{\bar{A},\bar{h}}\cap M_{i_0}$ stronger than $r$. We have
assumed that each ${\cal I}^{\bar{A},\bar{h},\theta}_\xi$ (for $\xi<\kappa$)
is dense in $\bQ^{S,\theta}_{\bar{A},\bar{h}}$, so we have a condition
$p^+\in\bQ^{S,\theta}_{\bar{A},\bar{h}}$ stronger than $p$ and such that
$M_\delta \cap\kappa\in \dom(p^+)$. In the next steps, the strategy for $\com$
will have the property that for each $i\geq i_0-1$ it says $\com$ to play a
condition $p_i\in \bQ^{S,\theta}_{\bar{A},\bar{h}}$ such that  
\begin{description}
\item[$(\boxdot)_i$]\qquad $Z_{\theta^+}\cap M_{i+1}\subseteq\dom(p_i)$ and
$p_i\rest Z_{\theta^+}=p^+\rest (Z_{\theta^+}\cap M_{i+1})$.
\end{description}
So, first the player $\com$ chooses a condition $p_{i_0-1}\in M_{i_0}\cap
\bQ^{S,\theta}_{\bar{A},\bar{h}}$ stronger than $p$ and such that
\[Z_{\theta^+}\cap M_{i_0}\subseteq\dom(p_{i_0})\quad\mbox{ and }\quad
p_{i_0}\rest (Z_{\theta^+}\cap M_{i_0})=p^+\rest(Z_{\theta^+}\cap
M_{i_0}).\]
Why is it possible? We know that  
\[\|Z_{\theta^+}\cap M_{i_0}\|\leq\theta^+<M_0\cap\kappa\leq\|M_{i_0}\|<\cf(
\delta_{i_0-1})\]
(where $\delta_{i_0-1}+1=\lh(\bar{N}^{i_0-1})$) and therefore $Z_{\theta^+}
\cap M_{i_0}\subseteq N^{i_0-1}_\varepsilon$ for some $\varepsilon<\delta_{
i_0}$. Taking possibly larger $\varepsilon$ we may have $\dom(p)\subseteq
N^{i_0-1}_{\varepsilon}$ too. Let $p'\in N^{i_0-1}_{\varepsilon+1}\cap\bQ^{S,
\theta}_{\bar{A},\bar{h}}$ be such that $p\leq p'$ and $N^{i_0-1}_\varepsilon
\cap\kappa\subseteq\dom(p')$. Let 
\[p_{i_0}=p'\rest(\dom(p')\setminus Z_{\theta^+}\cup p^+\rest Z_{\theta^+}\cap
N^{i_0-1}_\varepsilon.\]
Note that $p_{i_0}:\dom(p_{i_0})\longrightarrow\theta$ is a well defined
function such that $p_{i_0}\in N^{i_0-1}_{\varepsilon+1}$ (for the last
remember \ref{ruler}(2)(c)($\gamma$): we are sure that $Z_{\theta^+}\cap
N^{i_0-1}_\varepsilon\in N^{i_0-1}_{\varepsilon+1}$ and $p^+\rest
(Z_{\theta^+}\cap N^{i_0-1}_\varepsilon)\in N^{i_0-1}_{\varepsilon+1}$, as
$\|N^{i_0-1}_\varepsilon\|^{\delta\cdot\theta^+}+1\subseteq N^{i_0-1}_{
\varepsilon+1}$). Finally, to check that $p_{i_0}$ is a condition in
$\bQ^{S,\theta}_{\bar{A},\bar{h}}$ suppose that $\gamma\in S\cap(\dom(p_{i_0-
1})+1)$. If $A_\gamma\subseteq Z_{\theta^+}$ then $p_{i_0}\rest A_\gamma=p^+
\rest A_\gamma$ and the requirement of \ref{forcing}(2) is satisfied. If
$A_\gamma$ is not contained in $Z_{\theta^+}$ then necessarily $Z_{\theta^+}
\cap\gamma$ is bounded in $\gamma$ and we use the fact that $p_{i_0-1}\rest
(A_\gamma\setminus Z_{\theta^+})=p'\rest (A_\gamma\setminus Z_{\theta^+})$,
$p'\in \bQ^{S,\theta}_{\bar{A},\bar{h}}$.\\  
At a stage $i\in [i_0,\delta)$ of the game the player $\com$ applies a similar
procedure, but first he looks at the union $p^*_i=\bigcup\limits_{j<i}
\bigcup\limits_{\varepsilon<\delta_j}q_{j,\varepsilon}$ of all conditions
played by his opponent so far. If $i$ is not limit then, directly from
\ref{cl10}, we know that $p^*_i$ is a condition in $\bQ^{S,\theta}_{\bar{A},
\bar{h}}$ (stronger than $r$). But what if $i$ is limit? In this case the
demands $(\boxdot)_j$ for $j<i$ help. The only possible trouble could come
from $A_{M_i\cap\kappa}$ when $M_i\cap\kappa\in S$. But then the set
$Z_{\theta^+}$ contains $A_{M_i\cap\kappa}$ and, by $(\boxdot)_j$ for $j<i$,
$p^*_i\rest A_{M_i\cap\kappa}=p^+\rest A_{M_i\cap\kappa}$. This implies that
the set 
\[\{\xi\in A_{M_i\cap\kappa}: h_{M_i\cap\kappa}(\xi)\neq p^*_i(\xi)\}\]
is bounded in $M_i\cap\kappa$. Hence easily $p^*_i\in\bQ^{S,\theta}_{\bar{A},
\bar{h}}$. Next, player $\com$ extends the condition $p^*_i$ to $p_i\in
N^i_{\varepsilon+1}\cap \bQ^{S,\theta}_{\bar{A},\bar{h}}$ (for some
$\varepsilon<\delta_i$) such that the demand $(\boxdot)_i$ is satisfied,
applying a procedure similar to the one described for getting $p_{i_0}$.

Why is the strategy described above the winning one? Suppose that $\langle
p_i: i_0-1\leq i<\delta\rangle$ is a sequence constructed by $\com$ during a
play in which he uses this strategy. As it is an increasing sequence of
conditions and $\bigcup\limits_{i<\delta}\dom(p_i)=M_\delta\cap\kappa$, the
only thing we should check is that the set 
\[\{\xi\in A_{M_\delta\cap\kappa}: h_{M_\delta\cap\kappa}(\xi)\neq
(\bigcup_{i<\delta}p_i)(\xi)\}\]
is bounded in $M_\delta\cap\kappa$. But by the choice of $Z\subseteq
Z_{\theta^+}$, and by keeping the demand $(\boxdot)_i$ (for $i<\delta$) we
know that   
\[\{\xi\in A_{M_\delta\cap\kappa}: h_{M_\delta\cap\kappa}(\xi)\neq
(\bigcup_{i<\delta}p_i)(\xi)\}\subseteq\{\xi\in A_{M_\delta\cap\kappa}:
h_{M_\delta\cap\kappa}(\xi)\neq p^+(\xi)\},\]
so the choice of $p^+$ works.\\
This finishes the proof of the claim and that of the proposition. \QED
\bigskip

Now, let us turn to the applications for Abelian groups (i.e.~the
forcing notions needed for \ref{Gques}). We continue to use Hypothesis
\ref{hipoteza}. 

\begin{definition}
Assume that $G$ is a strongly $\kappa$--free Abelian group and $h:H
\stackrel{\rm onto}{\longrightarrow}G$ is a homomorphism onto $G$ with kernel
$K$ of cardinality $<\kappa$. We define a forcing notion $\bP_{h,H,G}$:

\noindent {\bf a condition} in $\bP_{h,H,G}$ is a function $q$ such that
\begin{description}
\item[(a)] $\dom(q)$ is a subgroup of $G$ of size $<\kappa$,
\item[(b)] $G/\dom(q)$ is $\kappa$--free,
\item[(c)] $q$ is a lifting for $\dom(q)$ and $h:H\longrightarrow G$;
\end{description}

\noindent {\bf the order} $\leq_{\bP_{h,H,G}}$ of $\bP_{h,H,G}$ is the
inclusion (extension).
\end{definition}

\begin{hypothesis}
\label{hypB561}
Let $\bar{G}=\langle G_i: i<\kappa\rangle$ be a filtration of $G$,
$\Gamma[G]\subseteq S$ (modulo the club filter on $\kappa$). So wlog
$\gamma[\bar{G}]\subseteq S$ and $\gamma[\bar{G}]$ is a set of limit
ordinals. Let $h^{-1}[G_i]=H_i$.
\end{hypothesis}

\begin{proposition}
For each $\alpha<\kappa$ the set
\[\bP_{h,H,G}^*\stackrel{\rm def}{=}\{q\in \bP_{h,H,G}: (\exists i<\kappa)(
\dom(q)=G_{i+1}\ \&\ i\geq\alpha\}\]
is dense in $\bP_{h,H,G}$.
\end{proposition}

\Proof Let $q\in \bP_{h,H,G}$ and let $i<\kappa$ be such that $\dom(q)
\subseteq G_i$ and $i\geq \alpha$. Then $G/dom(q)$ is $\kappa$--free and so
$G_{i+1}/\dom(q)$ is free. Now consider the mapping $x\mapsto x+\dom(q):
G_{i+1}\longrightarrow G_{i+1}/\dom(q)$. So by \ref{getL} we get that
$G_{i+1}=\dom(q)+L$ for some free $L$ ($L\cong G_{i+1}/\dom(q)$). 
Consequently, there is a lifting $f$ of $L$ and now $\langle f,q\rangle$
extends $q$ and it is in $\bP_{h,H,G}'$ \QED 

\begin{proposition}
\label{propB571}
The forcing notion $\bP_{h,H,G}$ is strongly complete for $\hat{\baza}_0$.
\end{proposition}

\Proof The two parts, not adding bounded subsets of $\kappa$ and completeness
for $\hat{\baza}_0$, are similar to those for uniformization, so we do just
the second.

Assume now that $\chi$ is a regular large enough cardinal, $N_i\prec\Hchi$,
$\bar{N}=\langle N_i:i\leq\delta\rangle$, $\bar{N}\rest (i+1)\in N_{i+1}$,
$N_i\cap\kappa\in \kappa$ is limit, $\bar{N}$ obeys $\bar{a}\in\baza_0^S$ and
$\bar{p}=\langle p_i:i<\delta \rangle$ is generic for $\bar{N}$ with error $n$,
$\dom(p_i)=G_{\gamma_i+1}$. In particular, $p_i\in N_{i+1}$ and as $\gamma_i$
is computable from $\bar{G},p_i$ we know that $\gamma_i\in N_{i+1}$. 

Let $\beta_i=\sup(N_i\cap\kappa)$ (so the sequence $\langle\beta_i:i\leq\delta
\rangle$ is increasing continuous). Note that 
\[p_{i+n}\in\bigcap\{{\cal I}\in N_{i+1}: {\cal I}\subseteq\bP_{h,H,G}'\mbox{
is open dense\/}\}\]
and $N_i,\beta_i\in N_{i+1}$. Moreover, the set 
\[{\cal I}_{\beta_i}=\{q\in \bP_{h,H,G}': \dom(q)\supseteq G_{\beta_i}\}\]
is open dense in $\bP_{h,H,G}'$. So $p_{i+n}\in {\cal I}_{\beta_i+1}\in N_{i+
1}$ and $\gamma_{i+n}>\beta_i$. Now, $\dom(\bigcup\limits_{i<\delta} p_i)=
G_{\bigcup\limits_{i<\delta}(\gamma_i+1)}$ and $\bigcup\limits_{i<\delta}(
\gamma_i+1)=\bigcup\limits_{i<\delta}\beta_i=N_\delta\cap\kappa$. Since
$N_\delta\cap\kappa\notin S$ and $S\supseteq\Gamma[G]$ we conclude
$N_\delta\cap\kappa\notin \Gamma[G]$, and thus $G_{N_\delta\cap\kappa+1}/
G_{N_\delta\cap\kappa}$ is free. So we can complete to a condition. \QED

\begin{proposition}
\label{forforgoe}
The forcing notion $\bP_{h,H,G}$ is complete for $(\baza_0^S,\baza_1^S)$.
\end{proposition}

\Proof Suppose that $\bar{M}=\langle M_i:i\leq\delta\rangle$ is ruled by
$(\baza^S_0,\baza^S_1)$. So $M_i\cap \kappa=a_i$ and $M_{i+1}=\bigcup\limits_{
\zeta<\cf(a_{i+1})} N^i_\zeta$ and $(\bar{N}^i,\bar{b}^i)$ is an
$\baza^S_0$--complementary pair, $\bar{b}^i\in\baza^S_0$ (also for $i=-1$). 

We are dealing with the case $\delta\in M_0$. Recall:

\begin{claim}
There is $G^+_i$ such that $G_i\subseteq G^+_i\subseteq G_{i+1}$, $\|G^+_i\|
\leq \|G_i\|+\aleph_0$ and $G_{i+1}/G^+_i$ is free. (Of course if $i\notin S$
is non-limit then $G^+_i/G_i$ is free.)
\end{claim}

\noindent{\em Proof of the claim:}\hspace{0.15in} Since $G_{i+1}$ is free we
may fix a basis $\langle x_{i,\varepsilon}:\varepsilon<\varepsilon_{i+1}
\rangle$ of it. Choose $A_i\subseteq \varepsilon_{i+1}$ such that $\|A_i\|
\leq\|G_i\|$ and $G_i\subseteq\langle\{x_{i,\varepsilon}:\varepsilon\in A_i
\}\rangle_G$ (and call the last group $G^+_i$). Then $G_{i+1}/G^+_i$ is freely
generated by $\{x_{i,\varepsilon}+G^+_i:\varepsilon\in\varepsilon_{i+1}
\setminus A_i\}$. The claim is proved.
\medskip

Let $H^+_\alpha=h^{-1}[G^+_\alpha]$. 

Thus if $i<j$ then $G_j/G^+_i$ is free. All action will be in $G^+_{a_i}/
G_{a_i}$ for limit $i\leq\delta$. Necessarily $a_i$ is a singular cardinal of
small cofinality ($\leq\delta<a_0$). [Remember $a_i=M_i\cap\kappa$ and
$\sup(M_i\cap\kappa)$ is a limit cardinal. Why? If not then there is a
cardinal $\lambda$ such that $\lambda<\sup(M_i\cap\kappa)<\lambda^+$, so there
is $\gamma\in M_i\cap\kappa$ such that $\lambda<\gamma<\sup(M_i\cap\kappa)<
\lambda^+$. Hence $\lambda^+=\|\gamma\|^+\in M_i$, a contradiction.]

We may have ``a difficulty'' in defining $p\rest H_{a_i}$, so we should
``think'' about it earlier. This will mean defining $p\rest G_{a_{j+1}}$,
$j<i$. The player $\com$ can give only a condition in $M_{j+1}$, and we will
arrange that our ``prepayments'' are of ``size'' $a_j$ (so bounded in
$M_{j+1}$ and thus included in some $N^j_\zeta$, $\zeta<\cf(a_{j+1})$; the
will even belong to it).

Let $r\in\bP_{h,H,G}'\cap M_0$. [Remember: $G_{a_0}/\dom(r)$ is free, so
there is a lifting.] Let $\inc$ choose non-limit $i_0<\delta$ and $p'_0\in
M_{i_0}\cap \bP_{h,H,G}$ above $p$, and $\bar{q}_0=\langle q_{i_0,\zeta}
:\zeta<\delta_{i_0-1}\rangle$ {\em generic} for some end segment of
$\bar{N}_{i_0-1}$.

We choose by induction on $i\leq\delta$ models $\cB_i\prec\Hchi$ such that
\begin{itemize}
\item $\bar{G},\bar{M},\langle\bar{N}^i:i\leq\delta\rangle,\langle H_i:i\leq
\delta\rangle,\ldots\in \cB_i$, 
\item the sequence $\langle\cB_i:i\leq\delta\rangle$ is increasing (but not
continuous),
\item $\|\cB_i\|=a_i$, $a_i+1\subseteq\cB_i$ and $\langle\cB_j:j<i\rangle\in
\cB_i$,
\item $\cB_i\cap M_j\in M_{j+1}$ if $i<j$.
\end{itemize}
(But see for additional requirements later.)\\
The rest of the moves are indexed by $i\in [i_0-1,\delta)$ and in the $i^{\rm
th}$ move $\com$ chooses $p_i\in M_i$ and $\inc$ plays $\bar{q}^i=\langle
q^i_\zeta:\zeta<\delta_i\rangle$ as in the definition of the game.

Now $\com$ will choose on a side also $f_i\in \bP_{h,H,G}$ such that
additionally:
\begin{description}
\item[$(*)_1$] $f_i\in \bP_{h,H,G}$ is a function with domain $\cB_i\cap
G^+_{a_\delta}$, increasing with $i$,
\item[$(*)_2$] $a_i\subseteq\dom(f_i)$,
\item[$(*)_3$] $f_i\rest a_{i+1}=f_i\rest (a_{i+1}\cap\cB_i)$ belongs to
$\bP_{h,H,G}$ and is below $p_i$.
\end{description}
Note that
\begin{description}
\item[$(\oplus)$] $\cB_j\cap M_{j+1}\cap G^+_{a_\delta}$ is a subset of
$G_{M_{j+1}\cap\kappa}$ of cardinality $\|M_j\|<\cf(\delta_j)$, hence it
belongs to $M_{j+1}$.
\end{description}
For $i=i_0+1$ let $f_0\in \bP_{h,H,G}$ be above $p$ and have domain
$G^+_{a_\delta}$, and let $p_{i_0-1}=f_0\rest \cB_{i_0-1}\cap
M_{i_0}$. Clearly $p_{i_0-1}$ is a function extending $p$, its domain belongs
to $M_{i_0}$ and it is a subgroup of $G^+_{a_\delta}$. Consequently,
$p_{i_0-1}$ has a lifting. Now, $G^+_{a_\delta}/\dom(p_{i_0-1})$ is free as
$G_{a_{i_0}}$ is free, $G_{a_{i_0}}\in \cB_{i_0-1}$. Also $G_{a_{i_0}}/\cB_{
i_0-1}\cap G_{a_{i_0}}$ is free (manipulate bases or see \cite{Sh:52}). 

\noindent For $i=j+1\geq i_0$ we have $f_j$, $p_j$ and $\bar{q}'_j=\langle
q^j_\zeta:\delta_j\rangle$. Let $q'_j=\bigcup\limits_{\zeta<\delta_j}q_\zeta$.
So as $\dom(q'_j)=a_{j+1}=a_i\notin S$ (by the choice of $\hat{\baza}_1$),
clearly $q'_j\in \bP_{h,H,G}$. We have to find $p_i\in \bP_{h,H,G}\cap
M_{i+1}$ above $q'_j$ and $f_j\rest M_{i+1}$ (and then choose $f_i$). Clearly
the domains of $q'_j$, $f_j\rest M_{i+1}$ are pure subgroups, and $p_i$,
$f_j\rest M_j$ agree on their intersection (which is $\cB_j\cap M_{j+1}$). 
Hence there is a common extension $p'_i$, a homomorphism from $G_{a_i}+(\cB_j
\cap M_{i+1})$ to $H$, which clearly is a lifting. Does $p'_i\in \bP_{h,H,G}$?
For this it suffices to show that the group $G^+_{a_\delta}/\dom(p'_i)$ is
free. But $G^+_{a_\delta}/G^+_{a_{j+1}}$ is free, hence $(G^+_{a_\delta}/
G_{a_{j+1}})/\cB_j\cap (G^+_{a_\delta}/G_{a_{j+1}})$ is free (see
\cite{Sh:52}). Therefore $G^+_{a_\delta}/(\cB_j\cap G^+_{a_\delta}+
G_{a_{j+1}})$ is free. Also $(\cB_j\cap G^+_{a_\delta}+G_{a_{j+1}})/(\cB_j
\cap G_{a_{i+1}}+G_{a_{j+1}})$ is free (see \cite{Sh:52}). Together,
$G^+_{a_\delta}/(\cB_j\cap G_{a_{i+1}}+ G_{a_{j+1}})$ is free as required.

We are left with the case of limit $i$. Let $q'_i=\bigcup\{q'_j: i_0-1\leq
j<i\}$. Then $q'_i$ is a lifting for $G_{a_i}$. Now clearly $f'_i=\bigcup\{
f_j: i_0-1\leq j<i\}$ is a lifting for $G^+_{a_\delta}\cap\bigcup\limits_{
j<i}\cB_j$, also $G^+_{a_\delta}/\dom(f'_i)$ is free (see \cite{Sh:52}) and
$G^+_{a_i}\in\cB_0$, $\|G^+_{a_i}\|=\|G_{a_i}\|\subseteq \bigcup\limits_{j<i}
\cB_j$. Hence $G^+_{a_i}\subseteq \bigcup\limits_{j<i}\cB_j$ and therefore
$G^+_{a_i}\subseteq \dom(f'_i)$ and we can proceed to define $p_i$ as above.

Having finished the play, again $\bigcup\{f_j: i_0-1\leq j<\delta\}\in
\bP_{h,H,G}$ (as in the limit case) is an upper bound as required. \QED

\begin{remark}
{\em 
In this section, we can replace $\hat{\baza}_1$ by any $\hat{\baza}^S_1$
defined below (or any subset which is rich enough):
\[\begin{array}{ll}
\baza^S_1=\big\{\bar{\alpha}=\langle\alpha_i:i\leq\delta\rangle:&\bar{\alpha}
\mbox{ is an increasing continuous sequence}\\
\ &\mbox{of ordinals from }\kappa,\quad a_{i+1}\notin S,\ \cf(a_{i+1})>a_i\
\mbox{ and}\\
\ &S\cap a_{i+1}\mbox{ not stationary}\big\}.
  \end{array}\]
}
\end{remark}

\section{The iteration theorem for inaccessible $\kappa$}
In this section we prove the preservation theorem needed for our present
case. Like in {\bf Case A}, we will use trees of conditions. So, our way to
prove the iteration theorem will be parallel to that of {\bf Case A}.

\begin{proposition}
\label{beforeonestep}
Assume that $\baza\in\ck$ is closed and $\bar{\bQ}=\langle\bP_\alpha,
\nbQ_\alpha:\alpha<\gamma\rangle$ is a $(<\kappa)$--support iteration of
forcing notions which are strongly complete for $\baza$. Let $\tree=(T,<,\rk)$
be a standard $(w,\alpha_0)^\gamma$--tree (see \ref{standleg}),
$\|T\|<\kappa$, $w\subseteq\gamma$, $\alpha_0$ an ordinal, and let $\bar{p}=
\langle p_t: t\in T\rangle\in\ftr'(\bar{\bQ})$. Suppose that ${\cal I}$ is an
open dense subset of $\bP_\gamma$.\\
Then there is $\bar{q}=\langle q_t:t\in T\rangle\in\ftr'(\bar{\bQ})$ such that
$\bar{p}\leq \bar{q}$ and for each $t\in T$
\begin{enumerate}
\item $q_t\in\{q\rest\rk(t): q\in{\cal I}\}$,\quad and
\item for each $\alpha\in\dom(q_t)$, either $q_t(\alpha)=p_t(\alpha)$ or
$\forces_{\bP_\alpha} q_t(\alpha)\in\nbQ_\alpha$ (not just in the completion
$\hat{\nbQ}_\alpha$). 
\end{enumerate}
\end{proposition}

\Proof Let $\langle t_i: i<i(*)\rangle$ be an enumeration of $T$ such that
\[(\forall i,j<i(*))(t_i<t_j\quad\Rightarrow\quad i<j).\]
We are proving the proposition by induction on $i(*)$.

\noindent{\sc Case 1:}\ \ \ $i(*)=1$.\\
In this case $T=\{\langle\rangle\}$ and we have to choose $q_{\langle\rangle}$
only, but this is easy, as the set $\{q\rest\rk(\langle\rangle): q\in {\cal
I}\}$ is open dense in $\bP_{\rk(\langle\rangle)}$. 

\noindent{\sc Case 2:}\ \ \ $i(*)=i_0+1>1$.\\
Let $T^*=\{t_i:i<i_0\}$. Then $T^*$ is a standard $(w,\alpha_0)^\gamma$--tree
to which we may apply the inductive hypothesis. Consequently we find $\langle
q^*_t: t\in T^*\rangle\in\ftr'(\bar{\bQ})$ such that for each $t\in T^*$:
\begin{enumerate}
\item $p_t\leq q_t^*\in\{q\rest\rk(t): q\in{\cal I}\}$,\quad and
\item for each $\alpha\in\dom(q_t^*)$, either $q_t^*(\alpha)=p_t(\alpha)$ or
$\forces_{\bP_\alpha} q_t^*(\alpha)\in\nbQ_\alpha$.
\end{enumerate}
Let $q^0=\bigcup\{q^*_s: s<t_{i_0}\}$ (note that $s<t_{i_0}\ \Rightarrow\ s\in
T^*$; easily $q^0\in\bP_{\rk(t_{i_0})}'$). Clearly $q^0$ and $p_{t_{i_0}}$ are
compatible (actually $q^0$ is stronger then the suitable restriction of
$p_{t_{i_0}}$) and therefore we may find a condition $q_{t_{i_0}}\in
\bP_{\rk(t_{i_0})}$ (note: no primes now) such that $q_{t_{i_0}}\in \{q\rest
\rk(t_{i_0}): q\in {\cal I}\}$ and $q_{t_{i_0}}$ stronger than both $q^0$ and
$p_{t_{i_0}}$. Next, for each $t\in T^*$ let
\[q_t\stackrel{\rm def}{=} q_{t_{i_0}}\rest\rk(t\cap t_{i_0})\cup q^*_t\rest
[\rk(t\cap t_{i_0}),\gamma)\geq q^*_t\geq p_t.\]
One easily checks that $\bar{q}=\langle q_t: t\in T\rangle$ is as required. 

\noindent{\sc Case 3:}\ \ \ $i(*)$ is a limit ordinal.\\
Let $\theta=\cf(i(*))$ and let $\langle i_\zeta: \zeta\leq\theta\rangle$ be an
increasing continuous sequence, $i_0=0$, $i_\theta=i(*)$. For $\alpha<
\gamma$, let $\name{x}_\alpha$ be a $\bP_\alpha$--name for a witness that
$\nbQ_\alpha$ is (forced to be) strongly complete for $\baza$ and let $x=
\langle\name{x}_\alpha: \alpha<\gamma\rangle$. Take an $\baza$--complementary
pair $(\bar{N},\bar{a})$ of length $\theta$ such that  $\langle i_\zeta:
\zeta<\theta\rangle,\bar{p},\bar{\bQ},\baza,x,T\in N_0$ and $\|T\|\subseteq
N_0$ (exists as $\baza\in\ck$ is closed: first take a complementary pair of
length $\|T\|^+$ and then restrict it to the interval
$[\|T\|+1,\|T\|+\theta]$).  

By induction on $\zeta\leq\theta$ we define a sequence $\langle \bar{
q}^\zeta:\zeta\leq\theta\rangle$:
\begin{quotation}
\noindent $\bar{q}^\zeta=\langle q^\zeta_t: t\in T\rangle$ is the
$<^*_\chi$--first sequence $\bar{r}=\langle r_t: t\in T\rangle\in\ftr'(\bar{
\bQ})$ such that 
\begin{description}
\item[(i)$_\zeta$] \ \ for every $t\in T$:\ \ $p_t\leq r_t$ and $(\forall \xi<
\zeta)(q^\xi_t\leq r_t)$ and if $\alpha\in\dom(r_t)$, $p_t(\alpha)\neq
r_t(\alpha)$ then $r_t(\alpha)$ is a name for an element of $\nbQ_\alpha$ (not
the completion),
\item[(ii)$_\zeta$] \ if $i_\zeta\leq i<i(*)$ and $\sup\{\rk(t_j): j<i_\zeta\
\&\ t_j<t_i\}\leq \alpha<\rk(t_i)$ then $r_{t_i}(\alpha)=p_{t_i}(\alpha)$,
\item[(iii)$_\zeta$]  if $i<i_\zeta$ then 
\[r_{t_i}\in \bigcap\big\{{\cal J}\in N_\zeta: {\cal J}\subseteq\bP_{\rk(t_i)}
\mbox{ is open dense\/}\big\}.\] 
\end{description}
\end{quotation}
To show that this definition is correct we have to prove that arriving at a
stage $\zeta\leq\theta$ of the construction we may find $\bar{r}$ satisfying
{\bf (i)$_\zeta$}--{\bf (iii)$_\zeta$}. Note that once we know that we may
define $\bar{q}^\xi$ for $\xi\leq\zeta$, we are sure that $\langle\bar{q}^\xi:
\xi\leq\zeta\rangle\in N_{\zeta+1}$ (remember $\bar{N}\rest(\zeta+1)\in
N_{\zeta+1}$). Similarly, arriving at a limit stage $\zeta<\theta$ we are sure
that $\langle\bar{q}^\xi: \xi<\zeta\rangle\in N_{\zeta+1}$. 

\noindent{\sc Stage}\ \ \ $\zeta=0$.\\
Look at $\bar{r}=\bar{p}$: as $i_0=0$, the clause {\bf (iii)$_0$} is empty and
{\bf (i)$_0$}, {\bf (ii)$_0$} are trivially satisfied.

\noindent{\sc Stage}\ \ \ $\zeta=\xi+1$.\\
Let $T^*=\{t_i:i<i_\zeta\}$, $\bar{p}^*=\langle q^\xi_t: t\in T^*\rangle$. We
may apply the inductive hypothesis to $T^*$, $p^*$ and 
\[{\cal I}^*\stackrel{\rm def}{=}\bigcap\big\{{\cal J}\in N_\zeta: {\cal J}
\subseteq\bP_\gamma\mbox{ is open dense\/}\big\}\]
(remember $i_\zeta<i(*)$ and $\bP_\gamma$ does not add new
$<\kappa$--sequences of ordinals, see \ref{nextiter}, so ${\cal I}^*$ is open
dense). Consequently we find $\bar{s}=\langle s_t: t\in T^*\rangle\in
\ftr'(\bar{\bQ})$ such that for each $t\in T^*$
\begin{itemize}
\item $q_t^\xi\leq s_t\in\{q\rest\rk(t): q\in{\cal I}^*\}$,\quad and
\item for each $\alpha\in\dom(s_t)$, either $q_t^\xi(\alpha)=s_t(\alpha)$ or
$\forces_{\bP_\alpha} s_t(\alpha)\in\nbQ_\alpha$. 
\end{itemize}
For $t\in T\setminus T^*$ let $\alpha_t=\sup\{\rk(t_i): i<i_\zeta\ \&\
t_i<t\}$. Note that, for $t\in T\setminus T^*$, $\bigcup\{s_{t_i}: i<i_\zeta\
\&\ t_i<t\}$ is a condition in $\bP_{\alpha_t}'$ stronger than $q^\xi_t\rest
\alpha_t$. So let 
\[r_t=\bigcup\{s_{t_i}: i<i_\zeta\ \&\ t_i<t\}\cup q^\xi_t\rest [\alpha_t,
\gamma)=\bigcup\{s_{t_i}: i<i_\zeta\ \&\ t_i<t\}\cup p_t\rest [\alpha_t,
\gamma)\]
for $t\in T\setminus T^*$ and $r_t=s_t$ for $t\in T^*$. It should be clear
that $\bar{r}=\langle r_t: t\in T\rangle\in \ftr'(\bar{\bQ})$ satisfies the
demands {\bf (i)$_\zeta$}--{\bf (iii)$_\zeta$}. 

\noindent{\sc Stage}\ \ \ $\zeta$ is a limit ordinal.\\
As we noted before, we know that $\langle\bar{q}^\varepsilon:\varepsilon\leq\xi
\rangle\in N_{\xi+1}$ for each $\xi<\zeta$. Hence, as $T\subseteq N_0$
(remember $\|T\|\subseteq N_0$ and $T\in N_0$), we have $\langle
q^\varepsilon_t:\varepsilon\leq\xi\rangle\in N_{\xi+1}$ for each $t\in T$ and
$\xi<\zeta$. Fix $i<i_\zeta$ and let $\xi<\zeta$ be such that $i<i_\xi$. Look
at the sequence $\langle q^\varepsilon_{t_i}:
\xi\leq\varepsilon<\zeta\rangle$. By the choice of $\bar{q}^\varepsilon$ (see
demands {\bf (i)$_\varepsilon$} and {\bf (iii)$_\varepsilon$}) we have that it
is an increasing $(\bar{N}\rest[\xi,\zeta),\bP_{\rk(t_i)})^*$--generic sequence
(note no primes; if we we are not in $\nbQ_\alpha$ then the value is
fixed). By \ref{nextiter} the forcing notion $\bP_{\rk(t_i)}$ is complete for
$\baza$ (and $N_\xi$ contains the witness), so $\langle q^\varepsilon_{t_i}:
\xi\leq\varepsilon<\zeta\rangle$ has an upper bound in
$\bP_{\rk(t_i)}$. Moreover, for each $\alpha<\rk(t_i)$, if $q\in\bP_\alpha$ is
an upper bound of $\langle q^\varepsilon_{t_i}\rest\alpha:\varepsilon<\zeta
\rangle$ then
\[q\forces_{\bP_\alpha}\mbox{``\/the sequence }\langle q^\varepsilon_{t_i}(
\alpha):\varepsilon<\zeta\rangle\mbox{ has an upper bound in
$\nbQ_\alpha$\/''.}\]
Now, for $t\in T$ we may let $\dom(r_t)=\bigcup\limits_{\varepsilon<\zeta}
\dom(q^\varepsilon_t)$ and define inductively $r_t(\alpha)$ for $\alpha\in
\dom(r_t)$ by
\begin{quotation}
\noindent if $(\forall\varepsilon<\zeta)(q^\varepsilon_t(\alpha)=p_t(\alpha))$
then $r_t(\alpha)=p_t(\alpha)$ and otherwise

\noindent $r_t(\alpha)$ is the $<^*_\chi$--first $\bP_\alpha$--name for an
element of $\nbQ_\alpha$ such that 
\[r_t\rest\alpha\forces_{\bP_\alpha}(\forall\varepsilon<\zeta)(
q^\varepsilon_t(\alpha)\leq_{\nbQ_\alpha} r_t(\alpha)).\]
\end{quotation}
It is a routine to check that $\bar{r}=\langle r_t: t\in T\rangle\in\ftr'(
\bar{\bQ})$ and it satisfies {\bf (i)$_\zeta$}--{\bf (iii)$_\zeta$}.

Thus our definition is correct and we may look at the sequence
$\bar{q}^\theta$. Since ${\cal I}\in N_0$ it should be clear that it is as
required. This finishes the inductive proof of the proposition. \QED
\medskip

Our next proposition corresponds to \ref{onesteptree}. However, note that the
meaning of $*$'s is slightly different now. The difference comes from another
type of the game involved and it will be more clear in the proof of theorem
\ref{pseudoiter} below.

\begin{proposition}
\label{CEonestep}
Assume that $\baza\in\ck$ is closed and $\bar{\bQ}=\langle\bP_\alpha,
\nbQ_\alpha:\alpha<\gamma\rangle$ is a $(<\kappa)$--support iteration and
$x=\langle\name{x}_\alpha:\alpha<\gamma\rangle$ is such that 
\[\forces_{\bP_\alpha}\mbox{`` }\nbQ_\alpha\mbox{ is strongly complete for
$\baza$ with witness $\name{x}_\alpha$''}\]
(for $\alpha<\gamma$). Further suppose that 
\begin{description}
\item[($\alpha$)] $(\bar{N},\bar{a})$ is an $\baza$--complementary pair,
$\bar{N}=\langle N_i:i\leq\delta\rangle$, and $x,\baza,\bar{\bQ}\in N_0$,
\item[($\beta$)]  $\tree=(T,<,\rk)\in N_0$ is a standard
$(w,\alpha_0)^\gamma$--tree, $w\subseteq \gamma\cap N_0$, $\|w\|<\cf(\delta)$,
$\alpha_0$ is an ordinal, $\alpha_1=\alpha_0+1$, $0\in w$,
\item[($\gamma$)] $\bar{p}=\langle p_t: t\in T\rangle\in\ftr'(\bar{\bQ})\cap
N_0$, $w\in N_0$, (of course $\alpha_0\in N_0$),  
\item[($\delta$)] $\|N_i\|^{\|w\|+\|T\|}\subseteq N_{i+1}$ for each
$i<\delta$, 
\item[($\varepsilon$)] for $i\leq\delta$, $\tree_i=(T_i,<_i,\rk_i)$ is such
that  

$T_i$ consists of all sequences $t=\langle t_\zeta: \zeta\in\dom(t)\rangle$
such that $\dom(t)$ is an initial segment of $w$, and
\begin{itemize}
\item each $t_\zeta$ is a sequence of length $\alpha_1$,
\item $\langle t_\zeta\rest\alpha_0: \zeta\in\dom(t)\rangle\in T$,
\item for each $\zeta\in\dom(t)$, either $t_\zeta(\alpha_0)=*$ or $t_\zeta(
\alpha_0)\in N_i$ is a $\bP_\zeta$--name for an element of $\nbQ_\zeta$ and 

if $t_\zeta(\alpha)\neq*$ for some $\alpha<\alpha_0$ then $t_\zeta(\alpha_0)
\neq *$, 
\end{itemize}
$\rk_i(t)=\min(w\cup\{\gamma\}\setminus\dom(t))$ and $<_i$ is the extension
relation. 
\end{description}
Then 
\begin{description}
\item[(a)] each $\tree_i$ is a standard $(w,\alpha_1)^\gamma$--tree, $\|T_i\|
\leq\|T\|\cdot\|N_i\|^{\|w\|}$, and if $i<\delta$ then $T_i\in N_{i+1}$,
\item[(b)] $\tree$ is the projection of each $\tree_i$ onto $(w,\alpha_0)$,
\item[(c)] there is $\bar{q}=\langle q_t: t\in T_\delta\rangle\in\ftr'(\bar{
\bQ})$ such that
\begin{description}
\item[(i)]   $\bar{p}\leq_{\proj^{\tree_\delta}_{\tree}}\bar{q}$,
\item[(ii)]  if $t\in T_\delta\setminus\{\langle\rangle\}$ then the condition
$q_t\in\bP_{\rk_\delta(t)}'$ is an upper bound of an $(\bar{N}\rest [i_0,
\delta],\bP_{\rk_\delta(t)})^*$--generic sequence (where $i_0<\delta$ is such
that $t\in T_{i_0}$), and for every $\beta\in\dom(q_t)=N_\delta\cap\rk_\delta
(t)$, $q_t(\beta)$ is a name for the least upper bound in $\hat{\nbQ}_\beta$
of an $(\bar{N}[\name{G}_\beta]\rest [\xi,\delta),\nbQ_\beta)^*$--generic
sequence (for some $\xi<\delta$),  

[Note that, by \ref{bazapres}, the first part of the demand on $q_t$ implies
that if $i_0\leq\xi$ then $q_t\rest\beta$ forces that $(\bar{N}[\name{G
}_\beta]\rest [\xi,\delta],\bar{a}\rest [\xi,\delta])$ is an
$\hat{\family}$--complementary pair.] 
\item[(iii)] if $t\in T_\delta$, $t'=\proj^{\tree_\delta}_{\tree}(t)\in T$,
$\zeta\in\dom(t)$ and $t_\zeta(\alpha_0)\neq *$ then
\[q_t\rest \zeta\forces_{\bP_\zeta}\mbox{``}p_{t'}(\zeta)\leq_{\hat{
\nbQ}_\zeta}t_\zeta(\alpha_0)\ \Rightarrow\ t_\zeta(\alpha_0)\leq_{\hat{
\nbQ}_\zeta} q_t(\zeta) \mbox{''},\]
\item[(iv)]  $q_{\langle\rangle}=p_{\langle\rangle}$.
\end{description}
\end{description}
\end{proposition}

\Proof Clauses (a) and (b) should be clear.\\
(c)\ \ \ One could try to use directly \ref{beforeonestep} for $\bigcap\{{\cal
I}\in N_\delta: {\cal I}\subseteq\bP_\gamma$ open dense\/$\}$ and suitably
``extend'' $\bar{p}$ (see e.g. the successor case below). However, this would
not guarantee the demand {\bf (ii)}. This clause is the reason for the
assumption that $\|w\|<\cf(\delta)$.   

By induction on $i<\delta$ we define a sequence $\langle\bar{q}^i:i<\delta
\rangle$:
\begin{quotation}
\noindent $\bar{q}^i=\langle q^i_t: t\in T_i\rangle$ is the $<^*_\chi$--first
sequence $\bar{r}=\langle r_t: t\in T_i\rangle\in\ftr'(\bar{\bQ})$ such that 
\begin{description}
\item[(i)$_i$] \ \ $\bar{p}\leq_{\proj^{\tree_i}_{\tree}}\bar{r}$ and $(\forall
j<i)(\forall t\in T_j)(q^j_t\leq_{\bP_{\rk(t)}'} r_t)$,
\item[(ii)$_i$] \ if $t=\langle t_\zeta:\zeta\in\dom(t)\rangle\in T_i$ and 
$t'=\proj^{\tree_i}_{\tree}(t)\in T$ then 
\begin{itemize}
\item $(\forall\alpha\in\dom(r_t))(p_{t'}(\alpha)=r_t(\alpha)\ \mbox{ or }\
\forces_{\bP_\alpha}r_t(\alpha)\in\nbQ_\alpha)$, and 
\item $r_t\in\bigcap\{{\cal I}\in N_i: {\cal I}\subseteq\bP_{\rk_i(t)}$ is
open dense\/$\}$\quad and
\item for every $\zeta\in\dom(t)$ such that $t_\zeta(\alpha_0)\neq *$,
\[r_t\rest\zeta\forces_{\bP_\zeta}\mbox{``\/}p_{t'}(\zeta)\leq_{\nbQ_\zeta}
t_\zeta(\alpha_0)\quad\Rightarrow\quad t_\zeta(\alpha_0)\leq_{\nbQ_\zeta}
r_t(\zeta)\mbox{\/'',}\]
\end{itemize}
\item[(iii)$_i$]  $r_{\langle\rangle}=p_{\langle\rangle}$.
\end{description}
\end{quotation}
We have to verify that this definition is correct, i.e. that for each
$i<\delta$ there is an $\bar{r}$ satisfying {\bf (i)$_i$}--{\bf (iii)$_i$}. So
suppose that we arrive to a non-limit stage $i<\delta$ and we have defined
$\langle \bar{q}^j:j<i\rangle$. Note that necessarily $\langle\bar{q}^j:j<i
\rangle\in N_i$ (remember $i$ is non-limit). Let $i=j+1$ and, if $j=-1$, let
$q^{-1}_t=p_{\proj^{\tree_0}_{\tree}(t)}$ for $t\in T_0$ and let
$T_{-1}=\{\langle\rangle\}$. For $t\in T_i$ we define $s_t\in\bP_{\rk_i(t)}$ as
follows.
\begin{description}
\item If $t\in T_j$ then $s_t=q^j_t$.
\item If $t\in T_i\setminus T_j$ and $\zeta^*\in w$ is the first such that
$t\rest (\zeta^*+1)\notin T_j$ then we let $\dom(s_t)=\dom(q^j_{t\rest
\zeta^*})\cup\dom(p_{t'})\cup\dom(t)$, where $t'=\proj^{\tree_i}_{\tree}(t)$. 
Next we define $s_t(\zeta)$ by induction on $\zeta\in\dom(s_t)$:

if $\zeta\in\dom(s_t)\cap\zeta^*$ then $s_t(\zeta)=q^j_{t\rest\zeta^*}(\zeta)$,

if $\zeta\in\dom(t)\setminus\zeta^*$ and $t_\zeta(\alpha_0)\neq *$ then $s_t(
\zeta)$ is the $<^*_\chi$--first $\bP_\zeta$--name for an element of
$\nbQ_\zeta$ such that  
\[s_t\rest\zeta\forces_{\bP_\zeta}\mbox{``\/} p_{t'}(\zeta)\leq s_t(\zeta)
\quad \mbox{ and }\quad p_{t'}(\zeta)\leq t_\zeta(\alpha_0)\ \Rightarrow\
t_\zeta(\alpha_0)\leq s_t(\zeta)\mbox{''}\]
and otherwise it is $p_{t'}(\zeta)$.
\end{description}
It should be clear that $\bar{s}=\langle s_t: t\in T_i\rangle\in\ftr'(\bar{
\bQ})$. Now we apply \ref{beforeonestep} to $T_i$, $\bar{s}$ and 
\[{\cal I}^*\stackrel{\rm def}{=}\bigcap\big\{{\cal I}\in N_i: {\cal
I}\subseteq\bP_\gamma\mbox{ open dense\/}\big\}\]
and we find $\bar{r}=\langle r_t: t\in T_i\rangle\in\ftr'(\bar{\bQ})$ such
that $\bar{s}\leq\bar{r}$ and for each $t\in T_i$
\[r_t\in\{q\rest\rk_i(t)\!: q\in {\cal I}^*\}\quad\mbox{and}\quad (\forall
\alpha\!\in\!\dom(r_t))(s_t(\alpha)\!=\! r_t(\alpha)\ \mbox{or}\
\forces_{\bP_\alpha} r_t(\alpha)\!\in\!\nbQ_\alpha).\]
One easily checks that this $\bar{r}$ satisfies demands {\bf (i)$_i$}--{\bf
(iii)$_i$}.\\
Now suppose that we have successfully defined $\bar{q}^j$ for $j<i$,
$i<\delta$ limit ordinal. Fix $t\in \bigcup\limits_{j<i}T_j$, say $t\in
T_{j_0}$, $j_0<i$. We know that $T_{j_0}\subseteq N_{j_0+1}$ (remember the
assumption $(\delta)$ and the assertion (a)) and that for each $j<i$, $\langle
\bar{q}^\varepsilon:\varepsilon\leq j\rangle\in N_{j+1}$. Consequently,
\[(\forall j\in [j_0,i))(\langle q^\varepsilon_t: j_0\leq\varepsilon\leq j
\rangle\in N_{j+1}).\]
By the demand {\bf (ii)$_\varepsilon$} we have that $\langle q^\varepsilon_t:
j_0\leq \varepsilon<i\rangle$ is an $(\bar{N}\rest [j_0,i],\bP_{
\rk_{j_0}(t)})^*$--generic sequence. As $\bP_{\rk_{j_0}(t)}$ is complete for
$\baza$ (see \ref{nextiter}) and $N_0$ contains all witnesses we conclude that
the sequence $\langle q^\varepsilon_t: j_0\leq\varepsilon <i\rangle$ has an
upper bound in $\bP_{\rk_{j_0}(t)}$. Moreover, if $\alpha<\rk_{j_0}(t)$, and
$q\in\bP_\alpha$ is an upper bound of the sequence $\langle q^\varepsilon_t
\rest\alpha:j_0\leq\varepsilon<i\rangle$ then 
\[q\forces_{\bP_\alpha}\mbox{``}\langle q^\varepsilon_t(\alpha): j_0\leq
\varepsilon <i\rangle\mbox{ has an upper bound in $\nbQ_\alpha$''}\]
(see the proof of \ref{nextiter}). Now we let $\dom(s_t)=\bigcup\{\dom(
q^\varepsilon_t): j_0\leq\varepsilon<i\}$ and we define inductively
\begin{quotation}
\noindent $s_t(\alpha)$ is the $<^*_\chi$--first $\bP_\alpha$--name for an
element of $\nbQ_\alpha$ such that 
\[s_t\rest\alpha\forces_{\bP_\alpha}(\forall\varepsilon\in [j_0,i))(
q^\varepsilon_t(\alpha)\leq_{\nbQ_\alpha} s_t(\alpha)).\]
\end{quotation}
This defines $\bar{s}^0=\langle s_t: t\in \bigcup\limits_{j<i}T_j\rangle$. 
Clearly $\bigcup\limits_{j<i}T_j$ is a standard $(w,\alpha_1)^\gamma$--tree
and $\bar{s}\in\ftr'(\bar{\bQ})$. Now suppose that $t\in T_i\setminus\bigcup
\limits_{j<i} T_j$ and let $\zeta^*$ be the first such that $t\rest\zeta^*
\notin\bigcup\limits_{j<i}T_j$ (so necessarily $\dom(t)\cap\zeta^*$ is cofinal
in $\zeta^*$ and $\cf(\otp(\dom(t)\cap\zeta^*))=\cf(i)$). Then
$\bigcup\{s_{t\rest\zeta}:\zeta<\zeta^*\}\in \bP_{\zeta^*}$. Now define 
\[s_t=\bigcup\{s_{t\rest\zeta}:\zeta<\zeta^*\}\cup p_{t'}\rest
[\zeta^*,\gamma),\]
where $t'=\proj^{\tree_i}_{\tree}(t)$. Note that $\bar{s}=\langle s_t: t\in
T_i\rangle\in\ftr'(\bar{\bQ})$ and if $t\in T_i$, $\alpha\in\dom(s_t)$ then
either $s_t(\alpha)=p_{t'}(\alpha)$ or $\forces_{\bP_\alpha} s_t(\alpha)\in
\nbQ_\alpha$. Now we proceed like in the successor case: we apply
\ref{beforeonestep} to $\bar{s}$, $T_i$ and 
\[{\cal I}^*\stackrel{\rm def}{=}\bigcap\big\{{\cal I}\in N_i: {\cal
I}\subseteq\bP_\gamma\mbox{ open dense\/}\big\}\]
and as a result we get $\bar{r}=\langle r_t: t\in T_i\rangle\in\ftr'(\bar{\bQ}
)$ such that for each $t\in T_i$:
\[\begin{array}{l}
s_t\leq r_t\in\{q\rest\rk_i(t): q\in {\cal I}^*\}\quad\mbox{ and}\\
(\forall\alpha\in\dom(r_t))(s_t(\alpha)= r_t(\alpha)\ \mbox{ or }\
\forces_{\bP_\alpha}r_t(\alpha)\!\in\!\nbQ_\alpha).
  \end{array}\] 
Now one easily checks that $\bar{r}$ satisfies the requirements {\bf
(i)$_i$}--{\bf (iii)$_i$}. 

Thus our definition is the legal one and we have the sequence
$\langle\bar{q}^i: i<\delta\rangle$. We define $\bar{q}=\bar{q}^\delta$
similarly to $\bar{s}$ from the limit stages $i<\delta$, but we replace ``the
$<^*_\chi$--first upper bound in $\nbQ_\alpha$'' by ``the least upper bound in
$\hat{\bQ}_\alpha$''. So suppose $t\in T_\delta$. Since $\|w\|<\cf(\delta)$ we
know that $t\in T_{j_0}$ for some $j_0<\delta$. We declare $\dom(q_t)=\bigcup
\{\dom(q^\varepsilon_t): j_0\leq\varepsilon<\delta\}$ and inductively define
$q_t(\alpha)$ for $\alpha\in\dom(q_t)$:
\begin{quotation}
\noindent $q_t(\alpha)$ is the $<^*_\chi$--first $\bP_\alpha$--name  such that
\[\begin{array}{ll}
q_t\rest\alpha\forces_{\bP_\alpha}&\mbox{``\/$q_t(\alpha)$ is the least upper
bound of the sequence}\\
\ &\ \ \langle q^\varepsilon_t(\alpha): j_0\leq\varepsilon<\delta\rangle\mbox{
in }\hat{\nbQ}_\alpha\mbox{\/''.}
  \end{array}
\] 
\end{quotation}
Like in the limit case of the construction, the respective upper bounds exist,
so $\bar{q}=\langle q_t: t\in T_\delta\rangle$ is well defined. Checking that
it has the required properties is straightforward. \QED

\begin{theorem}
\label{pseudoiter}
Suppose $(\baza_0,\baza_1)\in\ckp$ (so $\baza_0\in\ck$) and $\bar{\bQ}=\langle
\bP_\alpha,\nbQ_\alpha:\alpha<\gamma\rangle$ is a $(<\kappa)$--support
iteration such that for each $\alpha<\kappa$
\[\forces_{\bP_\alpha}\mbox{``\/$\nbQ_\alpha$ is complete for $(\baza_0,
\baza_1)$''}\]
Then 
\begin{description}
\item[(a)] $\forces_{\bP_\gamma}(\baza_0,\baza_1)\in\ckp$, moreover
\item[(b)] $\bP_\gamma$ is complete for $(\baza_0,\baza_1)$. 
\end{description}
\end{theorem}

\Proof We need only part (a) of the conclusion, so we concentrate on it.
Let $\chi$ be a large enough regular cardinal, $\name{x}$ be a name for
an element of ${\cal H}(\chi)$ and $p\in\bP_\gamma$. Let $\name{x}_\alpha$ be
a $\bP_\alpha$--name for the witness that $\nbQ_\alpha$ is (forced to be)
complete for $(\baza_0,\baza_1)$, and let $\bar{x}=\langle\name{x}_\alpha:
\alpha<\gamma\rangle$. Since $(\baza_0,\baza_1)\in\ckp$ we find $\bar{M}=
\langle M_i: i\leq\delta\rangle$ which is ruled by $(\baza_0,\baza_1)$ with an
$\baza_0$--approximation $\langle\bar{N}^i: -1\leq i<\delta\rangle$ and such
that $p,\bar{\bQ},\name{x},\bar{x},\baza_0,\baza_1\in M_0$ (see \ref{ruler}). 
Let $\bar{N}^i=\langle N^i_\varepsilon:\varepsilon\leq\delta_i\rangle$ and
$\bar{a}^i\in\baza_0$ be such that $(\bar{N}^i,\bar{a}^i)$ is an
$\baza_0$--complementary pair. Let $w_i=\{0\}\cup\bigcup\limits_{j<i}(\gamma
\cap M_j)$ (for $i\leq\delta$). By the demands of \ref{ruler} we know that
$\|w_i\|<\cf(\delta_i)$.  

By induction on $i\leq\delta$ we define standard $(w_i,i+1)^\gamma$--trees
$\tree_i\in M_{i+1}$ and $\bar{p}^i=\langle p^i_t: t\in T_i\rangle\in
\ftr'(\bar{\bQ})\cap M_{i+1}$ such that $\|T_i\|\leq \|M_{i}\|^{\|w_i\|}=
\|M_{i+1}\|$, and if $j<i\leq\delta$ then $\tree_j=\proj^{(w_i,i+1)}_{(w_j,
j+1)}(\tree_i)$ and $\bar{p}^j\leq_{\proj^{\tree_i}_{\tree_j}}\bar{p}^i$. 

\noindent{\sc Case 1:}\ \ \ $i=0$.\\
Let $T^*_0$ consists of all sequences $\langle t_\zeta:\zeta\in\dom(t)\rangle$
such that $\dom(t)$ is an initial segment of $w_0$ and $t_\zeta=\langle
\rangle$ for $\zeta\in\dom(t)$. Thus $T^*_0$ is a standard $(w_0,0
)^\gamma$--tree, $\|T^*_0\|=\|w_0\|$ . For $t\in T^*_0$ let $p^{*0}_t= 
p\rest \rk^*_0(t)$. Clearly the sequence $\bar{p}^{*0}=\langle p^{*0}_t: t\in
T^*_0\rangle$ is in $\ftr'(\bar{\bQ})\cap N^{-1}_0$. Apply \ref{CEonestep} to
$\baza_0,\bar{\bQ},\bar{N}^{-1},\tree^*_0,w_0$  and $\bar{p}^{*0}$ (note that
$\|N^{-1}_\varepsilon\|^{\|w_0\|}\leq \|N^{-1}_\varepsilon\|^{\|M_0\|}$ for
$\varepsilon<\delta_0$). As a result we get a $(w_0,1)^\gamma$--tree $\tree_0$
(the one called $\tree_{\delta_0}$ there) and $\bar{p}^0=\langle p^0_t: t\in
T_0\rangle\in\ftr'(\bar{\bQ})\cap M_1$ (the one called $\bar{q}$ there)
satisfying clauses \ref{CEonestep}($\varepsilon$), \ref{CEonestep}(c)(i)--(iv)
and such that $\|T_0\|\leq\|N^{-1}_{\delta_0}\|^{\|w_0\|}=\|M_0\|^{\|w_0\|}=
\|M_0\|$ (remember $\cf(\delta_0)>2^{\|M_0\|}$). So, in particular, if $t\in
T_0$, $\zeta\in\dom(t)$ then $t_\zeta(0)\in M_1$ is either $*$ or a
$\bP_\zeta$--name for an element of $\nbQ_\zeta$.\\
Moreover, we additionally require that $(\tree_0,\bar{p}^0)$ is the
$<^*_\chi$--first with all these properties, so $\tree_0,\bar{p}^0\in M_1$.

\noindent{\sc Case 2:}\ \ \ $i=i_0+1$.\\
We proceed similarly to the previous case. Suppose we have defined
$\tree_{i_0}$ and $\bar{p}^{i_0}$ such that $\tree_{i_0},\bar{p}^{i_0}\in
M_{i_0+1}$, $\|T_{i_0}\|\leq\|M_{i_0+1}\|$. Let $\tree^*_i$ be a standard
$(w_i,i)^\gamma$--tree such that 
\begin{quotation}
\noindent $T^*_i$ consists of all sequences $\langle t_\zeta:\zeta\in\dom(t)
\rangle$ such that $\dom(t)$ is an initial segment of $w_i$ and 
\[\hspace{-0.5cm}\langle t_\zeta:\zeta\in\dom(t)\cap w_{i_0}\rangle\in
T_{i_0}\ \mbox{ and }\ (\forall \zeta\in\dom(t)\setminus w_{i_0})(\forall
j\leq i_0)(t_\zeta(j)=*).\] 
\end{quotation}
Thus $\tree_{i_0}=\proj^{(w_i,i)}_{(w_{i_0},i_0)}(\tree^*_i)$ and $\|T^*_i\|
\leq\|M_i\|$. Let $p^{*i}_t=p^{i_0}_{t'}\rest\rk^*_i(t)$ for $t\in T^*_i$,
$t'=\proj^{\tree_i}_{\tree_{i_0}}(t)$. Now apply \ref{CEonestep} to $\baza_0$,
$\bar{\bQ}$, $\bar{N}^{i_0}$, $\tree^*_i$, $w_i$ and $\bar{p}^{*i}$ (check that
the assumptions are satisfied). So we get a standard $(w_i,i+1)^\gamma$--tree
$\tree_i$ and a sequence $\bar{p}^i$ satisfying \ref{CEonestep}($\varepsilon$),
\ref{CEonestep}(c)(i)--(iv), and we take the $<^*_\chi$--first pair $(\tree_i,
\bar{p}^i)$ with these properties. In particular we will have $\|T_i\|\leq
\|M_{i_0}\|\cdot\|N^{i_0}_{\delta_i}\|^{\|M_{i_0}\|}=\|M_{i_0+1}\|$, and
$\bar{p}^i,\tree_i \in M_{i+1}$.
  
\noindent{\sc Case 3:}\ \ \ $i$ is a limit ordinal.\\
Suppose we have defined $\tree_j$, $\bar{p}^j$ for $j<i$ and we know that
$\langle (\tree_j,\bar{p}^j): j<i\rangle\in M_{i+1}$ (this is the consequence
of taking ``the $<^*_\chi$--first such that$\ldots$''). Let $\tree^*_i=\inver(
\langle \tree_j: j<i\rangle)$. Now, for $t\in T^*_i$ we would like to define
$p^{*i}_t$ as the limit of $p^j_{\proj^{\tree^*_i}_{\tree_j}(t)}$. However,
our problem is that we do not know if the limit exists. Therefore we restrict
ourselves to these $t$ for which the respective sequence has an upper
bound. To be more precise, for $t\in \tree^*_i$ we apply the following
procedure. 
\begin{description}
\item[$(\otimes)$] Let $t^j=\proj^{\tree^*_i}_{\tree_j}(t)$ for $j<i$. Try to
define inductively a condition $p^{*i}_t\in\bP_{\rk^*_i(t)}$ such that
$\dom(p^{*i}_t)=\bigcup\{\dom(p^j_{t^j})\cap\rk^*_i(t): j<i\}$. Suppose we
have successfully defined $p^{*i}_t\rest\alpha$, $\alpha\in\dom(p^{*i}_t)$, in
such a way that $p^{*i}_t\rest\alpha\geq p^j_{t^j}\rest\alpha$ for all
$j<i$. We know that
\[p^{*i}_t\rest\alpha\forces_{\bP_\alpha}\mbox{``\/the sequence }\langle p^j_{
t^j}(\alpha): j<i\rangle\mbox{ is
$\leq_{\hat{\nbQ}_\alpha}$--increasing\/''.}\] 
So now, if there is a $\bP_\alpha$--name $\name{\tau}$ for an element of
$\nbQ_\alpha$ such that 
\[p^{*i}_t\rest\alpha\forces_{\bP_\alpha}(\forall j<i)(p^j_{t^j}(\alpha)
\leq_{\hat{\nbQ}_\alpha}\name{\tau}),\]
then we take the $<^*_\chi$--first such a name as $p^{*i}_t(\alpha)$, and we
continue. If there is no such $\name{\tau}$ then we decide that $t\notin
\tree^+_i$ and we stop the procedure.
\end{description}
Now, let $\tree^+_i$ consist of these $t\in T^*_i$ for which the above
procedure resulted in a successful definition of $p^{*i}_t\in\bP_{\rk^*_i(t)
}$. It might be not clear at the moment if $T^+_i$ contains anything more than
$\langle\rangle$, but we will see that this is the case. Note that 
\[\|T^+_i\|\leq\|T^*_i\|\leq\prod_{j<i}\|T_j\|\leq\prod_{j<i}\|M_j\|\leq
2^{\|M_i\|}\leq \|N^i_2\|.\]
Moreover, for $\varepsilon>2$ we have $\|N^i_\varepsilon\|^{\|w_i\|+\|T^+_i\|}
\leq \|N^i_\varepsilon\|^{\|N^i_2\|}\subseteq N^i_{\varepsilon+1}$ and
$\tree^+_i, \bar{p}^{*i}\in M_{i+1}$. Let $\tree_i=\tree^*_i$, $\bar{p}^i=
\bar{p}^{*i}$ (this time no need to take the $<^*_\chi$--first pair as the
process leaves no freedom). 

After the construction is carried out we continue in a similar manner as in
\ref{iterweak} (but note slightly different meaning of the $*$'s here).\\
So we let $\tree_\delta=\inver(\langle\tree_i: i<\delta\rangle)$. It is a
standard $(w_\delta,\delta)^\gamma$--tree. By induction on $\alpha\in
w_\delta\cup\{\gamma\}$ we choose $q_\alpha\in\bP_\alpha'$ and
a $\bP_\alpha$--name $\name{t}_\alpha$ such that
\begin{description}
\item[(a)] $\forces_{\bP_\alpha}$``$\name{t}_\alpha\in T_\delta\ \&\
\rk_\delta(\name{t}_\alpha)=\alpha$'', and $i^\alpha_0=\min\{i<\delta:\alpha
\in M_i\}<\delta$, 
\item[(b)] $\forces_{\bP_\alpha}$``$\name{t}_\beta=\name{t}_\alpha\rest
\beta$'' for $\beta<\alpha$,
\item[(c)] $\dom(q_\alpha)=w_\delta\cap\alpha$,
\item[(d)] if $\beta<\alpha$ then $q_\beta=q_\alpha\rest \beta$,
\item[(e)] $p^i_{\proj^{\tree_\delta}_{\tree_i}(\name{t}_\alpha)}$ is well
defined and $p^i_{\proj^{\tree_\delta}_{\tree_i}(\name{t}_\alpha)}\rest\alpha
\leq q_\alpha$ for each $i<\delta$, 
\item[(f)] for each $\beta<\alpha$
\[\begin{array}{ll}
q_\alpha\forces_{\bP_\alpha}&\mbox{``\/}(\forall i<\delta)((\name{t}_{\beta+1}
)_\beta(i)=*\ \Leftrightarrow\ i<i^\beta_0)\mbox{ and the sequence}\\
\ &\langle i^\beta_0, p^{\name{i}^\beta_0}_{\proj^{\tree_\delta}_{\tree_{
i^\beta_0}}(\name{t}_{\beta+1})}(\beta),\langle(\name{t}_{\beta+1})_\beta(i),
p^i_{\proj^{\tree_\delta}_{\tree_i}(\name{t}_{\beta+1})}(\beta):\ i^\beta_0
\leq i<\delta\rangle\rangle\\
\ &\mbox{is a result of a play of the game }{\cal G}^\spadesuit_{\bar{M}[
\name{G}_\beta], \langle \bar{N}^i[\name{G}_\beta]: i<\delta\rangle}(
\nbQ_\beta,{\bf 0}_{\nbQ_\beta}),\\
\ &\mbox{won by player $\com$'',}
\end{array} \]
\item[(g)] the condition $q_\alpha$ forces (in $\bP_\alpha$) that

``the sequence $\bar{M}[\name{G}_{\bP\alpha}]\rest [i_\alpha,\delta]$ is ruled
by $(\baza_0,\baza_1)$ and $\langle\bar{N}^i[\name{G}_{\bP_\alpha}]:
i^\alpha_0\leq i<\delta\rangle$ is its $\baza_0$--approximation''.
\end{description}
(Remember: $\hat{\baza}_1$ is closed under end segments.) This is done
completely parallely to the last part of the proof of \ref{iterweak}. 

Finally look at the condition $q_\gamma$ and the clause {\bf (g)} above. \QED

\begin{proposition}
\label{easycc}
Suppose that $\mu^*=\kappa$, $\baza\in\ck$ is closed and $\bar{\bQ}=\langle
\bP_\alpha,\nbQ_\alpha: \alpha<\gamma\rangle$ is a $(<\kappa)$--support
iteration such that for each $\alpha<\gamma$
\[\forces_{\bP_\alpha}\mbox{``\/$\nbQ_\alpha$ is strongly complete for $\baza$
and $\|\nbQ_\alpha\|\leq\kappa$\/''.}\]
Then $\bP_\gamma$ satisfies $\kappa^+$--cc (even more: it satisfies the
$\kappa^+$--Knaster condition). 
\end{proposition}

\Proof For $\alpha<\gamma$ choose $\bP_\alpha$--names $\name{x}_\alpha$ and
$\name{h}_\alpha$ such that
\[\begin{array}{ll}
\forces_{\bP_\alpha}&\mbox{``\/$\name{x}_\alpha$ witnesses that $\nbQ_\alpha$
is complete for $\baza$ and}\\ 
\ &\ \/\name{h}_\alpha:\nbQ_\alpha\stackrel{1 - 1}{\longrightarrow}\kappa
\mbox{ is one--to--one''.}
  \end{array}\]
Since $\baza\in\ck$, for each $p\in\bP_\gamma$ we find an
$\baza$--complementary pair $(\bar{N}^p,\bar{a}^p)$ such that $\bar{N}^p= 
\langle N^p_i:i\leq\omega\rangle$ and $p,\bar{\bQ},\baza,\langle
\name{x}_\alpha:\alpha<\gamma\rangle, \langle\name{h}_\alpha: \alpha<\gamma
\rangle\in N_0^p$. Next choose an increasing sequence $\bar{q}^p=\langle
q^p_i:i<\omega\rangle$ of conditions from $\bP_\gamma$ such that for each
$i<\omega$: 
\begin{description}
\item[$(\alpha)$] $p\leq q^p_0$,\quad $\bar{q}^p\rest (i+1)\in N^p_{i+1}$,
\item[$(\beta)$]  $q^p_i\in\bigcap\{{\cal I}\in N_i: {\cal I}\subseteq
\bP_\gamma$ open dense\/$\}$.
\end{description}
[Why possible? Remember \ref{nextiter} and particularly \ref{cl9}.] So the
condition $q^p_i$ is generic over $N_i$ and therefore it decides the values of
$\name{h}_\alpha(q^p_j(\alpha))$ for each $j<i$, $\alpha\in \dom(q^p_j)$
(remember: if $j<i$ then $q^p_j\in N_i$ and thus $\dom(q^p_j)\subseteq
N_i$). Let $\varepsilon^{p,\alpha}_j<\kappa$ be such that for each $i>j$
(remember $\bar{q}^p$ is increasing)
\[q^p_i\rest\alpha\forces_{\bP_\alpha} \name{h}_\alpha(q^p_j(\alpha))=
\varepsilon^{p,\alpha}_j.\]
Suppose now that $\langle p^\zeta:
\zeta<\kappa^+\rangle\subseteq\bP_\gamma$. For $\zeta<\kappa^+$ let
$A_\zeta=\bigcup\limits_{i<\omega}\dom(q^{p^\zeta}_i)$ (so $A_\zeta\in
[\gamma]^{<\kappa}$). Applying the $\Delta$--system lemma (remember $\kappa$
is strongly inaccessible) we find ${\cal X}\subseteq \kappa^+$, $\|{\cal X}
\|=\kappa^+$ such that $\{A_\zeta:\zeta\in {\cal X}\}$ forms a
$\Delta$--system and for each $\zeta,\xi\in {\cal X}$: 
\begin{itemize}
\item $\|A_\zeta\|=\|A_\xi\|$,
\item if $\alpha\in A_\zeta\cap A_\xi$ then
\[\min\{i<\omega:\alpha\in\dom(q^{p^\zeta}_i)\}=\min\{i<\omega:\alpha\in
\dom(q^{p^\xi}_i)\}=i_0,\]
and for each $i<\omega$
\[\otp(\alpha\cap\dom(q^{p^\zeta}_i))=\otp(\alpha\cap\dom(q^{p^\xi}_i))\quad
\mbox{ and }\quad\varepsilon^{p^\zeta,\alpha}_i=\varepsilon^{p^\xi,\alpha}_i\]
(the last for $i\geq i_0$).
\end{itemize}
We are going to show that for each $\xi,\zeta\in {\cal X}$ the conditions
$p^\zeta,p^\xi$ are compatible. To this end we define a common upper bound $r$
of $p^\zeta,p^\xi$. First we declare that
\[\dom(r)=A_\zeta\cup A_\xi\]
and then we inductively define $r(\alpha)$ for $\alpha\in\dom(r)$:
\begin{quotation}
\noindent if $\alpha\in A_\zeta$ then $r(\zeta)$ is a $\bP_\alpha$--name such
that
\[\begin{array}{ll}
r\rest\alpha\forces_{\bP_\alpha}&\mbox{``\/$r(\alpha)$ is the upper bound of
}\langle q^{p^\zeta}_i(\alpha): i<\omega\rangle\\
\ &\ \/\mbox{with the minimal value of }\name{h}_\alpha(r(\alpha))\mbox{''}
  \end{array}\]
and otherwise (i.e. if $\alpha\in A_\xi\setminus A_\zeta$) it is a
$\bP_\alpha$--name such that
\[\begin{array}{ll}
r\rest\alpha\forces_{\bP_\alpha}&\mbox{``\/$r(\alpha)$ is the upper bound of
}\langle q^{p^\xi}_i(\alpha): i<\omega\rangle\\
\ &\ \/\mbox{with the minimal value of }\name{h}_\alpha(r(\alpha))\mbox{''.}
  \end{array}\]
\end{quotation}
By induction on $\alpha\in\dom(r)\cup\{\gamma\}$ we show that 
\[q^{p^\zeta}_i\rest\alpha\leq_{\bP_\alpha} r\rest\alpha\quad\mbox{and}\quad
q^{p^\xi}_i\rest\alpha\leq_{\bP_\alpha} r\rest\alpha\quad\mbox{ for all
}i<\omega.\]
Note that, by \ref{bazapres}, this implies that the respective upper bounds
exist and thus $r(\alpha)$ is well defined then. There is nothing to do at
non-successor stages, so suppose that we have arrived to a stage
$\alpha=\beta+1$.\\ 
If $\beta\in A_\zeta$ then, by the definition of $r(\beta)$, we have
\[r\rest\beta\forces_{\bP_\beta}(\forall i<\omega)(q^{p^\zeta}_i(\beta)\leq
r(\beta)).\]
Similarly if $\beta\in A_\xi\setminus A_\zeta$ and we consider
$q^{p^\xi}_i(\beta)$. Trivially, no problems can happen if $\beta\in
A_\zeta\setminus A_\xi$ and we consider $q^{p^\xi}_i(\beta)$ or if $\beta\in
A_\xi\setminus A_\zeta$ and we consider $q^{p^\zeta}_i(\beta)$. So the only
case we may worry about is that $\beta\in A_\zeta\cap A_\xi$ and we want to
show that $r(\beta)$ is (forced to be) stronger than all
$q^{p^\xi}_i(\beta)$. But note: by the inductive hypothesis we know that
$r\rest\beta$ is an upper bound to both $\langle q^{p^\xi}_i\rest\beta:
i<\omega\rangle$ and $\langle q^{p^\zeta}_i\rest\beta:i<\omega\rangle$ and
therefore 
\[r\rest\beta\forces_{\bP_\beta}\mbox{`` }\name{h}_\beta(q^{p^\xi}_i(\beta))
=\varepsilon^{p^\xi,\beta}_i\ \ \&\ \ \name{h}_\beta(q^{p^\zeta}_j(\beta))=
\varepsilon^{p^\zeta,\beta}_j\mbox{ ''},\]
whenever $i,j<\omega$ are such that $\beta\in\dom(q^{p^\xi}_i)$,
$\beta\in\dom(q^{p^\zeta}_j)$. But now, by the choice of ${\cal X}$ we have:
\[\beta\in\dom(q^{p^\xi}_i)\quad\Leftrightarrow\quad\beta\in\dom(q^{p^\zeta
}_i),\qquad\mbox{ and } \varepsilon^{p^\zeta,\beta}_i=\varepsilon^{p^\xi,
\beta}_i.\]
Since $\name{h}_\beta$ is (forced to be) a one--to--one function, we conclude
that 
\[r\rest\beta\forces_{\bP_\beta}(\forall i<\omega)(q^{p^\zeta}_i(\beta)=
q^{p^\xi}_i(\beta)),\]
so taking care of the $\zeta$'s side we took care of the $\xi$'s side as
well. This finishes the proof of the proposition. \QED

\section{The Axiom and its applications}
\label{aksitsapp}
\begin{definition}
Suppose that $(\baza_0,\baza_1)\in\ckp$ and $\theta$ is a regular
cardinal. Let $\Axkt(\baza_0,\baza_1)$, {\em the forcing axiom for $(\baza_0,
\baza_1)$ and $\theta$}, be the following sentence:
\begin{quotation}
\noindent If $\bQ$ is a complete for $(\baza_0,\baza_1)$ forcing notion of
size $\leq\kappa$ and $\langle {\cal I}_i: i<i^*<\theta\rangle$ is a sequence
of dense subsets of $\bQ$,

\noindent then there exist a directed set $H\subseteq\bQ$ such that
\[(\forall i<i^*)(H\cap {\cal I}_i\neq\emptyset).\]
\end{quotation}
\end{definition}

\begin{theorem}
\label{getaxiom}
Assume that $\mu^*=\kappa$, $(\baza_0,\baza_1)\in\ckp$ and 
\[\kappa<\theta=\cf(\theta)\leq\mu=\mu^\kappa.\]
Then there is a strongly complete for $\baza_0$ forcing notion $\bP$ of
cardinality $\mu$ such that 
\begin{description}
\item[($\alpha$)] $\bP$ satisfies the $\kappa^+$--cc,
\item[($\beta$)]  $\forces_{\bP}(\baza_0,\baza_1)\in\ckp$ and even more:
\item[($\beta^+$)] if $\baza_1^*\subseteq\baza_1$ is such that $(\baza_0,
\baza^*_1)\in\ckp$ then $\forces_{\bP}(\baza_0,\baza_1^*)\in\ckp$,
\item[($\gamma$)] $\forces_{\bP}\Axkt(\baza_0,\baza_1)$.
\end{description}
\end{theorem}

\Proof The forcing notion $\bP$ will be the limit of a $(<\kappa)$--support
iteration $\langle\bP_\alpha,\nbQ_\alpha: \alpha<\alpha^*\rangle$ (for some
$\alpha^*<\mu^+$) such that
\begin{description}
\item[(a)] for each $\alpha<\alpha^*$
\[\forces_{\bP_\alpha}\mbox{``\/$\nbQ_\alpha$ is a partial order on $\kappa$
complete for $(\baza_0,\baza_1)$\/''}.\]
\end{description}
By \ref{easycc} we will be sure that $\bP=\bP_{\alpha^*}$ satisfies
$\kappa^+$--cc. Applying \ref{pseudoiter} we will see that $\forces_{\bP_{
\alpha^*}}(\baza_0,\baza_1)\in\ckp$ (also $\bP_{\alpha^*}$ is complete for
$(\baza_0,\hat{\baza}_1)$).  The iteration $\langle\bP_\alpha,\nbQ_\alpha:
\alpha<\alpha^*\rangle$ will be built by a bookkeeping argument, but we do not
determine in advance its length $\alpha^*$.

Before we start the construction, note that if $\bQ$ is a $\kappa^+$--cc
forcing notion of size $\leq\mu$ then there are at most $\mu$ $\bQ$--names for
partial orders on $\kappa$ (up to isomorphism). Why? Remember $\mu^\kappa=\mu$
and each $\bQ$--name for a poset on $\kappa$ is described by a
$\kappa$--sequence of maximal antichains of $\bQ$. By a similar argument we
will know that each $\bP_\alpha$ has a dense subset of size $\leq\mu$ (for
$\alpha\leq\alpha^*$). Consequently there are, up to an isomorphism, at most
$\mu$ $\bP_{\alpha^*}$--names for partial orders on $\kappa$.   

Let ${\frak K}$ consist of all $(<\kappa)$--support iterations $\bar{\bQ}= 
\langle\bP_\alpha,\nbQ_\alpha: \alpha<\alpha_0\rangle$ of length $<\mu^+$
satisfying the demand {\bf (a)} above (with $\alpha_0$ in place of
$\alpha^*$). Elements of ${\frak K}$ are naturally ordered by
\[\bar{\bQ}^0\leq_{\frak K}\bar{\bQ}^1\quad\mbox{ if and only if }\quad
\bar{\bQ}^0=\bar{\bQ}^1\rest \lh(\bar{\bQ}^0).\]
Note that every $\leq_{\frak K}$--increasing sequence of length $<\mu^+$ has
the least upper bound in $({\frak K},\leq_{\frak K})$. By what we said before,
we know that if $\langle \bP_\alpha,\nbQ_\alpha:\alpha<\alpha_0\rangle\in
{\frak K}$, then $\bP_{\alpha_0}$ contains a dense subset of size $\leq\mu$,
satisfies $\kappa^+$--cc and forces that $(\baza_0,\baza_1)\in\ckp$. Moreover, 
\begin{description}
\item[$(\circledast_{\frak K})$] if $\bar{\bQ}^0=\langle \bP_\alpha^0,
\nbQ_\alpha^0:\alpha<\alpha_0\rangle\in{\frak K}$ and $\nbQ$ is a
$\bP_{\alpha_0}^0$--name for a forcing notion on $\kappa$ then
\begin{description}
\item[$(\oplus_1)$] {\em either\/} there is no $\bar{\bQ}^1=\langle
\bP_\alpha^1,\nbQ_\alpha^1: \alpha<\alpha_1\rangle\in {\frak K}$ such that
$\bar{\bQ}^0\prec_{\frak K}\bar{\bQ}^1$ and
\[\forces_{\bP^1_{\alpha_1}}\mbox{``\/$\nbQ$ is complete for $(\baza_0,
\baza_1)$\/''}\]
\item[$(\oplus_2)$] {\em or\/} there is $\bar{\bQ}^1=\langle \bP_\alpha^1,
\nbQ_\alpha^1: \alpha<\alpha_1\rangle\in {\frak K}$ such that $\bar{\bQ}^0
\prec_{\frak K}\bar{\bQ}^1$ and
\[\begin{array}{ll}
\forces_{\bP^1_{\alpha_1}}&\mbox{``\/there is a directed set $H\subseteq\nbQ$
which meets all}\\
\ &\mbox{\ \/dense subsets of $\nbQ$ from }\V^{\bP^0_{\alpha_0}}\mbox{\/''.}
  \end{array}\]
\end{description}
\end{description}
[Why? Suppose that $(\oplus_1)$ fails and it is exemplified by
$\bar{\bQ}^1$. Take $\bar{\bQ}^1*\nbQ$.]\\
Consequently, as ${\frak K}$ is closed under increasing $<\mu^+$--sequences,
we have 
\begin{description}
\item[$(\circledast_{\frak K}^+)$] for every $\bar{\bQ}\in {\frak K}$ there is 
$\bar{\bQ}^0=\langle \bP_\alpha^0,\nbQ_\alpha^0:\alpha<\alpha_0\rangle\in 
{\frak K}$ such that $\bar{\bQ}\prec_{\frak K}\bar{\bQ}^0$ and for every 
$\Lim(\bar{\bQ})$--name $\nbQ$ for a forcing notion on $\kappa$ one of the
following conditions occurs: 
\begin{description}
\item[$(\oplus_1)$] there is no $\bar{\bQ}^1=\langle\bP_\alpha^1,
\nbQ_\alpha^1: \alpha<\alpha_1\rangle\in {\frak K}$ such that $\bar{\bQ}^0
\prec_{\frak K}\bar{\bQ}^1$ and
\[\forces_{\bP^1_{\alpha_1}}\mbox{``\/$\nbQ$ is complete for $(\baza_0,
\baza_1)$\/''}\]
\item[$(\oplus_2^+)$] $\forces_{\bP^0_{\alpha_0}}$`` there is a directed set
$H\subseteq\nbQ$ which meets all dense subsets of $\nbQ$ from $\V^{\Lim(
\bar{\bQ})}$ ''.
\end{description}
\end{description}
[Why? Remember that there is at most $\mu$ $\Lim(\bar{\bQ})$--names for
partial orders on $\kappa$.]

Using these remarks we may build our iteration in the following way. We choose
a $\prec_{\frak K}$--increasing continuous sequence $\langle \bar{\bQ}^\zeta:
\zeta\leq\theta^+\rangle\subseteq {\frak K}$ such that 
\begin{description}
\item[(b)] for every $\zeta<\theta^+$, $\bar{\bQ}^{\zeta+1}$ is given by
$(\circledast^+_{\frak K})$ for $\bar{\bQ}^\zeta$.
\end{description}
Now it is a routine to check that $\bP=\bP^{\theta^+}_{\alpha_{\theta^+}}$ is
as required. \QED 
\bigskip

In \ref{weacom} below remember about our main case: $S^*\subseteq\kappa$ is
stationary co-stationary and $\baza_0$ consists of all increasing continuous
sequences $\bar{a}=\langle a_i: i\leq\alpha\rangle$ such that $a_i\in\kappa
\setminus S^*$ (for $i\leq\alpha$). In this case the forcing notion $\bR$ is
the standard way to make the set $S^*$ non-stationary (by adding a club of
$\kappa$; a condition gives an initial segment of the club). Since forcing
with $\bR$ preserves stationarity of subsets of $\kappa\setminus S^*$, the
conclusion of \ref{weacom} below gives us
\begin{description}
\item[($*$)] in $\V^{\Lim{\bar{\bQ}}}$, every stationary set $S\subseteq\kappa
\setminus S^*$ reflects in some inaccessible.
\end{description}

\begin{proposition}
\label{weacom}
Suppose that $(\baza_0,\baza_1)\in\ckp$ (so $\baza_0\in\ck$), $\mu^*=\kappa$
(for simplicity) and $\bar{\bQ}=\langle \bP_\alpha,\nbQ_\alpha:\alpha<\gamma
\rangle$ is a $(<\kappa)$--support iteration such that for each
$\alpha<\kappa$ 
\[\forces_{\bP_\alpha}\mbox{``\/$\nbQ_\alpha$ is complete for $(\baza_0,
\baza_1)$ and $\|\nbQ_\alpha\|\leq\kappa$.\/''}\]
Further assume that:
\begin{description}
\item[(a)] $\baza_0$ is reasonably closed: it is closed under subsequences and
if $\bar{a}=\langle a_i:i\in\delta\rangle\in\baza_0$ and $\bar{b}^i=\langle
b^i_\alpha:\alpha\leq\alpha_i\rangle\in\baza_0$ are such that $b^i_0=a_i$,
$b^i_{\alpha_i}=a_{i+1}$ (for $i<\delta$) then the concatenation of all
$\bar{b}^i$ (for $i<\delta$) [e.g.~$\baza_0$ is derived from $S\subseteq
\kappa$ like in \ref{forcing}],
\item[(b)] $\bR=(\baza_0,\vartriangleleft)$,
\item[(c)] in $\V^{\bR}$ and even in $\V^{\bR*\cohen_\kappa}$, $\kappa$ is a
weakly compact cardinal (or just: stationary subsets of $\kappa$ reflect in
inaccessibles).  
\end{description}
Then, in $\V^{\bP_\gamma*\bR}$, $\kappa$ is weakly compact (or just:
stationary subsets of $\kappa$ reflect in inaccessibles).  
\end{proposition}

\Proof First note that the forcing with $\bR$ does not add new sequences of
length $<\kappa$ of ordinals. [Why? Suppose that $\name{x}$ is an $\bR$--name
for a function from $\theta$ to $\V$, $\theta<\kappa$ is a regular cardinal
and $r\in\bR$. Take an $\baza_0$--complementary pair $(\bar{N},\bar{a})$ such
that $\bar{N}=\langle N_i:i\leq\theta\rangle$ and $r,\name{x}\in N_0$ and the
error is, say, $n$. Now build inductively an increasing sequence $\langle r_i:
i\leq\theta\rangle\subseteq\bR$ such that for every $i\leq\theta$:
\begin{itemize}
\item $r_0=r$, the condition $r_{i+1}$ decides the value of $\name{x}(i)$,
\item if $i=\gamma+k+1$, $\gamma$ is a non-successor, $k<\omega$ then
$r_i\in N_{\gamma+(2k+2)(n+1)}$ and if $r_i=\langle a^i_\xi:\xi\leq\alpha_i
\rangle$ then $a^i_{\alpha_i}=a_{\gamma+(2k+1)(n+1)}$,
\item if $i<\theta$ is limit then $\langle r_j:j<i\rangle\in N_{i+1}$ and
$r_i$ is the least upper bound of $\langle r_j:j<i\rangle$ (so $r_i\in
N_{i+1}$).
\end{itemize}
The construction is straightforward. If we have defined $r_i\in N_{\gamma+(
2k+2)(n+1)}$ then we first take the $<_\chi$--first condition $r^*_i=\langle
a^*_\xi:\xi\leq\alpha^*\rangle$ stronger than $r_i$ and deciding the value of
$\name{x}(i)$ (so $r^*_i\in N_{\gamma+(2k+2)(n+1)}$). We know that
$a^*_{\alpha^*}\subseteq a_{\gamma+(2k+2)(n+1)+n}\in N_{\gamma+(2k+2)(n+1)
+2n+1}$. Let $r_{i+1}=r^*_i\conc\langle a_{\gamma+(2k+2)(n+1)+n}\rangle$. 
Clearly $r_{i+1}\in N_{\gamma+(2k+4)(n+1)}$. By the choice of ``the
$<_\chi$--first'' conditions we are sure that, arriving to a limit stage
$i<\theta$, we have $\langle r_j:j<i\rangle\in N_{i+1}$. Now use the
assumption {\bf (b)} on $\baza_0$ to argue that the sequence $\langle r_i:
i<\theta\rangle$ has the least upper bound $r_\theta$ -- clearly this
condition decides the name $\name{x}$.]

Without loss of generality we may assume that, for each $\alpha<\gamma$
\[\forces_{\bP_\alpha}\mbox{``\/$\nbQ_\alpha$ is a partial order on
$\kappa$\/''.}\]
For a forcing notion $\bQ$ let $\check{\bQ}$ stand for the completion of $\bQ$
with respect to increasing $<\kappa$--sequences (i.e. it is like $\hat{\bQ}$
but we consider only increasing sequences of length $<\kappa$). Note that
$\bQ$ is dense in $\check{\bQ}$ and if $\|\bQ\|\leq\kappa$ then
$\|\check{\bQ}\|\leq\kappa$ ($\kappa$ is strongly inaccessible!). Now, let
$\langle \bP_\alpha',\nbQ_\alpha':\alpha<\gamma\rangle$ be the iteration of
the respective $<\kappa$--completions of the $\nbQ_\alpha$'s. Thus each
$\bP_\alpha$ is a dense subset of $\bP_\alpha'$  (see \ref{uzupelnienie}). We
may assume that each $\check{\nbQ}_\alpha$ is a $\bP_\alpha$--name for a
partial order on $\kappa+\kappa$ (for $\alpha<\gamma$). Now, for $\alpha\leq
\gamma$, let 
\[\begin{array}{ll}
\bP_\alpha''=\big\{p\in\bP_\alpha':&\mbox{there is a sequence }\langle
\bar{p}^\alpha\!:\alpha\!\in\!\dom(p)\rangle\mbox{ such that for some }\delta
\!<\!\kappa,\\
\ &\mbox{each $\bar{p}^\alpha=\langle p^\alpha_\zeta:\zeta<\delta\rangle$ is a
$\delta$--sequence of ordinals $<\kappa$ and}\\
\ &p(\alpha)\mbox{ is (the $\bP'_\alpha$--name of) the minimal (as an
ordinal)}\\
\ &\mbox{least upper bound of $\bar{p}^\alpha$ by $\leq_{\nbQ_\alpha}$}\big\}. 
  \end{array}\]

\begin{claim}
\label{cl13}
For each $\alpha\leq\gamma$, $\bP_\alpha''$ is a dense subset of
$\bP_\alpha'$.  
\end{claim}

\noindent{\em Proof of the claim:}\hspace{0.15in} Let $p\in\bP_\alpha'$. By
\ref{nextiter} we know that $\bP_\alpha$ is strongly complete for $\baza_0$. 
Let $(\bar{N},\bar{a})$ be an $\baza_0$--complementary pair such that $\bar{N} 
=\langle N_i: i<\omega\rangle$ and $p,\bar{\bQ},\bP_\alpha',\baza_0\ldots\in
N_0$. Take an increasing sequence $\langle q_i:i<\omega\rangle\subseteq
\bP_\alpha$ such that $q_i\in N_{i+1}$ is generic over $N_i$ and such that
$p\leq_{\bP_\alpha'}q_0$. Now let $q\in\bP_\alpha'$ be defined by $\dom(q)=
N_\omega\cap\alpha$ and:
\[\begin{array}{ll}
q\rest\beta\forces_{\bP_\beta}&\mbox{``\/}q(\beta)\mbox{ is the minimal (as an
ordinal) least upper bound in }\check{\nbQ}_\beta\\
\ &\mbox{\ \/of the sequence }\langle q_i(\beta): i<\omega\rangle\mbox{\/''.}
  \end{array}\]
By \ref{cl8} (actually by its proof) we know that the above definition is
correct. Now it a routine to check that $q\in\bP_\alpha''$ is as required,
finishing the proof of the claim. 
\medskip

\noindent One could ask what is the point of introducing $\bP_\alpha''$. The
main difference between $\bP_\alpha'$ and $\bP_\alpha''$ is that in the first,
$q(\beta)$ is a least upper bound of an increasing sequence of conditions from
$\nbQ_\alpha$, but we know the name for the sequence only. In $\bP_\alpha''$,
we have the representation of $q(\alpha)$ as the least upper bound of a
sequence of ordinals from $\V$! This is of use if we look at the iteration in
different universes. If we look at $\bar{\bQ}$ (defined as an iteration in
$\V$) in $\V^{\bR}$, then it does not have to be an iteration anymore: let
$\alpha<\gamma$. Forcing with $\bR$ may add new maximal antichains in
$\bP_\alpha$ thus creating new names for elements of $\nbQ_\alpha$. However

\begin{claim}
\label{cl14}
For each $\alpha\leq\gamma$, in $\V^{\bR}$, $\langle\bP_\alpha'',\nbQ_\beta:
\alpha\leq\gamma, \beta<\gamma\rangle$ is a $(<\kappa)$--support iteration.
\end{claim}

\noindent{\em Proof of the claim:}\hspace{0.15in} Easy induction on $\alpha$.

\begin{claim}
\label{cl15}
For each $\alpha<\gamma$
\[\forces_{\bP_\alpha*\bR}\mbox{``\/$\check{\nbQ}_\alpha$ is isomorphic to
$\cohen_\kappa$\/''.}\]
\end{claim}

\noindent{\em Proof of the claim:}\hspace{0.15in} Working in $\V^{\bP_\alpha}$
choose an increasing continuous sequence $\bar{N}=\langle N_i:i<\kappa\rangle$
of elementary submodels of $\Hchi$ such that for each $i<\kappa$
\[\bar{N}\rest(i+1)\in N_{i+1},\quad N_i\cap\kappa\in\kappa,\quad\mbox{ and
}\quad \|N_i\|<\kappa.\]
Now, passing to $\V^{\bP_\alpha*\bR}$, we can find an increasing continuous
sequence $\bar{j}=\langle j_\zeta:\zeta<\kappa\rangle\subseteq\kappa$ such
that 
\[(\forall\varepsilon<\kappa)(\langle N_{j_\zeta}\cap\kappa:\zeta\leq
\varepsilon\rangle\in\baza_0).\]
[Why? Forcing with $\bR$ adds an increasing continuous sequence $\bar{\beta}=
\langle\beta_\zeta:\zeta<\kappa\rangle$ such that $\bar{\beta}\rest (\zeta+1)
\in\baza_0$ for each $\zeta<\kappa$. Now let $\bar{j}$ be the increasing
enumeration of $\{j<\kappa: N_j\cap\kappa=j\ \&\ (\exists\zeta<\kappa)(j=
\beta_{\omega\cdot\zeta})\}$; remember that $\baza_0$ is closed under
subsequences.]\\
Now, for $p\in\check{\nbQ}_\alpha$ let
\[\begin{array}{ll}
j(p)=\sup\big\{j<\kappa:&j=0\mbox{ or }j\in\{j_\zeta:\zeta<\kappa\}\mbox{
and}\\
\ &p\in\bigcap\{{\cal I}\in N_j:{\cal I}\subseteq\nbQ_\alpha\mbox{ is open
dense in }V^{\bP_\alpha}\}\big\}
  \end{array}\]
and $k(p)=\min\{j<\kappa: p\in N_j\}$. Now we finish notifying that
\begin{enumerate}
\item if $\bar{p}=\langle p_\varepsilon:\varepsilon<\delta\rangle$ is
increasing in $\check{\nbQ}_\alpha$ and such that $(\forall\varepsilon<\delta)
(k(p_\varepsilon)<j(p_{\varepsilon+1}))$

then the sequence $\bar{p}$ has an upper bound in $\check{\nbQ}_\alpha$;
\item for every $j<\kappa$ the set $\{p\in\check{\nbQ}_\alpha: j(p)>j\}$ is
open dense in $\check{\nbQ}_\alpha$.
\end{enumerate}
This finishes the proof of the claim and the proposition. 

Alternatively, first prove that wlog $\gamma<\kappa^+$ and then show that
$\bP_\gamma'$ becomes $\kappa$--Cohen in $\V^{\bR}$. \QED

\begin{conclusion}
\label{gobelcon}
Assume that
\begin{description}
\item[$\V_0\models\quad$] $\kappa$ is weakly compact and GCH holds (for
simplicity),
\item[$\V_1$\ \quad]    is a generic extension of $\V_0$ making ``$\kappa$
weakly compact'' indestructible by $\cohen_\kappa$ (any member of
$\kappa$--Cohen),
\item[$\V_2=\V_1^{\bR_0}$] where $\bR_0$ adds a stationary non-reflecting
subset $S^*$ of $\kappa$ by initial segments.
\end{description}
Further, in $\V_2$, let $\baza_0=\baza_0[S^*]$, $\baza_1=\baza_1[S^*]$ be as
in the main example for the current case (see e.g. \ref{forcing}), both in
$\V_2$. Suppose that $\bar{\bQ}$ is a $(<\kappa)$--support iteration of
forcing notions on $\kappa$, say of length $\gamma^*$, complete for
$(\baza_0,\baza_1)$. Let $\V_3=\V_2^{\Lim(\bar{\bQ})}$ and let $\bR$ be the
forcing notion killing stationarity of $S^*$ in $\V_3$, but also in $\V_2$;
see \ref{weacom}.\\ 
Then
\[\V_2^{\bR}=\V_1^{\bR_0*\bR}=\V_1^{\cohen}\models \mbox{``\/$\kappa$ is weakly
compact indestructible by $\cohen$\/'',}\]
and in $\V_2^{\bR}$, the forcing notion $\Lim(\bar{\bQ})$ is adding
$\kappa$--Cohen. Consequently in $\V^{\bR}_3=(\V_2^{\bR})^{\lim(\bar{\bQ})}$,
$\lim(\bar{\bQ})$ is adding Cohens and hence $\kappa$ is weakly compact in
$\V^{\bR}_3$. 
\end{conclusion}

\begin{conclusion}
\begin{enumerate}
\item Let $\V={\bf L}$ and let $\kappa$ be a weakly compact cardinal,
$\chi^\kappa=\chi$. Then for some forcing notion $\bP$ we have, in $\V^{\bP}$:
\begin{description}
\item[(a)] there are almost free Abelian groups in $\kappa$,
\item[(b)] all almost free Abelian groups in $\kappa$ are Whitehead.
\end{description}
\item If $\V\models$ GCH then $\V^{\bP}\models$ GCH.
\item We can add:
\begin{description}
\item[(c)] the forcing does not collapse any cardinals nor changes
cofinalities, and it makes $2^\kappa=\chi$, $\chi=\|\bP\|$,
\item[(d)] if $\kappa<\theta=\cf(\theta)\leq\chi$ then we can add
$\Axkt(\baza_0,\baza_1)$.
\end{description}
If $\kappa$ is $\kappa$--Cohen indestructible weakly compact
cardinal (or every stationary set reflects)  then we may add:
\begin{description}
\item[(e)] the forcing adds no bounded subsets to $\kappa$.
\end{description}
\end{enumerate}
\end{conclusion}

\Proof 1)\quad Let $\V_0=\V$ and let $\V_1,\V_2,\bR_0$ be defined as in
\ref{gobelcon}, just $\bR_0$ adds a non-reflecting stationary subset of
$\{\delta<\kappa: \cf(\delta)=\aleph_0\}$. Working in $\V_2$ define
$\bar{\bQ}=\langle \bP_\alpha,\nbQ_\alpha:\alpha<\alpha^*\rangle$,
$\alpha^*<\chi^+$ be as in the proof of the consistency of
$\Axkt(\baza_0,\baza_1)$ in \ref{getaxiom}. The desired universe is
$\V_3=\V_2^{\bP_{\alpha^*}}$. 

Clearly, as every step of the construction is a forcing extension, we have
$\V_3=\V^{\bP}$ for some forcing notion $\bP$. The forcing notion $\bR_0\in
\V_1$ adds a non-reflecting stationary subset $S$ to $\kappa$. As
$\bP_{\alpha^*}$ preserves $(\baza_0^S,\baza_1^S)\in\ckp$ (by
\ref{pseudoiter}) the set $S$ is stationary also in $\V_3$. Since $(\forall
\delta\in S)(\cf(\delta)=\aleph_0)$ we may use $S$ to build an almost free
Abelian group in $\kappa$, so clause (a) holds. Let us prove the demand (b). 

Suppose that $G$ is an almost free Abelian group in $\kappa$ with a filtration
$\bar{G}=\langle G_i:i<\kappa\rangle$. Thus the set $\gamma(\bar{G})=\{i<
\kappa: G/G_i$ is not $\kappa$--free $\}$ is stationary. Now we consider two
cases. 

\noindent{\sc Case 1:}\qquad the set $\gamma(\bar{G})\setminus S$ is
stationary.\\
By \ref{weacom} we know that after forcing with $\bR$ (defined as there) the
cardinal $\kappa$ is still weakly compact (or hust all its stationary subsets
reflect in inaccessibles). But this forcing preserves the stationarity of
$\gamma(\bar{G})\setminus S$ (and generally any stationary subset of $\kappa$
disjoint from $S$, as $S$ does not reflect). Consequently, in $\V_3$, the set 
\[\begin{array}{ll}
\Gamma'=\{\kappa':&\kappa'\mbox{ is strongly inaccessible and}\\
\ &(\gamma(\bar{G})\setminus S)\cap\kappa'\mbox{ is a stationary subset of
}\kappa'
  \end{array}\]
is stationary in $\kappa$. Hence for some $\kappa'\in\Gamma'$ we have
$(\forall i<\kappa')(\|G_i\|<\kappa')$ and therefore the filtration $\langle
G_i:i<\kappa'\rangle$  of $G_{\kappa'}$ shows that $G_{\kappa'}$ is not free, 
contradicting ``$G$ is almost free in $\kappa$''. 

\noindent{\sc Case 2:}\qquad the set $\gamma(\bar{G})\setminus S$ is not
stationary.\\
By renaming, wlog $\gamma(\bar{G})\subseteq S$. We shall prove that $G$ is
Whitehead. So let $H$ be an Abelian group extending ${\Bbb Z}$ and let $h:
H\stackrel{\rm onto}{\longrightarrow}G$ be a homomorphism such that $\Ker(h)
={\Bbb Z}$. By \ref{forforgoe} the forcing notion $\bP=\bP_{h,H,G}$ is well
defined and it is complete for $(\baza_0^S,\baza_1^S)$ and has cardinality
$\kappa$ (and for each $\alpha<\kappa$ the set ${\cal I}_\alpha=\{p\in\bP:
G_\alpha\subseteq p\}$ is dense in $\bP$). Since $\V_3\models\Axkt$, there is
a directed set ${\cal G}\subseteq\bP$ such that ${\cal G}\cap{\cal
I}_\alpha\neq\emptyset$ for each $\alpha<\kappa$. Thus $f=\bigcup{\cal G}$ is
a lifting as required (and $G$ is Whitehead).
 
2)\quad Implicit in the proof above. \QED

\bibliographystyle{literal-plain}
\bibliography{listb,lista,listc,listf,listx}
\end{document}